\documentclass[reqno,11pt,twoside]{article}
\usepackage[hmargin=1.in, vmargin=1.in, marginparwidth=0.8in, marginparsep=0.1in, a4paper, centering]{geometry}
\usepackage[lining]{ebgaramond}
\usepackage{amsmath,amsthm}
\usepackage[ebgaramond]{newtxmath}
\usepackage[scr=boondox]{mathalpha} %needed for small calligraphic p, o, e
\usepackage[font=small, labelfont=bf, labelsep=colon, justification=centerlast, margin=0.5in]{caption}
\usepackage{array,booktabs}
\usepackage{graphicx}
\usepackage[dvipsnames,svgnames]{xcolor}%
\usepackage{mleftright} %%% Fix Spacing Around "\left(\right)" 
\mleftright % redefine \left as \mleft and \right as \mright.
\usepackage{enumerate} 
\newcommand{\eps}{\varepsilon}
\newcommand{\ext}{\mathrm{ext}}
\newcommand{\Omegaext}{\Omega^\ext}
\newcommand{\Dext}{\mathcal{D}^\ext}
\newcommand{\Dir}{\mathrm{D}}
\newcommand{\Neu}{\mathrm{N}}
\newcommand{\dr}{\mathrm{d}}
\newcommand{\ir}{\mathrm{i}}
\newcommand{\er}{\mathrm{e}}
\newcommand{\myscal}[1]{\left( #1\right)}
\newcommand{\mydotp}[1]{\left\langle #1\right\rangle}
\usepackage{tikz}
\usetikzlibrary{shapes,arrows}
\usetikzlibrary{positioning}
\tikzstyle{block} = [draw, fill=blue!20, rectangle,  minimum height=3em, minimum width=6em, align=center, rounded corners=0.1cm]
\theoremstyle{plain}
\newtheorem{theorem}{Theorem}[section]
\newtheorem{proposition}[theorem]{Proposition} 
\newtheorem{corollary}[theorem]{Corollary}
\newtheorem{lemma}[theorem]{Lemma}
\newtheorem{open}[theorem]{Open problem}
\theoremstyle{remark}
\newtheorem{remark}[theorem]{Remark}
\theoremstyle{definition}
\newtheorem{definition}[theorem]{Definition}
\newtheorem{example}[theorem]{Example}
\newcommand{\R}{\mathbb{R}}
\newcommand{\C}{\mathbb{C}}
\newcommand{\N}{\mathbb{N}}
\newcommand{\Sp}{\mathbb{S}}
\usepackage[nobottomtitles,pagestyles]{titlesec}
\titleformat{\section}
{\normalfont\large\bfseries}
{\filcenter\S\thesection.}{1ex}{\filcenter}
\titleformat{\subsection}
{\normalfont\bfseries}
{\filcenter\S\thesubsection.}{1ex}{\filcenter}
\renewcommand\footnotemark{}
\usepackage[colorlinks,allcolors=blue,pagebackref]{hyperref}
\newcommand{\mydoi}[1]{\href{https://doi.org/#1}{doi: #1}}
\newcommand{\myarXiv}[1]{\href{https://arxiv.org/abs/#1}{arXiv: #1}}
\numberwithin{equation}{section}
\usepackage{pifont}
\renewcommand*{\backrefalt}[4]{%
\ifcase #1 %
No citations%
\or
\ding{43}~p.~#2%
\else
\ding{43}~pp.~#2%
\fi}
\usepackage{etoolbox}
\newcommand{\addQEDstyle}[2]{\AtBeginEnvironment{#1}{\pushQED{\qed}\renewcommand{\qedsymbol}{#2}}\AtEndEnvironment{#1}{\popQED}}
\addQEDstyle{remark}{$\blacktriangleleft$}
\addQEDstyle{example}{$\blacktriangleleft$}
\addQEDstyle{definition}{$\blacktriangleleft$}
\newcommand*{\MR}[1]{MR~\href{https://mathscinet.ams.org/mathscinet-getitem?mr=#1}{#1}}
\renewcommand*{\MR}[1]{}
\newcommand*{\Zbl}[1]{Zbl~\href{https://zbmath.org/?q=an:#1}{#1}}
\renewcommand*{\Zbl}[1]{}
\newenvironment{ack}{\section*{Acknowledgements}}{}
\newenvironment{funding}{\section*{Funding}}{}
\title{The exterior Steklov problem for Euclidean domains%
}
\author{\qquad Lukas Bundrock\thanks{\textbf{L. B.}: Department of Mathematics, The University of Alabama, 505 Hackberry Ln, Tuscaloosa, AL 35401, USA; \href{mailto:lukas.bundrock@gmail.com}{lukas.bundrock@gmail.com};  ORCID: 0009-0001-2475-5665%
}
\and
Alexandre Girouard\thanks{\textbf{A. G.}: D\'epartement de math\'ematiques et de statistique, Pavillon Alexandre-Vachon, Universit\'e Laval, Qu\'ebec QC, G1V 0A6, Canada; \href{mailto:Alexandre.Girouard@mat.ulaval.ca}{Alexandre.Girouard@mat.ulaval.ca}; \url{https://agirouard.mat.ulaval.ca}; ORCID: 0000-0001-8823-831X%
}
\and
Denis S. Grebenkov\qquad\thanks{\textbf{D. S. G.}: Laboratoire de Physique de la Mati\`ere Condens\'ee, CNRS -- \'Ecole Polytechnique, Institut Polytechnique de Paris, 91120 Palaiseau, France; \href{mailto:denis.grebenkov@polytechnique.edu}{denis.grebenkov@polytechnique.edu}; \url{https://pmc.polytechnique.fr/pagesperso/dg/}; ORCID: 0000-0002-6273-9164%
}
\and
Michael Levitin\thanks{%
\textbf{M. L.: }Department of Mathematics and Statistics, University of Reading, 
Pepper Lane, Whiteknights, Reading RG6 6AX, UK;
\href{mailto:M.Levitin@reading.ac.uk}{M.Levitin@reading.ac.uk}; \url{https://www.michaellevitin.net}; ORCID: 0000-0003-0020-3265%
}
\and
Iosif Polterovich\thanks{%
\textbf{I. P.: }D\'e\-par\-te\-ment de math\'ematiques et de statistique, Univer\-sit\'e de Mont\-r\'eal, 
CP 6128 succ Centre-Ville, Mont\-r\'eal QC  H3C 3J7, Canada;
\href{mailto:iosif.polterovich@umontreal.ca}{\nolinkurl{iosif.polterovich@umontreal.ca}}; \url{https://www.dms.umontreal.ca/\~iossif}; ORCID: 0009-0007-0052-6589%
}
}
\date{\small arXiv:2511.09490v4, to appear in Journal of Spectral Theory}
\begin{document}
\maketitle
%\pagestyle{headings}
%%%%%%%%%%%%%%%%%%%%%%%%%%%%%%%%%%%%%%%%%%%%%%%%%%%%%%%%%%%%%%%%%
\begin{abstract} 
We investigate the Steklov eigenvalue problem in an exterior Euclidean domain. First, we present
several formulations of this problem and establish the equivalences between them. Next, we examine various properties of the exterior Steklov eigenvalues and eigenfunctions. One of our main findings is an Escobar-type lower bound for the first exterior Steklov eigenvalue on convex domains in dimensions three and higher. This bound is expressed in terms of the principal curvatures of the boundary and is sharp, with equality attained for a ball. Moreover, it implies the existence of a sequence of convex domains with fixed volume and the first exterior Steklov eigenvalues tending to infinity. This contrasts with the interior case, as well as with the two-dimensional exterior case, for which we show that an analogue of the Weinstock isoperimetric inequality holds.
\end{abstract}
%%%%%%%%%%%%%%%%%%%%%%%%%%%%%%%%%%%%%%%%%%%%%%%%%%%%%%%%%%%%%%%%%
{\small \tableofcontents}
%%%%%%%%%%%%%%%%%%%%%%%%%%%%%%%%%%%%%%%%%%%%%%%%%%%%%%%%%%%%%%%%%
\section{Introduction and main results}

\subsection{Statement of the problem}\label{sec:statement}

Let $\Omega \subset \mathbb{R}^n$, $n \ge 2$,  be a bounded open set 
with Lipschitz boundary 
$\partial\Omega$. The classical Steklov eigenvalue problem in $\Omega$ is given by
\begin{equation}\label{eq:PDEint}
\begin{cases}
\Delta u  = 0\qquad&\text{ in } \Omega, \\
\partial_{\overline{\nu}}u  = \sigma u\qquad&\text{ on } \partial\Omega,
\end{cases}
\end{equation}	
where $\partial_{\overline{\nu}} u = \mydotp{ \nabla u, \overline{\nu} }$ is the normal derivative of $u$ in the direction of the  unit normal $\overline{\nu}$ pointing towards the exterior of $\Omega$, and $\sigma$ is a spectral parameter. 
Since the trace operator $H^{1}(\Omega) \to L^2(\partial\Omega)$ is compact, the Steklov spectrum is discrete. It consists of a sequence of eigenvalues
\[
0 = \sigma_1(\Omega) \leq \sigma_2(\Omega) \leq \sigma_3(\Omega)\leq \ldots\nearrow+\infty,
\]
accumulating at infinity. The Steklov eigenvalues can also be understood as the eigenvalues of the Dirichlet-to-Neumann operator $\mathcal{D}: H^\frac{1}{2}(\partial\Omega) \to H^{-\frac{1}{2}}(\partial\Omega)$, defined by $\mathcal{D}: f\mapsto  \partial_{\overline{\nu}} (\mathcal{H}f)$, where $\mathcal{H}f\in H^1(\Omega)$ is the harmonic extension of a function $f\in H^\frac{1}{2}(\partial\Omega)$ from the boundary to the interior.
The eigenfunctions of the Dirichlet-to-Neumann operator are the boundary traces of the Steklov eigenfunctions,
and they can be chosen to form an orthonormal basis of $L^2(\partial\Omega)$.
For a comprehensive survey and detailed analysis of the Steklov spectrum on bounded domains (as well as on compact Riemannian manifolds with boundary) we refer to \cite{girouard2017spectral, colbois2024some} as well as to \cite[Chapter 7]{levitin2023topics}.  

The Steklov problem and the Dirichlet-to-Neumann map have many physical applications, notably, to the modelling of diffusion processes with a particular behaviour of the particles when they reach the boundary surface, see \cite{grebenkov2020paradigm, grebenkov2024steklov}. 
Problem \eqref{eq:PDEint} arises if the diffusion occurs in the interior of $\Omega$. However, in many situations it is of interest to consider the diffusion in the {\em exterior} domain
\[
\Omegaext := \mathbb{R}^n\setminus \overline{\Omega}.
\]
It is therefore natural to look for an analogue of the Steklov problem \eqref{eq:PDEint} on an {\em unbounded} domain $\Omegaext$. 
Problems of this type also arise in scattering theory,  when modelling  the interaction of waves with bounded obstacles, see, e.g., \cite{CDMM}. Another motivation for the exterior Steklov problem comes from the study of the exterior Robin eigenvalue problem  \cite{krejcirik2016optimisation, krejvcivrik2020optimisation, krejcirik2023optimisation, bundrock2024optimizing}, which is in a sense dual to the Steklov problem, see \S\ref{subsec:Robin}.

Let us fix some basic notations. A point $x\in \R^n$ will be represented in Cartesian coordinates as $x=(x_1,\dots,x_n)$ and in polar coordinates as $x=(r, \theta)$, where $r=|x|$ and $\theta\in\Sp^{n-1}$. We will denote by 
\[
B_R=\left\{x\in\R^n: |x|<R\right\}
\]
the ball of radius $R$ centred at the origin, and set
\[
\Omegaext_R:=\Omegaext \cap B_R,
\]
for each
\[
R>R_0(\Omega) := \inf\,\left\{R > 0 : \overline{\Omega} \subset B_R \right\}.
\]
We also denote
\[
H^1_\mathrm{loc}\left(\Omegaext\right):=\left\{u:\Omegaext\to\mathbb{R}: u\in H^1\left(\Omegaext_R\right)\text{ for all }R>R_0(\Omega)\right\}.
\]

Simply  replacing  $\Omega$  by $\Omegaext$ in \eqref{eq:PDEint} does not lead to a well-defined exterior Dirichlet-to-Neumann map. 
Indeed,  a harmonic extension of a function to  an unbounded domain is not unique unless certain conditions are imposed at infinity. This can be immediately seen by considering the function $f \equiv 0$ on the boundary of the unit disk $B_1\subset\R^2$: it can be extended harmonically to the exterior of the disk as $u \equiv 0$  or as $u=\log r$. One way to fix the choice of the harmonic extension is to assume that 
\begin{equation}\label{eq:inftyr}
u(x)=O\left(|x|^{2-n}\right)\qquad\text{ as } |x| \to \infty, \quad x \in \Omegaext\subset\mathbb{R}^n.
\end{equation}

\begin{theorem}[see {\cite[Theorem 3.1]{KNP25}} and {\cite[Theorem 8.10]{mclean2000strongly}}]\label{thm:harmext}
Let $\Omega\subset\mathbb{R}^n$ be a bounded open set with Lipschitz boundary and  with connected $\Omegaext$. For any $f\in H^{1/2}(\partial\Omega)$, there exists a unique harmonic extension $u=\mathcal{H}^\ext f\in H^1_\mathrm{loc}\left(\Omegaext\right)$ satisfying  \eqref{eq:inftyr}.
\end{theorem}

\begin{figure}[t]
\centering
\includegraphics{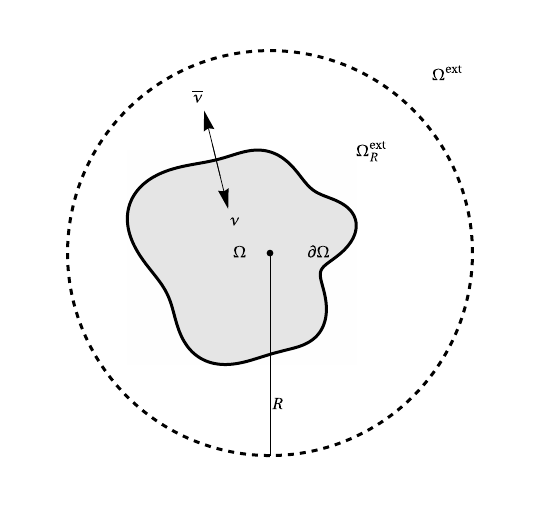}
\caption{The geometry of an exterior problem.}\label{fig:domains}
\end{figure}

\begin{definition}\label{defn:Steklovext}
The {\em exterior Steklov problem} on $\Omegaext \subset \mathbb{R}^n$ is defined as follows:  find $\sigma\in\mathbb{R}$ for which there exists a nonzero $u:\Omegaext \to\mathbb{R}$ such that $u\in H^1_\mathrm{loc}\left(\Omegaext\right)$, and
\begin{equation*}\label{eq:PDEext}\tag{$\mathrm{ES}$}
\begin{cases}
\Delta u  = 0\qquad&\text{ in } \Omegaext, \\
\partial_\nu u  = \sigma u\qquad&\text{ on } \partial\Omega,\\
\text{$u$ satisfies \eqref{eq:inftyr}}.
\end{cases}
\end{equation*}
In \eqref{eq:PDEext} and further on, $\partial_\nu u$ is the normal derivative of $u$ in the direction of the  unit normal $\nu=-\overline{\nu}$ pointing towards the exterior of $\Omegaext$ (and therefore towards the interior of $\Omega$), see Figure \ref{fig:domains}. 
\end{definition}

We denote by $\Dext: H^\frac{1}{2}(\partial\Omega) \to H^{-\frac{1}{2}}(\partial\Omega)$, $\Dext f = \partial_\nu\left(\mathcal{H}^\ext f\right)$,  the corresponding {\em exterior Dirichlet-to-Neumann operator}. The following basic result holds. 

\begin{theorem}\label{thm:basic}
Let $\Omega \subset \mathbb{R}^n$ be a bounded open set with Lipschitz boundary and  with connected $\Omegaext$.
The spectrum of the exterior Steklov problem \eqref{eq:PDEext} in $\Omegaext$ is discrete, and consists of a sequence of eigenvalues
\[
0  \leq \sigma_1\left(\Omegaext\right)  \leq  \sigma_2\left(\Omegaext\right) \leq \sigma_3\left(\Omegaext\right)\leq \ldots\nearrow+\infty,
\]
accumulating at infinity.
The  boundary traces $f_k=u_k|_{\partial\Omega}$ of the corresponding Steklov eigenfunctions can be chosen to form an orthonormal basis of $L^2(\partial\Omega)$. 
\end{theorem}
The proof of this theorem can be essentially deduced from the results of \cite{auchmuty2014representations, arendt2015dirichlet}, see \S\ref{subsec:basic}. 

\begin{remark}\label{condition}
Without loss of generality, it will always be assumed in what follows that $\Omegaext$ is a {\em connected} unbounded domain. Otherwise, $\Omegaext$  would be a disjoint union of finitely many bounded domains and a single connected unbounded domain. Hence, it  suffices to consider the exterior Steklov problem  on connected  $\Omegaext$ only.
\end{remark}

\begin{remark}\label{rem:firsteig}
Note that  for any {\em planar} open set $\Omega$ with connected $\Omegaext$, the function  $u_1=\operatorname{const}$ solves \eqref{eq:PDEext} with $\sigma=\sigma_1\left(\Omegaext\right)=0$. In contrast, $\sigma_1\left(\Omegaext\right)$ is strictly positive in dimensions $n\ge 3$, as will follow from Theorems \ref{thm:SpectralADAN} and  \ref{thm:spectral}. 
\end{remark} 

\subsection{Different approaches to the exterior Steklov problem and their equivalence}

One of the difficulties arising in the study of the exterior Steklov problem is that  the solutions of~\eqref{eq:PDEext} may no longer be square-integrable, like the constant function in dimension two, or the function $1/r$  in the exterior of a ball in $\mathbb{R}^3$, see \S\ref{examp:balln2} and \S\ref{examp:balln3} for further details. In particular, the space 
$H^1\left(\Omegaext\right)$ is not suitable for the variational characterisation of exterior eigenvalues.
 
One of the goals of the present paper is to unify different approaches to the exterior Steklov problem by showing that they are all equivalent in an appropriate sense to Definition \ref{defn:Steklovext}. 
We list  these approaches in an informal way below and provide references to precise formulations: 

\begin{itemize}
\item\textbf{Approach I: Finite energy spaces in dimensions $n\ge 3$.} Solve the exterior Steklov problem using variational characterisation \eqref{WeakFormulationExterior} in the space of functions of finite energy, i.e.\  functions with their gradients in $L^2\left(\Omegaext\right)$, and satisfying a certain decay condition at infinity, see \S\ref{subsec:fen}. This approach was developed in \cite{auchmuty2014representations, arendt2015dirichlet}.
\item\textbf{Approach II: Conformal mapping in dimension two. }
Use a conformal mapping  to transform the exterior Steklov problem to a weighted interior Steklov problem \eqref{eq:problemtransform}, see \S\ref{subsec:conformal}.
\item\textbf{Approach III: Truncated domains. } Given $R>R_0(\Omega)$, consider a mixed Steklov--Dirichlet eigenvalue problem \eqref{eq:PDEArendtElstDirichlet} on $\Omegaext_R$ with Steklov condition on $\partial\Omega$ and Dirichlet condition on $\partial B_R$, and take the limit of eigenvalues and eigenfunctions as $R \to \infty$, see \S\ref{subsec:truncatedomains}. This approach was explored  in \cite{arendt2015dirichlet} in dimension $n \ge 3$.
\item\textbf{Approach IV: Helmholtz equation. } Consider the exterior Steklov problem \eqref{eq:PDEGrebChai} for the Helmholtz equation $(-\Delta+\Lambda^2) u = 0$ with a given  $\Lambda>0$.  The  decaying solutions of  \eqref{eq:PDEGrebChai} belong to $H^1\left(\Omegaext\right)$, and one can use the variational characterisation of eigenvalues in the usual way. One can then take the limit of eigenvalues and eigenfunctions as $\Lambda \searrow 0$, see \S\ref{subsec:Greben}. This approach was developed in \cite{grebenkov2024steklov, christiansen2023low}. 
\item\textbf{Approach V: Layer potentials. } Finally, one can define the exterior Steklov problem using the layer potential approach, see \eqref{eq:Steklov_integral} for dimensions $n\ge 3$ and \eqref{eq:ext_2D_2} for $n=2$, in a similar way to how it is usually done for the interior Steklov problem. We refer to \S\ref{sec:layer} for details. An analogous approach was used in \cite{KNP25}.
\end{itemize}

\begin{figure}[htb]
\centering
\begin{tikzpicture}[>=latex']
\node [block] (CM) at (0,0) {Conformal\\mapping};
\node [block] (TD) at  (-4.5,0) {Truncated\\Dirichlet};
\node [block] (TN) at  (0,-2) {Truncated\\Neumann};
\node [block] (HE) at  (4.5,0)  {Helmholtz\\equation};
\node [block] (LP) at (4.5,-2) {Layer\\potentials};
\node [block] (ES) at (0,2) {Exterior\\Steklov};
\node [block] (FV) at (-4.5,2) {Vanishing\\flow};
\draw [double distance=1pt, <->] (CM) -- node[midway, above]{Theorem  \ref{theo:sigmaD*}} (TD);
\draw [double distance=1pt, <->] (CM) -- node[midway, above]{Theorem \ref{theo:lambda*}} (HE);
\draw [double distance=1pt, <->] (CM) -- node[yshift=3pt][midway, right]{Theorem \ref{thm:CMLP2}} (LP);
\draw [double distance=1pt, <->] (CM) -- node[midway, right]{Proof of Theorem \ref{thm:basic}} (ES);
\draw [double distance=1pt, <->] (CM) -- node[midway, left]{Theorem \ref{theo:convEVNeum}} (TN);
\draw [double distance=1pt, <->] (FV) -- node[midway, above]{Remark \ref{prop:2D-equivalence}} (ES);
\end{tikzpicture}\\\ 
\caption{Relations between approaches for $n=2$.}\label{fig:scheme2d}
\end{figure}
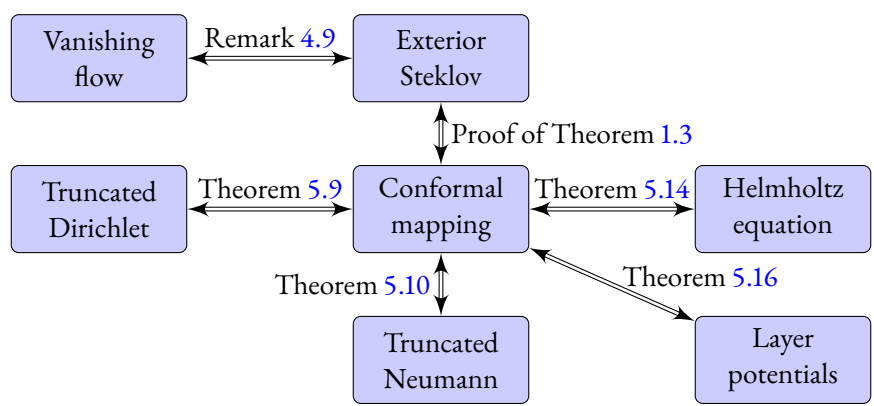

The equivalence of these formulations is understood as follows. We will show that for Approaches I, II, and V,  the corresponding  eigenvalues and eigenfunctions coincide with those of the exterior Steklov problem \eqref{eq:PDEext}, see Definition \ref{defn:Steklovext}. For Approaches III and IV, we will prove that  the $k$th eigenvalues  of the problems \eqref{eq:PDEArendtElstDirichlet} as $R \to \infty$ and  \eqref{eq:PDEGrebChai} as $\Lambda\searrow 0$ converge to $\sigma_k\left(\Omegaext\right)$ for all $k\ge 1$, and the corresponding eigenfunctions converge in an appropriate sense to the eigenfunctions  $u_k$ of the exterior Steklov problem. For a schematic relation between different approaches, see Figures \ref{fig:scheme3d} and \ref{fig:scheme2d}. For the precise statements, we refer to \S\ref{sec:layer} and \S\ref{sec:equivalence}.

\begin{remark} \label{rem:other}
There exist other possible approaches to the exterior Steklov problem. One approach, useful for numerical analysis, is to reduce \eqref{eq:PDEext} to a mixed problem \eqref{eq:SPsi} in a bounded domain $\Omegaext_\rho$, with the Steklov condition on $\partial\Omega$ and a pseudodifferential matching condition on $\partial B_\rho$, see \S\ref{sec:numerics}. 

Also, one can modify the truncated domains approach and
impose Neumann conditions   on $\partial B_R$ instead of the Dirichlet conditions. 
In two dimensions this formulation is equivalent to the others, see \S\ref{section:twod}.
At the same time, as was shown in \cite{arendt2015dirichlet}, in dimensions $n\ge 3$ the Neumann truncation does  not provide an equivalent formulation to the ones described above. In fact, this
can be seen from the fact that in dimensions $n\ge 3$ the first exterior Steklov eigenvalue is positive, while the Neumann truncation yields a zero eigenvalue. In these dimensions, the Neumann truncation corresponds instead to the exterior Steklov problem with a decay condition on the gradient (vanishing flow), as described in \S\ref{subsec:vanishflow}. 

An approach closely related to  \eqref{eq:PDEGrebChai}  has been introduced in  \cite{bernard2025magneticdirichletneumannoperator}, see also \cite{HKN25}. 
Instead of introducing a Helmholtz parameter $\Lambda$, one can  consider a magnetic Laplacian  and examine the spectral behaviour as the magnetic potential vanishes. For the exterior of the unit disc, this formulation is equivalent to \eqref{eq:PDEext}, see \cite[Proposition 1.2]{bernard2025magneticdirichletneumannoperator} and \S \ref{examp:balln2}. We expect that adapting the methods of \S\ref{sec:AHvGC} and \S\ref{sec:equivhelmn2} one can show that the same equivalence holds for arbitrary Euclidean domains. Note that compared to the Helmholtz regularisation, there is one additional difficulty: the magnetic quadratic form does not contain an explicit $L^2$-term, and bounded magnetic energy does not directly give the corresponding $H^1$-bound. However, 
this can be resolved using the diamagnetic inequality, which controls the Dirichlet energy of the modulus of a function by its magnetic energy.
\end{remark}

\begin{figure}[htb]
\centering
\begin{tikzpicture}[>=latex']
\node [block] (FE) at (0,0) {Finite\\energy};
\node [block] (TD) at  (-4.5,0) {Truncated\\Dirichlet};
\node [block] (TN) at  (0,-2) {Truncated\\Neumann};
\node [block] (HE) at  (4.5,0)  {Helmholtz\\equation};
\node [block] (LP) at (4.5,-2) {Layer\\potentials};
\node [block] (ES) at (0,2) {Exterior\\Steklov};
\node [block] (FV) at (-4.5,-2) {Vanishing\\flow};
\draw [double distance=1pt, <->] (FE) -- node[midway, above]{Theorem  \ref{thm:truncation.equiv.finiteenergy}} (TD);
\draw [double distance=1pt, <->] (FE) -- node[midway, above]{Theorem  \ref{theo:lambdamuallg}} (HE);
\draw [double distance=1pt, <->] (FE) -- node[yshift=3pt][midway, right]{Corollary \ref{cor:FELP}} (LP);
\draw [double distance=1pt, <->] (FE) -- node[midway, right]{Proof of Theorem \ref{thm:basic}} (ES);
\draw [double distance=1pt, <->] (FV) -- node[midway, above]{Theorem \ref{thm:FVTN}} (TN);
\end{tikzpicture}\\\ 
\caption{Relations between approaches for $n\ge 3$.}\label{fig:scheme3d}
\end{figure}

\begin{remark} Aside from the pure Steklov conditions on $\partial\Omega$, one can also consider exterior problems with mixed Steklov--Dirichlet--Neumann boundary conditions. This  setting naturally arises in some applications, see \cite{henrici1970sloshing,grebenkov2025steklov}.  The formulations of the exterior problem discussed above,  as well as of Theorem \ref{thm:basic} (with $L^2(\partial\Omega)$ replaced by  $L^2(\partial_S\Omega)$, where $\partial_S\Omega$ is the part of the boundary on which the Steklov condition is imposed), extend directly to this case. 
\end{remark}

\subsection{Spectral properties of the exterior Steklov problem}\label{subsec:spectral}
The approaches outlined in the previous section allow us to prove various  properties of  eigenvalues and eigenfunctions of the exterior Steklov problem.  

\subsubsection{Courant's theorem and multiplicity bounds}
Recall that a \emph{nodal domain} of a continuous function is a connected component of the complement of its \emph{nodal set} (that is, the set of zero values of the function). One of the fundamental properties of the Laplace and the (interior) Steklov eigenfunctions is Courant's nodal domain theorem.  It is easy to show that it holds for the exterior Steklov eigenfunctions as well. 

\begin{theorem}\label{thm:Courant}
Let $\Omega \subset \mathbb{R}^n$, $n \geq 2$, be a bounded open set with Lipschitz boundary and  with connected $\Omegaext$. An eigenfunction $u_k$ associated to $\sigma_k\left(\Omegaext\right)$ has at most $k$ nodal domains. 
\end{theorem}

The result in two dimensions follows immediately from the conformal approach and Courant's theorem for the interior problem, and in dimensions $n\ge 3$ one can use the finite energy approach and apply  the same  proof  as in the interior case, see  \cite{kuttler1969inequality}. 
The  proof of Theorem \ref{thm:Courant} is given in \S\ref{subsec:Courant}.

Courant's theorem implies that the first exterior Steklov eigenfunction does not change sign in $\Omegaext$. In dimension two this is clear, since the first eigenfunction is constant, however in higher dimensions it is a priori not obvious.
Moreover, we have the following

\begin{corollary}\label{corol:simple}
The first exterior Steklov eigenvalue $\sigma_1\left(\Omegaext\right)$ is simple. 
\end{corollary}

Once again, this is immediate in dimension two, and the proof 
in dimension $n\ge 3$ is presented in \S\ref{subsec:Courant}.

\begin{remark} 
Similarly to the interior problem, there is no direct analogue of Courant's theorem for the nodal domains of the eigenfunctions of the exterior Dirichlet-to-Neumann map (i.e.\ for the boundary nodal domains of Steklov eigenfunctions).  In two dimensions, Theorem \ref{thm:Courant} and elementary topological considerations yield a bound on the number of boundary sign changes, cf.\  \cite[Section 6]{girouard2017spectral}. Also note that using  the results on the exterior Dirichlet-to-Neumann semigroup obtained in \cite[Section 5]{arendt2015dirichlet} one can deduce that the first exterior Steklov eigenfunction on any Lipschitz domain can be chosen to be strictly positive on the boundary, cf.\ \cite{arendt2020strict}.
\end{remark}

\subsubsection{Escobar-type lower bounds for the first eigenvalue}

For the interior problem, it was shown by Payne  \cite{payne1970some} that the first nonzero Steklov eigenvalue on a planar  convex  bounded domain is estimated from below by the minimal curvature and from above by the maximal curvature of the boundary. Note that both estimates are sharp and attained on a disk. Escobar found a non-sharp extension  of Payne's lower bound to two-dimensional  compact manifolds with boundary \cite{escobar1997geometry}. He also conjectured a sharp bound that in the case of Euclidean domains of arbitrary dimension was proved in \cite{xia2024escobar}. It states that if all the principal  curvatures of the boundary  $\partial\Omega$ are bounded below by a constant $c$ then the first nonzero Steklov eigenvalue $\sigma_2(\Omega) \ge c>0$ with equality if and only if $\Omega$ is a ball. 

For the exterior Steklov problem in dimension $n\ge 3$ we use the finite energy space approach to  prove a stronger Escobar-type bound. We show that  the first exterior Steklov eigenvalue is bounded from below by the minimum of the logarithmic (and, hence, the geometric) mean of the curvatures, which itself is bounded from below by the  minimum of the  principal curvatures.

Let $k\in \mathbb{N}$, and $\alpha_1,\dots, \alpha_k$ be distinct positive real numbers. Following \cite{pittenger1985logarithmic}, we define their {\em logarithmic mean}   $L(\alpha_1,\ldots, \alpha_k)$ as
\begin{equation}\label{eq:LM}
L(\alpha_1,\ldots, \alpha_k) = \frac{1}{(k-1)\sum_{j=1}^{k} \frac{	\alpha_j^{k-2}}{\prod\limits_{i=1, \, i \neq j}^{k}  (\alpha_j-\alpha_i)} \log \alpha_j}. 
\end{equation}
If the $\alpha_i$ are not distinct, or  if an $\alpha_i$ vanishes, $L$ is defined as the limit of the right-hand side of \eqref{eq:LM}. Note that for $k>2$ there exist other definitions of the logarithmic mean in the literature. For $k=2$, however, they all agree and yield
\[
L(\alpha_1,\alpha_2) = \begin{cases}
\frac{\alpha_1-\alpha_2}{\log \alpha_1 - \log \alpha_2} &\text{ for } \alpha_1 \neq \alpha_2,\\
\alpha_1 &\text{ for } \alpha_1 = \alpha_2.
\end{cases}
\]

For $j=1, \ldots, n-1$ and $s \in \partial \Omega$, let $\kappa_j(s)$ denote the $j$th principal curvature of $\partial\Omega$, where defined. We adopt the convention that the principal curvatures are nonnegative if 
$\Omega$ is convex. 

\begin{theorem}\label{prop:low} 
Suppose that $\Omega \subset \mathbb{R}^n$, $n \geq 3$, is a bounded, convex domain  with $\partial\Omega \in C^{1,1}$. Then, 
\begin{equation}\label{eq:logmean}
\sigma_1\left(\Omegaext\right)   \geq (n-2) \inf_{s \in \partial\Omega}	 L\left(\kappa_1(s), \ldots, \kappa_{n-1}(s)\right)=:\beta(\partial\Omega).
\end{equation}
\end{theorem}

For $\partial \Omega \in C^{1,1}$, the principal curvatures are defined almost everywhere. Accordingly, the infimum in \eqref{eq:logmean} and in \eqref{eq:Ksgeommean} is taken over points where the curvatures are defined. 

Additionally, Theorem~\ref{prop:low} immediately yields a lower bound expressed in terms of the geometric mean of the curvatures. We expect that equality in~\eqref{eq:logmean} holds if and only if 
$\Omega$ is a ball, although we were only able to establish this in the setting of Corollary~\ref{coro:lowbound}.

\begin{corollary}\label{coro:lowbound}
Suppose that $\Omega \subset \mathbb{R}^n$, $n \geq 3$, is a bounded, convex domain with $\partial\Omega \in C^{1,1}$. Then,
\begin{equation}\label{eq:Ksgeommean}
\sigma_1\left(\Omegaext\right)   \geq (n-2) \inf_{s \in \partial\Omega}	 \sqrt[n-1]{\prod_{j=1}^{n-1}  \kappa_j(s)},
\end{equation}
with equality if and only if $\Omega$ is a ball.
\end{corollary}

The proofs of Theorem \ref{prop:low} and Corollary \ref{coro:lowbound} are given in \S\ref{sec:lowerbounds}.

\begin{remark}
A different lower bound for the first exterior Steklov eigenvalue was obtained in  \cite[Theorem 1]{xiong2023sharp}:  for any bounded open set $\Omega \subset \mathbb{R}^n$, $n \geq 3$, with $C^1$ boundary and with $0 \in \Omega$, the first exterior Steklov eigenvalue satisfies
\begin{equation}\label{eq:Xiongbound}
\sigma_1\left(\Omegaext\right) \geq (n-2) \min_{x \in \partial \Omega} \frac{\mydotp{x, -\nu }}{|x|^2}=:\beta_\mathrm{X}\left(\partial\Omega\right)
\end{equation}
(the bound is non-trivial only for star-shaped domains $\Omega$).
In Example \ref{examp:ellipsoid}, we compare \eqref{eq:Xiongbound} with our bound from Theorem \ref{prop:low} for various spheroids, showing that neither inequality implies the other one.
\end{remark}

\begin{remark}\label{rem:compesc}
Note that one cannot replace the minimum of principal curvatures by the minimum of their geometric mean in the original Escobar's conjecture for the first nonzero interior Steklov eigenvalue. Indeed, consider a prolate spheroid $\mathscr{p}_a:=\left\{\frac{x_1^2}{a^2} +\frac{x_2^2}{a^2} +x_3^2=1\right\}$, with $0<a<1$,  which divides $\mathbb{R}^3$ into a bounded domain $\mathcal{P}_a$ and an unbounded domain $\mathcal{P}^\ext_a$. One can check that the minimum of the geometric mean of the principal curvatures of $\mathscr{p}_a$  (i.e.\ the square root of the Gaussian curvature) is equal to one for any $a>0$. At the same time, taking the trial function $x_3$ in the Rayleigh quotient for $\sigma_2\left(\mathcal{P}_a\right)$ (note that $x_3$ is orthogonal to constants on $\mathscr{p}_a$) one can show that $\sigma_2\left(\mathcal{P}_a\right) \to 0$ as $a\searrow 0$. Therefore, for the interior problem, the geometric mean of the curvatures cannot be used to bound $\sigma_2$ from below.
\end{remark}

\begin{remark}
The upper bound in \cite{payne1970some} follows directly from the Gauss--Bonnet formula for curves together with Weinstock's inequality. Since an analogue of Weinstock's inequality also holds for the exterior problem in two dimensions (see Theorem~\ref{theo:weinstock} below), the same upper bound extends to the first nonzero exterior Steklov eigenvalue. In contrast, the proof of the lower bound established in \cite{payne1970some}, does not carry over immediately to the exterior case. It would be interesting to check whether such a lower bound remains valid for the exterior problem.
\end{remark}

\subsubsection{Isoperimetric upper bounds for the first eigenvalue}

Let us start with an analogue of the Weinstock inequality for  the first nonzero exterior Steklov eigenvalue  
$\sigma_2$ in dimension two.

\begin{theorem}\label{theo:weinstock}
Suppose $\Omega \subset \mathbb{R}^2$ is a bounded simply connected Lipschitz domain. Then
\begin{equation}\label{eq:weinstock}
 \sigma_2\left(\Omegaext\right) |\partial\Omega| \leq 2\pi.
\end{equation}
with  equality if and only if \,$\Omega$ is a disk.
\end{theorem}

Moreover, combining \eqref{eq:weinstock} with the isoperimetric inequality we obtain
\[
 \sigma_2\left(\Omegaext\right) |\Omega|^{1/2} \leq \sqrt{\pi},
\]
with equality if and only if $\Omega$ is a disk.

Inequality \eqref{eq:weinstock} follows from the conformal mapping approach and the usual Weinstock inequality, see the proof in \S\ref{sec:shapeopt}. 
The proof of the equality case contains additional steps compared to the interior Steklov problem, because one needs to keep track of the image of infinity under the conformal mapping. 

\begin{remark}
In fact, the conformal approach essentially allows one to extend any result that is valid for the weighted interior Steklov problem  in two dimensions to the exterior Steklov problem. In particular, assuming that $\Omega\subset\mathbb{R}^2$ is simply connected we obtain the analogues of   the Hersch--Payne--Schiffer inequalities \cite{girouard2010hersch},
\[
\sigma_{k+1}\left(\Omegaext\right) |\partial\Omega|  \le 2\pi k,\qquad k\in\mathbb{N},
\]
as well as the  eigenvalue asymptotics \eqref{eq:rozenas} below.
\end{remark}

Consider now the case $n\ge 3$. In the interior case, Brock's inequality \cite{brock2001isoperimetric}
implies that the first nonzero Steklov eigenvalue $\sigma_2$ attains its maximum on a ball among all domains of given volume. Moreover, there is a higher-dimensional version of Weinstock's inequality \cite{bucur2021weinstock} stating that among all {\em convex} domains with given surface area, $\sigma_2$ is maximised by a ball. It turns out that neither of these results  holds in the exterior case.

\begin{theorem}\label{theo:noisopineq}
Let $n\ge 3$. There exists a sequence of convex smooth bounded domains $\Omega_m \subset \mathbb{R}^n$  of fixed volume such that $\sigma_1\left(\Omegaext_m\right) \to +\infty$ as $m \to \infty$.
\end{theorem}

The proof of this theorem, see \S\ref{sec:shapeopt}, is obtained by applying Theorem~\ref{prop:low} to prolate spheroids. Moreover, by the classical isoperimetric inequality, the sequence in Theorem~\ref{theo:noisopineq} can also be normalised by surface area instead of by volume.

\subsubsection{Exterior Robin problem}\label{subsec:Robin}

As a consequence of Theorem \ref{theo:noisopineq}, we obtain the following result for the exterior Robin problem. For $\alpha \in \mathbb{R}$ and a bounded Lipschitz domain $\Omega \subset \mathbb{R}^n$ with connected $\Omegaext$, the lowest point of the spectrum of the Robin Laplacian $-\Delta^{\mathrm{Rob},\alpha}$ in $\Omegaext$ with the parameter $\alpha$ in the Robin condition  is given by
\[
\lambda_1^\alpha\left(\Omegaext\right) := \min\operatorname{Spec}\left(-\Delta^{\mathrm{Rob},\alpha}\right) = \inf_{u \in H^1\left(\Omegaext\right)} \frac{ \int_{\Omegaext} | \nabla u|^2 \, \dr x + \alpha \int_{\partial \Omega} |u|^2 \, \dr S}{\int_{\Omegaext} | u|^2 \, \dr x},
\]
see \cite{krejcirik2016optimisation,krejvcivrik2020optimisation}. For $n \geq 3$, it was shown in \cite[Theorem 1]{bundrock2024optimizing} that
\begin{equation}\label{eq:Robin}
\lambda_1^\alpha\left(\Omegaext\right) < 0 \qquad\text{if and only if}\qquad \alpha < - \sigma_1\left(\Omegaext\right),
\end{equation}
in which case it is an eigenvalue. Otherwise $\lambda_1^\alpha\left(\Omegaext\right) = 0 =  \min\operatorname{Spec}_\mathrm{ess}\left(-\Delta^{\mathrm{Rob},\alpha}\right)$ is the bottom of the essential spectrum. Note that in \cite{bundrock2024optimizing}, $\Omega$ is assumed to be a domain (and hence connected), but the proof of \eqref{eq:Robin} carries over verbatim to any open bounded set $\Omega$ with Lipschitz boundary and  with  connected $\Omegaext$.

Thus, for given $\alpha \in \mathbb{R}$ and $R>0$, Theorem~\ref{theo:noisopineq} ensures the existence of a smooth convex domain $\Omega \subset \mathbb{R}^n$ such that $|\Omega| = |B_R|$ and $\sigma_1(\Omega^{\ext}) \ge -\alpha$. Then by \eqref{eq:Robin}, $\lambda_1^\alpha\left(\Omegaext\right) = 0$. If we now only consider  the values $\alpha\in\left(-\infty, -\sigma_1\left(B_R^\ext\right)\right)$, then we have $\lambda_1^\alpha\left(B_R^\ext\right) < 0=\lambda_1^\alpha\left(\Omegaext\right)$. The preceding argument gives the following

\begin{corollary}\label{cor:Robin}
Let $n \ge 3$.  For any $R>0$ and $\alpha < -\sigma_1\left(B_R^\ext\right) = -\frac{n-2}{R}$, there exists a smooth convex domain $\Omega \subset \mathbb{R}^n$ with $|\Omega| = |B_R|$ satisfying $\lambda_1^\alpha\left(B_R^\ext\right) < \lambda_1^\alpha\left(\Omegaext\right)$. 
\end{corollary}

It was already observed in \cite{krejcirik2016optimisation} that the ball does not maximise $\lambda_1^\alpha$ as $\alpha \to -\infty$. Corollary \ref{cor:Robin} extends this observation to all relevant values of $\alpha$. The duality between the Robin and Steklov eigenvalues in \eqref{eq:Robin} extends to higher eigenvalues and also holds in dimension $n=2$ \cite[Corollary~1]{bundrock2024optimizing}. Specifically, for $n \geq 2$, for any bounded open Lipschitz set $\Omega \subset \mathbb{R}^n$ with connected $\Omegaext$, and for $k \in \mathbb{N}$,
\begin{equation}\label{eq:Robinhigher}
\lambda_k^\alpha\left(\Omegaext\right) < 0 \quad \text{ if and only if } \quad \alpha < - \sigma_k\left(\Omegaext\right).
\end{equation}
 
In dimension two, $\sigma_1\left(\Omegaext\right) = 0$ for any $\Omega$ and therefore $\lambda_1^\alpha\left(\Omegaext\right) < 0$ whenever $\alpha < 0$. This was also shown directly for the first Robin eigenvalue in \cite[Proposition~2]{krejcirik2016optimisation}. Moreover, in dimension $n=2$, for any $\alpha < 0$ the disk maximises $\lambda_1^\alpha\left(\Omegaext\right)$ among smooth bounded simply connected domains $\Omega$ with fixed perimeter or area \cite[Corollary~5]{krejvcivrik2020optimisation}, see also \cite[Theorem~4]{krejvcivrik2020optimisation}  for results regarding disconnected $\Omega$.
 
In \cite{krejcirik2023optimisation}, Krej\v{c}i\v{r}\'{i}k and Lotoreichik study the second Robin eigenvalue $\lambda_2^\alpha\left(\Omegaext\right) $ in dimension $n=2$. They observe the existence of a threshold value $\alpha_* \left(\Omegaext \right) < 0$ such that $\lambda_2^\alpha\left(\Omegaext\right) < 0$ if and only if $  \alpha < \alpha_* \left(\Omegaext \right)$. In view of \eqref{eq:Robinhigher}, this threshold satisfies $\alpha_* \left(\Omegaext \right) = - \sigma_2\left(\Omegaext\right)$.  They conjecture  \cite[Conjecture~1.2]{krejcirik2023optimisation} that $\lambda_2^\alpha\left(\Omegaext\right) < \lambda_2^\alpha(B_R^\mathrm{ext})$ for all bounded simply connected open sets $\Omega$ with the same perimeter or area as $B_R$. However, this conjecture remains open and instead they prove that the disk maximises $\lambda_2^\alpha\left(\Omegaext\right)$ among convex domains with prescribed \emph{elastic energy}, given by 
$\frac{1}{2}\int_{\partial \Omega} \kappa(s)^2 \, \mathrm{d}S$, 
and they derive an inequality for $\alpha_*\left(\Omegaext \right)$ in terms of the elastic energy of $\partial \Omega$. While Theorem \ref{theo:weinstock}  does not settle \cite[Conjecture 1.2]{krejcirik2023optimisation}, it constitutes a step towards this conjecture. Namely, it implies, for a  bounded simply connected domain $\Omega \subset \mathbb{R}^2$, with $|\Omega| = |B_R|$ and $\Omega \neq B_R$, and with $\alpha =  -\sigma_2\left(B_R^\ext\right)$, that
\[
\lambda_2^{\alpha} \left( \Omegaext \right) < 0 = \lambda_2^{\alpha}\left( B_R^\ext \right).
\]
It appears reasonable to expect that this inequality also holds for $\alpha$  in a neighbourhood of
$-\sigma_2\left(B_R^\ext\right)$. A rigorous justification of this, however, is beyond the scope of the present discussion.

\subsubsection{Eigenvalue asymptotics}

Denote by
\[
\mathcal{N}_{\Omegaext}(\sigma)=\#\left\{k\in\mathbb{N}:\sigma_k\left(\Omegaext\right) \le \sigma\right\}
\] 
the eigenvalue counting function for the exterior problem. Extending the methods of \cite{GKLP2022} to the exterior setting one can prove the following result.
\begin{proposition}\label{thm:Asymptotics}
Let $\Omega \subset \mathbb{R}^n$ be a bounded open set  with $C^{2,\alpha}$ boundary
for some $\alpha>0$, and assume that $\Omegaext$ is connected.  Then
\begin{equation}\label{eq:Weyln}
\mathcal{N}_{\Omegaext}(\sigma)=\frac{\omega_{n-1}}{(2\pi)^{n-1}}  |\partial\Omega| \sigma^{n-1} + O\left(\sigma^{n-2}\right) \quad \text{ as } \sigma \to \infty,
\end{equation}
where $\omega_{n-1}$ denotes the volume of the unit ball in $\mathbb{R}^{n-1}$. 
\end{proposition}
As in the interior case,  the remainder estimate in \eqref{eq:Weyln} is sharp and attained for a ball.

\begin{remark}\label{rem:asymp2}
If $\Omega$ is smooth, the layer potential approach yields in a standard way that the exterior Dirichlet-to-Neumann map on $\partial\Omega$ is an elliptic  pseudodifferential operator, and the result then follows from  the general theory of such operators \cite{hormander1968spectral}. Further improvements can be obtained for smooth domains in two dimensions 
using the conformal approach and the results of  \cite{Rozenblum, girouard2014steklov}.
In particular, if $\Omega$ is a smooth simply connected planar domain, we have
\begin{equation}\label{eq:rozenas}
\sigma_{2k+1}\left(\Omegaext\right) =\sigma_{2k}\left(\Omegaext\right)+O\left(k^{-\infty}\right) =\frac{2\pi k}{|\partial\Omega|}+O\left(k^{-\infty}\right)\qquad\text{as }k\to\infty,
\end{equation}
where $O\left(k^{-\infty}\right)$ means that the remainder decays faster than any  negative power of $k$.
\end{remark}

For  domains with Lipschitz boundary,  asymptotics \eqref{eq:Weyln} is expected to hold with the remainder term of order  $o\left(\sigma^{n-1}\right)$. We prove this in dimension two using the conformal approach and  the results of \cite{KLP2023}. It is likely that it could be proved in arbitrary dimension using the  methods of  \cite{Rozen2023}.

\begin{proposition}\label{prop:Weyl2n}
Let  $\Omega \subset \mathbb{R}^2$ be a bounded open set with  Lipschitz boundary such that $\Omegaext$ is connected. Then 
\begin{equation}\label{eq:Weyl2}
\mathcal{N}_{\Omegaext}(\sigma)=\frac{ |\partial\Omega| }{\pi}\sigma + o\left(\sigma\right) \quad \text{ as } \sigma \to \infty.
\end{equation}
\end{proposition}

Propositions \ref{thm:Asymptotics} and \ref{prop:Weyl2n} are proved in \S\ref{sec:asymptproofs}.

\subsection{Outline of the paper and discussion}
We begin in \S\ref{section:examples} with examples where separation of variables allows the explicit computation of exterior Steklov eigenvalues. These are useful in guiding us towards the various formulations of the exterior Steklov problem, some already found in the existing literature,  that we present in \S\ref{sec:survey}. For planar domains we present an approach based on conformal mappings. In \S\ref{sec:layer} we develop a formulation using boundary layer potentials. The equivalence of the different formulations is established in \S\ref{sec:equivalence}; see Figures~\ref{fig:scheme3d} and~\ref{fig:scheme2d}. In \S\ref{sec:properties}, we study geometric properties of the eigenvalues and eigenfunctions using the various formulations. This section also contains the proofs of the results stated in \S\ref{subsec:spectral}, notably a Weinstock-type upper bound in two dimensions and an Escobar-type lower bound in higher dimensions, showing that isoperimetric inequalities fail in the latter setting.

As was already indicated earlier, different approaches to the exterior problem are used to prove various results throughout the paper. We provide a brief roadmap here for the convenience of the reader. The spectral theorem, Theorem \ref{thm:basic}, is obtained using the finite-energy formulation in dimensions $n\geq 3$ and using the conformal approach in dimension $n=2$. The same is true also for Theorem~\ref{thm:Courant} and Corollary~\ref{corol:simple}. Theorem~\ref{prop:low} and Corollary~\ref{coro:lowbound}  are proved using the finite energy formulation.  Theorem~\ref{theo:noisopineq} and Corollary~\ref{cor:Robin} then follow from these results and the Robin--Steklov duality. The planar Weinstock-type inequality, Theorem~\ref{theo:weinstock}, and the two-dimensional Weyl law,  Proposition~\ref{prop:Weyl2n}, rely on the conformal approach.  Finally, Proposition~\ref{thm:Asymptotics} is proved using the truncated domains approach, through comparison with mixed Steklov--Dirichlet and Steklov--Neumann problems. 

The approaches based on the Helmholtz equation and the layer potentials are not used for the proofs of spectral results. At the same time, the Helmholtz formulation gives a rigorous justification of the limiting procedure considered in \cite{grebenkov2024steklov} and thus connects the numerical study carried out in that paper with  the present work. The layer-potential formulation allows one to recast the problem in terms of boundary integral operators, which is especially useful for numerical implementations. 

Finally, we believe that all the approaches discussed above could be useful for further applications and for the study of the exterior Steklov problem in more general settings, e.g. on complete Riemannian manifolds.

\section{Examples}\label{section:examples}

In this section  we discuss  the exterior Steklov problem on domains for which one can use separation of variables, such as balls and spherical shells (see also  \cite[Section 3.2]{bandle2023shape}).  As above, we consider the dimensions $n=2$ and $n \ge 3$ separately.

\subsection{Exterior problem for a disk}\label{examp:balln2}

We begin by stating a basic but useful property of harmonic functions, which, for example, follows from \cite{axler86}. 

\begin{lemma}\label{lemma:usefulbasic}
Let $\rho>0$, and consider the disk $ B_\rho \subset \mathbb{R}^2 $. 
Any harmonic function $u:B_\rho^\mathrm{ext} \to \mathbb{R}$ admits the representation
\[
u(r,\theta)
= a_0 + b_0 \log r
+ \sum_{k=1}^\infty \left( a_k r^k + b_k r^{-k} \right) \cos(k\theta)
+ \sum_{k=1}^\infty \left( c_k r^k + d_k r^{-k} \right) \sin(k\theta),
\]
where $r \in (\rho, \infty)$, $\theta \in\mathbb{S}^1$. This series converges absolutely in $B_\rho^\mathrm{ext}$ and uniformly on any compact subset.
\end{lemma}

For  $\Omegaext=B_\rho^\ext=\{(r, \theta): \theta\in\mathbb{S}^1, r>\rho\}$, it is easy to check, using Lemma \ref{lemma:usefulbasic},  that the functions
\begin{equation}\label{eq:diskeigenf}
\begin{gathered}
u_{(0,1)}(r,\theta):=1,\\ 
u_{(\ell,1)}(r,\theta):= r^{-\ell}\sin(\ell\theta),\quad u_{(\ell,2)}(r,\theta):=r^{-\ell} \cos (\ell\theta), \qquad \ell \in \mathbb{N},
\end{gathered}
\end{equation}
are bounded solutions of the exterior Steklov problem \eqref{eq:PDEext}  which satisfy the  condition \eqref{eq:inftyr}, and their boundary traces form a basis in $L^2(\partial B_\rho)$.
The corresponding eigenvalues are given by 
\begin{equation}\label{eq:diskeigenv}
\sigma_{(0)}\left(B_\rho^\ext\right)=0,\qquad\sigma_{(\ell)}\left(B_\rho^\ext\right)=\frac{\ell}{\rho},\qquad\ell\in\mathbb{N},
\end{equation} 
and have multiplicities 
\[
d_{2,\ell}=\begin{cases}
1,\qquad&\text{if }\ell=0,\\
2,\qquad&\text{if }\ell\in\mathbb{N}.
\end{cases}
\]
They can be re-ordered to form the single non-decreasing sequence $\sigma_{k}\left(B_\rho^\ext\right)$ of eigenvalues enumerated with multiplicities via
\[
\sigma_1:=\sigma_{(0)},\qquad  \sigma_{2\ell} = \sigma_{2\ell+1} := \sigma_{(\ell)}, \qquad \ell\in\mathbb{N}.
\]

We reiterate  that the functions  $u_{(\ell, \cdot)}(r,\theta)$ are not in $L^2\left(\Omegaext\right)$ for $\ell\in\{0,1\}$, and thus the usual Sobolev space $H^1\left(\Omegaext\right)$ does not provide a suitable functional setting for the exterior problem.

Let us now illustrate how the truncated domains approach works in this case. In fact, we will do it for both Dirichlet and Neumann truncations, which,  as mentioned in Remark \ref{rem:other}, yield equivalent formulations in dimension two.

For $R>\rho$, consider two mixed problems in the truncated domain $B_{\rho, R}^\ext=B_R \setminus \overline{B_\rho}$ with the Steklov condition on $\partial B_\rho$. 

We start by imposing the Dirichlet condition on $\partial B_{R}$.  Solving the mixed Steklov--Dirichlet problem in $B_{\rho, R}^\ext$ using the separation of variables, we get the eigenfunctions
\begin{gather*}
u_{(0,1)}^\Dir(r,\theta):=1 - \frac{\log r}{\log R},\\
u_{(\ell,1)}^\Dir(r,\theta):=\left(r^{-\ell}- \frac{r^\ell}{R^{2\ell}}\right)\sin (\ell\theta),\quad u_{(\ell,2)}^\Dir(r,\theta):=\left( r^{-\ell}- \frac{r^\ell}{R^{2\ell}} \right) \cos (\ell\theta), 
\quad  \ell \in \mathbb{N},
\end{gather*}
with the corresponding eigenvalues
\begin{equation}\label{eq:diskDeigenv}
\sigma_{(0)}^\Dir\left(B_{\rho, R}^\ext\right)=\frac{1}{\rho (\log R - \log \rho)},\quad 
\sigma_{(\ell)}^\Dir\left(B_{\rho, R}^\ext\right)= \frac{\ell}{\rho} \frac{1 + \left(\frac{\rho}{R}\right)^{2\ell}}{1- \left(\frac{\rho}{R}\right)^{2\ell}},  \quad  \ell \in \mathbb{N},
\end{equation} 
of multiplicity $d_{2,\ell}$.

On the other hand, imposing the Neumann condition on $\partial B_{R}$ and solving the mixed Steklov--Neumann problem in $B_{\rho, R}^\ext$, we obtain the eigenfunctions
\begin{gather*}
u_{(0,1)}^\Neu(r,\theta):=1,\\
u_{(\ell,1)}^\Neu(r,\theta):=\left(r^{-\ell}+\frac{r^\ell}{R^{2\ell}}\right)\sin (\ell\theta),\quad u_{(\ell,2)}^\Neu(r,\theta):=\left( r^{-\ell}+\frac{r^\ell}{R^{2\ell}} \right) \cos (\ell\theta), 
\quad  \ell \in \mathbb{N},
\end{gather*}
and the eigenvalues 
\begin{equation}\label{eq:diskNeigenv}
\sigma_{(0)}^\Neu\left(B_{\rho, R}^\ext\right)=0,\qquad 
\sigma_{(\ell)}^\Neu\left(B_{\rho, R}^\ext\right)= \frac{\ell}{\rho} \frac{1 - \left(\frac{\rho}{R}\right)^{2\ell}}{1+ \left(\frac{\rho}{R}\right)^{2\ell}},  \qquad  \ell \in \mathbb{N}.
\end{equation} 

It is easily seen from \eqref{eq:diskeigenv}, \eqref{eq:diskDeigenv}, and \eqref{eq:diskNeigenv}, that for all $R>\rho$ and all $\ell\in\{0\}\cup \mathbb{N}$, we have
\[
\sigma_{(\ell)}^\Dir\left(B_{\rho, R}^\ext\right) > \sigma_{(\ell)}\left(B_{\rho}^\ext\right)>\sigma_{(\ell)}^\Neu\left(B_{\rho, R}^\ext\right),
\]
and that
\[
\lim_{R\to\infty}\sigma_{(\ell)}^\Dir\left(B_{\rho, R}^\ext\right)=\lim_{R\to\infty}\sigma_{(\ell)}^\Neu\left(B_{\rho, R}^\ext\right)=\sigma_{(\ell)}\left(B_{\rho}^\ext\right).
\]
Moreover, the eigenfunctions $u_{(\ell, j)}^\Dir$ and $u_{(\ell, j)}^\Neu$ converge pointwise to eigenfunctions  \eqref{eq:diskeigenf} as $R\to\infty$.
  
Note that the exterior Steklov spectrum of a disk \eqref{eq:diskeigenv} coincides with its interior spectrum 
(see \S\ref{subsec:conformal} for an explanation using the conformal approach).

\begin{open}\label{openproblem1} Let  $\Omega$ be a bounded planar domain such that its exterior and interior Steklov spectra coincide with the account of multiplicities. Show that $\Omega$ is a disk.
\end{open}
Note that in dimensions   $n\ge 3$  the exterior and the interior Steklov spectra never coincide in view of Remark \ref{rem:firsteig}.

\subsection{Exterior problem for Euclidean balls in dimensions $n \ge 3$}\label{examp:balln3}

In order to fix the terminology, let us recall that the spectrum of the Laplace--Beltrami operator $-\Delta_{\mathbb{S}^{n-1}}$ on the unit sphere $\mathbb{S}^{n-1}\subset \mathbb{R}^n$, $n\ge 3$, consists of eigenvalues $\alpha_{(\ell)}:=\ell(\ell+n-2)$, $\ell \in \{ 0 \} \cup \mathbb{N}$, of multiplicity
\[
d_{n,\ell}:=\binom{\ell+n-1}{n-1}-\binom{\ell+n-3}{n-1}.
\]
The elements of the eigenspace of $-\Delta_{\mathbb{S}^{n-1}}$ corresponding to $\alpha_{(\ell)}$ are known as \emph{spherical harmonics} of degree $\ell$, and one can choose an $L^2\left(\Sp^{n-1}\right)$-orthonormal basis $\left\{Y_{\ell,i}\right\}_{i=1}^{d_{n,\ell}}$ in each such eigenspace. 

We solve the exterior Steklov problem \eqref{eq:PDEext} in $B_\rho^\ext\subset \mathbb{R}^n$, $n\ge 3$, by separation of variables in spherical coordinates using spherical harmonics, to obtain the eigenvalues
\begin{equation}\label{eq:sigmaextBall}
\sigma_{(\ell)}\left(B_\rho^\ext\right):=\frac{n+\ell-2}{\rho},\qquad \ell\in\{0\}\cup \mathbb{N},
\end{equation}
of multiplicity $d_{n,\ell}$, with the corresponding eigenfunctions
\begin{equation}\label{eq:usigmaextBall}
u_{(\ell,i)}( r, \theta) =  r^{2-n-\ell} \, {Y}_{\ell, i}(\theta),\qquad \ell\in \{0\}\cup \mathbb{N}, \quad i=1,\ldots, d_{n,\ell}, 
\end{equation}
where $r\in[\rho,\infty)$  and $\theta\in\mathbb{S}^{n-1}$.

Note that, although all  the functions $u_{(\ell,i)}$ decay at infinity,  the radial eigenfunction $u_{(0,1)}(r, \theta)=r^{2-n} \notin L^2\left(\Omegaext\right)$ for $n=3$ and  $n=4$. Once again, this reflects the need for  a space larger than $H^1\left(\Omegaext\right)$ in order to provide a functional setting for the exterior Steklov problem.

Consider now the truncated domains approach, see also \cite{colbois2021sharp}.
For a mixed Steklov--Dirichlet problem in the spherical shell  $B^\ext_{\rho,R}$, with the Steklov condition imposed on $\partial B_\rho$ and the Dirichlet condition on $\partial B_R$, we get the eigenvalues
\[
\sigma_{(\ell)}^\Dir\left(B_{\rho,R}^\ext\right) = 
\frac{n+\ell-2+\ell\left(\frac{\rho}{R}\right)^{n+2\ell-2}}{\rho\left(1-\left(\frac{\rho}{R}\right)^{n+2\ell-2}\right)},
\qquad \ell\in\{0\}\cup \mathbb{N}, 
\]
of multiplicity $d_{n,\ell}$, with the corresponding  eigenfunctions
\[
u_{(\ell,i)}^\Dir( r,  \theta) = \left(r^{2-n-\ell}-  R^{2-n-2\ell} r^{\ell} \right) {Y}_{\ell,i}(\theta),  \qquad\ell\in\{0\}\cup \mathbb{N},\quad  i=1,\ldots, d_{n,\ell}.
\]
Recalling that  $n \geq 3$ and taking the limit as  $R\to \infty$,  we observe the convergence of  eigenvalues,  
\[
\lim_{R\to\infty}\sigma_{(\ell)}^\Dir\left(B_{\rho,R}^\ext\right)=\sigma_{(\ell)}\left(B_{\rho}^\ext\right), \qquad \ell\in\{0\}\cup \mathbb{N},
\]
and  pointwise convergence of eigenfunctions,
\[
\begin{gathered}
\lim_{R\to\infty}u_{(\ell,i)}^\Dir(r,  \theta)=u_{(\ell,i)}( r, \theta),\\
\ell\in\{0\}\cup \mathbb{N},\quad i=1, \dots, d_{n,\ell}, \quad r\in[\rho,\infty), \quad\theta\in\mathbb{S}^{n-1},
\end{gathered}
\]
to those  of the exterior Steklov problem obtained above.  

For a mixed Steklov--Neumann problem in the spherical shell  $B^\ext_{\rho,R}$ with the Neumann condition on $\partial B_R$, we get the eigenvalues
\[
\sigma_{(0)}^\Neu\left(B_{\rho, R}^\ext\right)=0,\qquad 
\sigma_{(\ell)}^\Neu\left(B_{\rho, R}^\ext\right)=  \frac{\ell(n+\ell-2)\left(1 - \left(\frac{\rho}{R}\right)^{2\ell+n-2}\right)}{\rho\left(\ell+(n+\ell-2)\left(\frac{\rho}{R}\right)^{2\ell}\right)},\qquad  \ell \in \mathbb{N}.
\]
of multiplicity $d_{n,\ell}$, with the corresponding eigenfunctions
\begin{gather*}
u_{(0,1)}^\Neu( r, \theta) = 1,\\
u_{(\ell,i)}^\Neu( r, \theta) =  \left(r^{2-n-\ell}+\frac{n+\ell-2}{\ell}R^{2-n-2\ell} r^{\ell} \right) \, {Y}_{\ell, i}(\theta), \qquad  \ell \in \mathbb{N}, \quad  i=1, \ldots, d_{n,\ell}.
\end{gather*}
For $\ell\ge 1$ we observe, as in the Dirichlet case in dimensions $n\ge 2$ and in the Neumann case in dimension $n=2$, the convergence of eigenvalues and the pointwise convergence of eigenfunctions, as $R\to\infty$,   to those  of the exterior Steklov problem. However, it is not the case for  the first eigenvalue and eigenfunction (corresponding to $\ell=0$).

\begin{remark} 
For arbitrary domains, we consider eigenfunction convergence in a different sense, see \S\ref{sec:AHvAtE} for details. We note also that although the higher Steklov--Neumann eigenvalues for the ball $B^\ext_{\rho,R}$ converge to the corresponding ones for the exterior Steklov problem in $B^\ext_\rho$ as $R\to\infty$, this property is specific to balls; in general, Neumann truncation is expected to yield different higher eigenvalues in the limit.
\end{remark}

\subsection{Exterior Dirichlet-to-Neumann maps and harmonic extensions for disks and balls}

It is easily seen that for any $n\ge 2$, the functions
\[
\widetilde{Y}_{\ell, i}(\theta) := \rho^{(1-n)/2} Y_{\ell, i}(\theta),\qquad \ell \in \mathbb{N}, \quad  i=1, \ldots, d_{n,\ell},
\]
form an $L^2\left(\partial B_\rho\right)$-orthonormal basis of eigenfunctions of the exterior Dirichlet-to-Neumann map $\Dext$ on $\partial B_\rho$, where for uniformity we set in the planar case
\[
Y_{0, 1}(\theta):=\frac{1}{\sqrt{2\pi}},\qquad Y_{\ell, 1}(\theta):=\frac{1}{\sqrt{\pi}}\sin(\ell\theta),\quad  Y_{\ell, 2}(\theta):=\frac{1}{\sqrt{\pi}}\cos(\ell\theta),\qquad \ell\in\mathbb{N}.
\]
Therefore, $\Dext$ acts on any $f\in  H^{\frac{1}{2}}\left(\partial B_\rho\right)$ as
\begin{equation}\label{eq:Dextballs}
\Dext f = \sum_{\ell=0}^\infty \sum_{i=1}^{d_{n,\ell}} \sigma_{(\ell)}\left(B_\rho^\ext\right) \widehat{f}_{\ell,i} \widetilde{Y}_{\ell, i}\in H^{-\frac{1}{2}}\left(\partial B_\rho\right),\quad
%\text{where}\quad
\widehat{f}_{\ell, i}:=\myscal{f, \widetilde{Y}_{\ell, i}}_{L^2\left(\partial B_\rho\right)}.
\end{equation}

We additionally have
\begin{proposition}\label{prop:ballgrad}
Consider, for $n\ge 2$, a ball $B_\rho\in\mathbb{R}^n$, $\rho>0$, and let $f\in H^{1/2}\left(\partial B_\rho\right)$. Then $\nabla\left(\mathcal{H}^\ext f\right)\in L^2\left(B_\rho^\ext\right)$.
\end{proposition}

\begin{proof} We expand $f$ in the Fourier series in the basis of spherical harmonics on $\partial B_\rho$,  to get
\[
f =\sum_{\ell=0}^\infty \sum_{i=1}^{d_{n,\ell}}  \widehat{f}_{\ell,i} \widetilde{Y}_{\ell, i},
\]
where the condition $f\in  H^{1/2}\left(\partial B_\rho\right)$ is equivalent to 
\begin{equation}\label{eq:f12}
\sum_{\ell=0}^\infty \sum_{i=1}^{d_{n,\ell}} \sqrt{1+\frac{\alpha_{(\ell)}}{\rho^2}} \left|\widehat{f}_{\ell, i}\right|^2<\infty.
\end{equation}
The harmonic extension $u:=\mathcal{H}^\ext f$ is then 
\[
u(r,\theta) =\sum_{\ell=0}^\infty \sum_{i=1}^{d_{n,\ell}} \widehat{f}_{\ell,i} r^{2-n-\ell} \widetilde{Y}_{\ell, i},
\]
and explicit integration gives
\begin{equation}\label{eq:gradnormball}
\left\|\nabla u\right\|^2_{L^2\left(B_\rho^\ext\right)} = \sum_{\ell=0}^\infty \sum_{i=1}^{d_{n,\ell}}  \sigma_{(\ell)}\left(B_\rho^\ext\right) \left|\widehat{f}_{\ell, i}\right|^2,
\end{equation}
which is finite by \eqref{eq:f12} due to the fact that 
\[
\lim_{\ell\to\infty}\frac{\sigma_{(\ell)}\left(B_\rho^\ext\right)}{\sqrt{1+\frac{\alpha_{(\ell)}}{\rho^2}}}=1.
\]
\end{proof}

\subsection{Numerical examples}\label{sec:numerics}

\begin{figure}[t]
\centering
\includegraphics[width=0.95\textwidth]{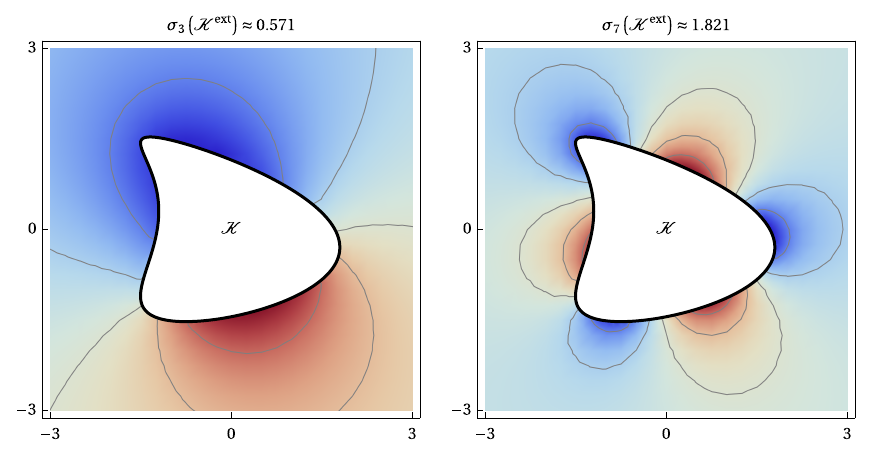}
\caption{Density plots of eigenfunctions of the Steklov problem in the exterior $\mathcal{K}^\ext$ of the kite.}\label{fig:kite}
\end{figure}

Let $\mathcal{K}\subset\mathbb{R}^2$ be an asymmetric bounded ``kite'' domain whose boundary is given parametrically by $\partial\mathcal{K}:=\left\{\left(1.5\cos t + 0.7\cos 2t - 0.4, 1.5\sin t - 0.3\cos t\right), t\in[0,2\pi)\right\}$. We show, in Figure \ref{fig:kite}, numerically computed eigenfunctions of the exterior Steklov problem in $\mathcal{K}^\ext$ corresponding to its third and seventh eigenvalues. %

To illustrate the differences between interior and exterior Steklov problems, we compare some numerically computed eigenvalues of  $\mathcal{K}$ and $\mathcal{K}^\ext$  in Table \ref{tab:K}; note that, as always in the planar case,  we have $\sigma_1\left(\mathcal{K}\right)=\sigma_1\left(\mathcal{K}^\ext\right)=0$.

\begin{table}[t]
\centering
\begin{tabular}{@{}>{$}r<{$}@{\hspace{1cm}}>{$}r<{$}@{\hspace{0.5cm}}>{$}r<{$}@{}}
\toprule
k & \sigma_k\left(\mathcal{K}\right) & \sigma_k\left(\mathcal{K}^\ext\right)\\
\midrule
2 & 0.403 & 0.545\\
3 & 0.524 & 0.571\\ 
4 & 1.183 & 1.130\\
5 & 1.384 & 1.309\\
6 & 1.721 & 1.746\\ 
7 & 2.018 & 1.821\\ 
8 & 2.201 & 2.293\\ 
9 & 2.706 & 2.450\\ 
10 & 2.785 & 2.903\\ 
\bottomrule\\
\end{tabular}
\caption{Numerically computed Steklov eigenvalues of  $\mathcal{K}$ and $\mathcal{K}^\ext$.}\label{tab:K}
\end{table}

\begin{figure}[t]
\centering
\includegraphics[width=0.95\textwidth]{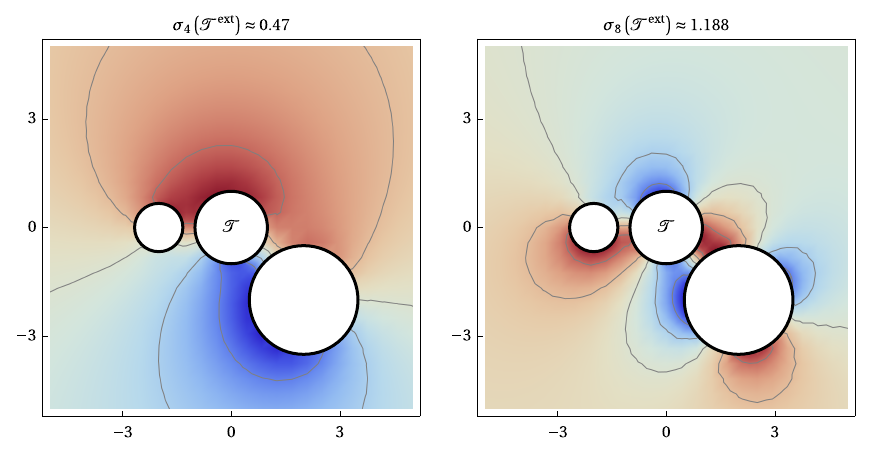}
\caption{Density plots of eigenfunctions of the Steklov problem in the exterior $\mathcal{T}^\ext$ of the disjoint union of three disks.}\label{fig:disks}
\end{figure}

Additionally, the numerically computed eigenfunctions of the Steklov problem (corresponding to its fourth and eighth eigenvalues)  in the exterior of the disjoint union $\mathcal{T}\subset \mathbb{R}^2$ of the unit disk centred at the origin, the disk of radius $\frac{2}{3}$ centred at $(-2,0)$, and the disk of radius $\frac{3}{2}$ centred at $(2,-2)$, are shown in Figure \ref{fig:disks}.

For both geometries, the exterior eigenvalues and eigenfunctions are computed using the conformal transform method of \S\ref{subsec:conformal} and the finite element calculations in FreeFEM \cite{freefem} for the corresponding interior domain problem, see the scripts cited in the Data availability statement on page \pageref{sec:data}. 

An alternative numerical approach to the exterior Steklov problem, see, e.g., \cite{grebenkov2024steklov, ChaGr24}, can be realised via the domain decomposition technique by reducing the problem to a mixed Steklov--matching condition problem in a bounded domain. Let $\Omega\subset\mathbb{R}^n$ be a bounded open set with Lipschitz boundary and connected $\Omegaext$, and fix any $\rho>R_0(\Omega)$. Then it is easy to see that the problem \eqref{eq:PDEext} is equivalent to 
\begin{equation*}\label{eq:SPsi}\tag{$\mathrm{S}\Psi_\rho$}
\begin{cases}
\Delta u  = 0\qquad&\text{ in } \Omegaext_\rho, \\
\partial_\nu u  = \sigma u\qquad&\text{ on } \partial\Omega,\\
\left(\partial_{\overline{\nu}} + \Dext\right)u =0\qquad&\text{ on } \partial B_\rho,
\end{cases}
\end{equation*}
where the exterior Dirichlet-to-Neumann map $\Dext$ appearing in the pseudodifferential boundary condition on $\partial B_\rho$ is given by \eqref{eq:Dextballs}. In addition to a possible numerical implementation, the problem \eqref{eq:SPsi} can be viewed as yet another equivalent formulation of the exterior Steklov problem.

\section{Different approaches to the exterior Steklov problem}\label{sec:survey}

 In this section we  introduce some of the existing approaches to the exterior Steklov problem. 

\subsection{The Steklov problem for finite-energy functions, $n \ge 3$}\label{subsec:AuchmutyHan}

We start with the global approach introduced by Arendt and ter Elst~\cite{arendt2015dirichlet} for the Dirichlet-to-Neumann operators,  as well as by Auchmuty and Han~\cite{auchmuty2014representations} for the weak formulation of the Steklov problem. We work directly in function spaces over the full exterior domain. It is in this setting that we will present the spectral theorem for the exterior problem under consideration. It will be considered as a reference problem, that will later be  compared to various other methods. The planar case is different and will be discussed separately.

\subsubsection{Dirichlet-to-Neumann operators}\label{sec:DtNs}
For $n\geq 3$, Arendt and ter Elst~\cite{arendt2015dirichlet} introduced two exterior Dirichlet-to-Neumann operators $A^\Neu$ and $A^\Dir$ on $L^2(\partial\Omega)$ using the following two function spaces,
\[
W\left(\Omegaext\right) := \left\{ u \in H^1_\mathrm{loc}\left(\Omegaext\right) : \int_{\Omegaext} |\nabla u|^2 \, \mathrm{d}x < \infty \right\}
\]
and 
\[
W^\Dir\left(\Omegaext\right) := W\left(\Omegaext\right) \cap L^p\left(\Omegaext\right),
\]
where $p:=\frac{2n}{n-2}>2$, and therefore $\frac{1}{p}=\frac{1}{2}-\frac{1}{n}$.

By \cite[Proposition 2.6]{arendt2015dirichlet},
\[
W^\Dir\left(\Omegaext\right)=\left\{u\in W\left(\Omegaext\right)\,:\,\langle u\rangle=0\right\},
\]
where
\begin{equation}\label{eq:averageu}
\langle u\rangle =\lim_{R\to \infty}\frac{1}{\left|\Omegaext_R\right|}\int_{\Omegaext_R}u\,\dr x
\end{equation}
is the average of $u$ over $\Omegaext$. The existence of this limit follows from a similar result in the full space $\R^{n}$ \cite[Theorem 1.1]{LuOu}, combined with the uniform boundedness of the extension operator. This is proved in detail in \cite[Section 2]{arendt2015dirichlet}, where  properties of the relevant spaces are discussed at length.
Because $p>2$,  it follows from the Sobolev embedding theorem that $H^1\left(\Omegaext\right)\subset W^\Dir\left(\Omegaext\right)$.  This inclusion is strict. Indeed, there are functions in $W^\Dir\left(\Omegaext\right)$ that are not in $L^2\left(\Omegaext\right)$ as their decay at infinity is only fast enough to guarantee that they belong to $L^p\left(\Omegaext\right)$.

The space $W\left(\Omegaext\right)$, provided with the norm
\[
\| u \|_{W\left(\Omegaext\right)} := \left( \int_{\Omegaext} |\nabla u|^2 \, \mathrm{d}x + |\langle u \rangle|^2 \right)^{1/2},
\] 
is complete, the trace operator $\boldsymbol{\tau}:W\left(\Omegaext\right) \to L^2(\partial\Omega)$, $\boldsymbol{\tau}: u\mapsto u|_{\partial\Omega}$,  is compact, and the norm defined above is equivalent to $u \mapsto \left(\|\nabla u\|_{L^2\left(\Omegaext\right)}^2+\|u\|_{L^2(\partial\Omega)}^2 \right)^{1/2}$, see~\cite[Proposition 2.5]{arendt2015dirichlet}.

\begin{definition}\label{def:ANAD}
The \emph{Dirichlet-to-Neumann operator with Neumann condition at infinity} ,
\[
A^\Neu: \operatorname{Dom}\left(A^\Neu\right)\subset L^2(\partial\Omega)\to L^2(\partial\Omega),
\] 
is the self-adjoint unbounded operator associated with $\left(W\left(\Omegaext\right),Q,\boldsymbol{\tau}\right)$, where the quadratic form \newline$Q: W\left(\Omegaext\right)\to\mathbb{R}$ is given by
\[
Q[u] := \int_{\Omegaext} |\nabla u|^2 \, \dr x.
\]
Similarly, the \emph{Dirichlet-to-Neumann operator with Dirichlet condition at infinity}, 
\[
A^\Dir: \operatorname{Dom}\left(A^\Dir\right)\subset L^2(\partial\Omega)\to L^2(\partial\Omega),
\] 
is the self-adjoint operator on $L^2(\partial\Omega)$ associated with $\left(W^\Dir\left(\Omegaext\right),Q,\boldsymbol{\tau}\right)$, see \cite{arendt2015dirichlet}.
\end{definition}

The following result follows from \cite[Theorem 5.9 and the remark before Proposition 5.3]{arendt2015dirichlet}.

\begin{theorem}\label{thm:SpectralADAN}
The self-adjoint operators $A^\Neu$ and $A^\Dir$ have the same domain  $\operatorname{Dom}\left(A^\Neu\right)=\operatorname{Dom}\left(A^\Dir\right)\subset L^2(\partial\Omega)$, and
\[
A^\Dir h=A^\Neu h+\frac{1}{b}\myscal{v,A^\Dir1}_{L^2(\partial\Omega)}(A^\Dir1),
\]
where $b:=\int_{\partial\Omega}A^\Dir1>0$. Thus, $A^\Dir$ is a bounded rank-one perturbation of $A^\Neu$.
Both operators have compact resolvents and their spectra are purely discrete and non-negative,  $0=\lambda_1\left(A^\Neu\right)\leq\lambda_2\left(A^\Neu\right)\leq\cdots\nearrow+\infty,$ and
$0<\lambda_1\left(A^\Dir\right)\leq\lambda_2\left(A^\Dir\right)\leq\cdots\nearrow+\infty$.  In both cases, the corresponding eigenfunctions can be chosen to form a complete orthonormal set in $L^2(\partial\Omega)$.
\end{theorem}

We additionally have the following interlacing result.

\begin{theorem}\label{thm:interlaceADAN} For any $k\in\mathbb{N}$, 
\[
\lambda_k\left(A^\Neu\right) \le \lambda_k\left(A^\Dir\right) \le \lambda_{k+1}\left(A^\Neu\right).
\]  
\end{theorem}

The proof of Theorem \ref{thm:interlaceADAN} (presented below) relies on the following abstract result generalising the finite-dimensional Cauchy interlacing theorem \cite[Chapter 7]{bellman1997introduction}. Although the result is elementary, we provide a proof, as we have been unable to locate a reference for unbounded operators.

\begin{proposition}\label{prop:interlace} 
Let $\mathcal{A}$ and $\mathcal{B}$ be two self-adjoint semi-bounded below operators in a Hilbert space $H$ with the same domain $D:=\operatorname{Dom}(\mathcal{A})= \operatorname{Dom}(\mathcal{B})$ and with discrete spectra, and suppose that there exists ${h}_0\in D$ such that 
\[
\mathcal{B}{h} = \mathcal{A}{h} + \myscal{h, h_0}_H {h}_0\qquad \text{for all }{h}\in D.
\]
Then the eigenvalues of $\mathcal{A}$ and $\mathcal{B}$ interlace,
\begin{equation}\label{eq:interlace}
\lambda_k( \mathcal{A}) \le \lambda_k( \mathcal{B}) \le \lambda_{k+1}( \mathcal{A})\qquad \text{for all }k\in\mathbb{N}.
\end{equation}
\end{proposition}

\begin{proof}[Proof of Proposition \ref{prop:interlace}] 
First, for any ${h}\in D$, 
\[
\myscal{\mathcal{B}{h},{h}}_H = \myscal{\mathcal{A}{h},{h}}_H + \myscal{{h},{h}_0}_H^2\ge \myscal{\mathcal{A}{h},{h}}_H,
\]
immediately establishing the first inequality in \eqref{eq:interlace}.

Let $\mathcal{R}_\mathcal{A}[{h}]$  and $\mathcal{R}_\mathcal{B}[{h}]$ be the Rayleigh quotients of $\mathcal{A}$ and $\mathcal{B}$ evaluated at ${h}\in D\setminus\{0\}$.  Fix $k\in\mathbb{N}$, and recall the variational principle
\[
\lambda_{k+1}( \mathcal{A}) = \min_{\substack{\mathcal{L}\subset D\\\dim\mathcal{L}=k+1}}\ \max_{\substack{{h}\in\mathcal{L}\\{h}\ne 0}}\ \mathcal{R}_\mathcal{A}[{h}],\qquad
\lambda_{k}(\mathcal{B}) = \min_{\substack{\mathcal{M}\subset D\\\dim\mathcal{M}=k}}\ \max_{\substack{{h}\in\mathcal{M}\\{h}\ne 0}}\ \mathcal{R}_\mathcal{B}[{h}].
\]

Take any subspace $\mathcal{L}\subset D$ with $\dim\mathcal{L}=k+1$, and set $\widetilde{\mathcal{L}}:=\left\{h\in \mathcal{L}: \myscal{h, h_0}_H=0\right\}$. Then $\dim\widetilde{\mathcal{L}}\ge k$. Choose now any subspace $\mathcal{M}\subset \widetilde{\mathcal{L}}$, $\dim \mathcal{M}=k$. Then restrictions of the operators $\mathcal{A}$ and $\mathcal{B}$ to $\mathcal{M}$ coincide, and therefore, as $\mathcal{M}\subset \mathcal{L}$, 
\begin{equation}\label{eq:RABM}
\max_{\substack{{h}\in\mathcal{M}\\{h}\ne 0}}\ \mathcal{R}_\mathcal{B}[{h}] = \max_{\substack{{h}\in\mathcal{M}\\{h}\ne 0}}\ \mathcal{R}_\mathcal{A}[{h}] \le  \max_{\substack{{h}\in\mathcal{L}\\{h}\ne 0}}\ \mathcal{R}_\mathcal{A}[{h}].  
\end{equation}

Let us now take, in the right- and left-hand sides of \eqref{eq:RABM}, the minimum over all subspaces $\mathcal{L}\subset D$, $\dim\mathcal{L}=k+1$, and all  $\mathcal{M}\subset\widetilde{\mathcal{L}}$, $\dim \mathcal{M}=k$. As the right-hand side of \eqref{eq:RABM} is independent of $\mathcal{M}$, we get 
\[
\min_{\substack{\mathcal{L}\subset D\\\dim\mathcal{L}=k+1\\\mathcal{M}\subset \widetilde{\mathcal{L}}\\\dim \mathcal{M}=k}}\ \max_{\substack{{h}\in\mathcal{L}\\{h}\ne 0}}\ \mathcal{R}_\mathcal{A}[{h}]=
\min_{\substack{\mathcal{L}\subset D\\\dim\mathcal{L}=k+1}}\ \max_{\substack{{h}\in\mathcal{L}\\{h}\ne 0}}\ \mathcal{R}_\mathcal{A}[{h}] = \lambda_{k+1}( \mathcal{A}).
\]
At the same time, minimising in the left-hand side over a wider choice of all subspaces  $\mathcal{M}\subset D$ with $\dim \mathcal{M} = k$, we get 
\[
\min_{\substack{\mathcal{L}\subset D\\\dim\mathcal{L}=k+1\\\mathcal{M}\subset \widetilde{\mathcal{L}}\\\dim \mathcal{M}=k}}\  \max_{\substack{{h}\in\mathcal{M}\\{h}\ne 0}}\ \mathcal{R}_\mathcal{B}[{h}] \ge 
\min_{\substack{\mathcal{M}\subset D\\\dim\mathcal{M}=k}}\ \max_{\substack{{h}\in\mathcal{M}\\{h}\ne 0}}\ \mathcal{R}_\mathcal{B}[{h}] = \lambda_{k}(\mathcal{B}), 
\]
which together imply the second inequality in \eqref{eq:interlace}.
\end{proof}

\begin{proof}[Proof of Theorem \ref{thm:interlaceADAN}] 
We apply Proposition \ref{prop:interlace} with $\mathcal{A}=A^\Neu$, $\mathcal{B}=A^\Dir$, $H=L^2(\partial\Omega)$, and ${h}_0=\frac{A^\Dir1}{\sqrt{b}}$, which we can do since $b = \myscal{A^\Dir1, 1}_{L^2(\partial\Omega)}\ge \lambda_1\left(A^\Dir\right)>0$, and 
\[
A^\Dir{h} = A^\Neu{h} + \myscal{{h},{h}_0}_{L^2(\partial\Omega)}{h}_0
\]
by Theorem \ref{thm:SpectralADAN}.
\end{proof} 

\subsubsection{Weak formulation of the Steklov problem}\label{subsec:fen}
It is often convenient to study the Steklov eigenvalue problem directly at the level of quadratic forms, in their weak formulation. For exterior domains, this approach was used by Auchmuty and Han in~\cite{auchmuty2014representations}, where they introduced
a seemingly different function space $E^1\left(\Omegaext\right)$ of real-valued finite-energy functions that decay at infinity~\cite[Section 3]{auchmuty2014representations}. It is the natural analogue, for exterior domains, of the space $D^1(\mathbb{R}^n)$ studied in  \cite[Chapter 8]{lieb2001analysis}, and consists of Lebesgue measurable functions $u: \Omegaext \to \mathbb{R}$ satisfying the following three conditions:
\begin{enumerate}[(i)]
\item $u \in L^1\left(\Omegaext_R\right)$ for any $R>R_0(\Omega)$,
\item $|\nabla u| \in L^2\left(\Omegaext\right)$,
\item $\left\{ x \in \Omegaext: |u(x)| \geq c \right\}$ has finite Lebesgue measure for any $c>0$.
\end{enumerate}
However, as a consequence of \cite[Corollary 3.4]{auchmuty2014representations}, this is the real version of $W^\Dir\left(\Omegaext\right)$,
\begin{equation}\label{eq:WD=E1}
E^1\left(\Omegaext\right)=\left\{f\in W^\Dir\left(\Omegaext\right)\,:\,f\text{ is real valued }\right\},
\end{equation}
for any bounded $\Omega \subset \mathbb{R}^n$, $n\geq 3$ with Lipschitz boundary and connected $\Omegaext=\mathbb{R}^n\setminus\overline{\Omega}$. Equipped with the inner product introduced in \cite[formula (3.5)]{auchmuty2014representations},
\[
\myscal{ f,g }_{E^1\left(\Omegaext\right)} := \int_{\Omegaext} \mydotp{\nabla f, \nabla g}\, \mathrm{d}x + \frac{1}{|\partial \Omega|} \int_{\partial \Omega} f g \, \mathrm{d}S,
\]
 $E^1\left(\Omegaext\right)$ becomes a Hilbert space and the associated norm is equivalent to the norm $\| \cdot \|_{W^\Dir\left(\Omegaext\right)}$.

The following theorem follows from the characterisation~\cite[Proposition 5.1]{arendt2015dirichlet} of the operator domain $\operatorname{Dom}\left(A^\Dir\right)$ and \cite[Section 9]{auchmuty2014representations}. 

\begin{theorem}\label{thm:spectral}
Let $n\geq 3$. Suppose $\Omega \subset \mathbb{R}^n$ is a bounded open set with Lipschitz boundary and with connected $\Omegaext$. A  function $f\in L^2(\partial\Omega)$ 
satisfies $A^\Dir f=\lambda f$ if and only $f$ is the trace $\boldsymbol{\tau}(u)$ of a function $u\in E^1\left(\Omegaext\right)$ such that
\begin{equation}\label{WeakFormulationExterior}
\int_{\Omegaext}\mydotp{\nabla u,\nabla v}\,\dr x=\lambda\int_{\partial\Omega}uv\,\dr S\qquad\text{for all }v\in E^1\left(\Omegaext\right).\tag{$\mathrm{FE}$}
\end{equation}
Moreover for each $k\geq 1$,
\begin{equation}\label{eq:mu_k}
\lambda_k\left(A^\Dir\right) = \min_{\substack{M \subset E^1\left(\Omegaext\right)\\ \dim(M)=k}} \max_{\substack{u \in M\\u\ne 0}} \frac{\int_{\Omegaext} | \nabla u |^2 \, \dr x}{\int_{\partial\Omega} |  u |^2 \, \dr S}.
\end{equation}
Additionally, the corresponding eigenfunctions $u_k$ can be chosen to form a complete orthonormal basis of the subspace
$\left\{ u \in E^1\left(\Omegaext\right): \Delta u=0\right\}$  and the restrictions of these eigenfunctions to $\partial\Omega$ form an orthonormal basis of $L^2(\partial\Omega)$.
\end{theorem}

The positivity of the first eigenvalue $\lambda_1$ that was claimed in Theorem~\ref{thm:SpectralADAN} also immediately follows from the fact that the constant function is not an element of $E^1\left(\Omegaext\right)$.

\begin{remark}\label{rem:not2d}
The approach presented in this section does not extend to dimension $ n = 2 $. In this case, the space $ E^1\left(\Omegaext\right) $ is no longer complete with respect to the $\| \cdot \|_{E^1\left(\Omegaext\right)}$-norm; see \cite[Appendix A]{auchmuty2014p}.

Even if the decay condition is weakened to $u \in L^\infty\left(\Omegaext\right)$, the resulting space is not complete in dimension two. For example, consider the sequence $\left(f_m\right)_{m \in \mathbb{N}} \subset E^1\left(B_1^\ext\right)$ defined by
\[
f_m(x) := \begin{cases}
0 \qquad&\text{if } |x| < \er, \\
\log \log |x| \qquad&\text{if } \er \leq |x| < m, \\
\log \log m  \qquad&\text{if } |x| \geq m,
\end{cases}
\]
which defines a Cauchy sequence in $ E^1\left(B_1^\ext\right) $, whose limit is not in $ L^\infty\left(\Omegaext\right) $.
\end{remark}

\begin{remark}
Another potential approach would be to use weighted Sobolev spaces to enforce some uniqueness of harmonic extensions,  see, for instance, \cite{AGG}.
\end{remark}

\subsection{Truncated domains in dimensions $n\geq 2$}\label{subsec:truncatedomains}

Arendt and ter Elst \cite{arendt2015dirichlet} proposed two ways  to approximate the Dirichlet-to-Neumann operator on exterior domains by truncated domains. For any $R>R_0(\Omega)$, recalling that $\Omegaext_R = \Omegaext \cap B_R$,  we consider the Dirichlet approximation,
\begin{equation*}\label{eq:PDEArendtElstDirichlet}
\begin{cases}
\Delta u_R  = 0\qquad&\text{ in } \Omegaext_R, \\
\partial_\nu u_R  = \sigma^\Dir u_R\qquad&\text{ on } \partial\Omega,\\
 u_R  = 0\qquad&\text{ on } \partial B_R,
\end{cases}\tag{$\mathrm{SD}_R$}
\end{equation*}
with the spectral parameter $\sigma^\Dir=\sigma^\Dir\left(\Omegaext_R\right)$, and the Neumann approximation,
\begin{equation*}\label{eq:PDEArendtElstNeumann}
\begin{cases}
\Delta u_R  = 0\qquad&\text{ in } \Omegaext_R, \\
\partial_\nu u_R  = \sigma^\Neu u_R\qquad&\text{ on } \partial\Omega,\\
\partial_\nu u_R  = 0\qquad&\text{ on } \partial B_R,
\end{cases}\tag{$\mathrm{SN}_R$}
\]
with the spectral parameter $\sigma^\Neu=\sigma^\Neu\left(\Omegaext_R\right)$. 

We can interpret~\eqref{eq:PDEArendtElstDirichlet} as a spectral problem for the self-adjoint operator $A_R^\Dir$ in $L^2(\partial\Omega)$ that is associated to the quadratic form 
\[
Q_R^\Dir[u] = \int_{\Omegaext_R} |\nabla u|^2 \, \mathrm{d}x, 
\]
defined on
\[
\operatorname{Dom}\left(Q_R^\Dir\right) = H_D^1\left(\Omegaext_R\right):= \left\{ u \in H^1\left(\Omegaext_R\right): u|_{\partial B_R} = 0 \right\},
\]
and to the compact trace operator $\boldsymbol{\tau}:H^1\left(\Omegaext_R\right)\to L^2(\partial\Omega)$. 
Similarly, we understand~\eqref{eq:PDEArendtElstNeumann} as the spectral problem for the self-adjoint operator $A_R^\Neu$ in $L^2(\partial\Omega)$ associated to the quadratic form 
\[
Q_R[u] = \int_{\Omegaext_R} |\nabla u|^2 \, \mathrm{d}x, \quad \operatorname{Dom}(Q_R) = H^1\left(\Omegaext_R\right)
\]
and the trace operator $\boldsymbol{\tau}$. 
These are well-known interior mixed problems:
the spectra of both operators are purely discrete and their eigenvalues can be characterised by
\begin{align}
\sigma_k^\Dir\left(\Omegaext_R\right) &= \inf_{\substack{M \subset H^1_D\left(\Omegaext_R\right) \\ \dim(M)=k}}\ \sup_{\substack{u\in M\\u\ne 0}} \frac{\int_{\Omegaext_R} | \nabla u |^2 \, \mathrm{d}x}{\int_{\partial\Omega} |  u |^2 \, \mathrm{d}S},\label{eq:sigma_kD}\\
\sigma_k^\Neu\left(\Omegaext_R\right) &= \inf_{\substack{M \subset H^1\left(\Omegaext_R\right) \\ \dim(M)=k}}\  \sup_{\substack{u\in M\\u\ne 0}} \frac{\int_{\Omegaext_R} | \nabla u |^2 \, \mathrm{d}x}{\int_{\partial\Omega} |  u |^2 \, \mathrm{d}S}.\label{eq:sigma_kN}
\end{align}

The following monotonicity results follow directly from~\eqref{eq:sigma_kD} and \eqref{eq:sigma_kN}.

\begin{lemma}\label{lemma:monotonicitymixed}
Let $n\geq 2$. Suppose $\Omega \subset \mathbb{R}^n$ is a bounded open set with Lipschitz boundary and with connected $\Omegaext$. 
For any  $R>R_0(\Omega)$ we have 
\[
\sigma_k^\Neu\left(\Omegaext_{R}\right)\le \sigma_k^\Dir\left(\Omegaext_{R}\right),\qquad k \in \mathbb{N},
\]
and for any $R_1>R>R_0(\Omega)$, we additionally have
\[
\sigma_k^\Neu\left(\Omegaext_{R_1}\right)\geq \sigma_k^\Neu\left(\Omegaext_{R}\right),\qquad
\sigma_k^\Dir\left(\Omegaext_{R_1}\right)\leq \sigma_k^\Dir\left(\Omegaext_{R}\right),\qquad k \in \mathbb{N}.
\] 
Therefore, both limits $\lim\limits_{R \to \infty}\sigma_k^\Neu\left(\Omegaext_R\right)$ and $\lim\limits_{R \to \infty}\sigma_k^\Dir\left(\Omegaext_R\right)$ exist.
\end{lemma}

For  dimensions $n\geq 3$, it was already proved in~\cite{arendt2015dirichlet} that the truncated domain approaches are equivalent to the approaches from \S\ref{sec:DtNs}. 

\begin{theorem}[{\cite[Theorems 5.5 and 5.6]{arendt2015dirichlet}}]\label{thm:resolvconv}
For $n \geq 3$, the resolvent of $A_R^\Dir$ and the resolvent of $A_R^\Neu$ converge uniformly to the resolvents of the operators $A^\Dir$ and $A^\Neu$,  respectively.
Thus, $\lim\limits_{R \to \infty}\sigma_k^\Dir\left(\Omegaext_R\right)=\lambda_k\left(A^\Dir\right)$ and $\lim\limits_{R \to \infty}\sigma_k^\Neu\left(\Omegaext_R\right)=\lambda_k\left(A^\Neu\right)$.
\end{theorem}

Therefore, by Theorem \ref{thm:spectral}, the truncated Dirichlet approach is equivalent to the energy space approach. For the analysis in the case of planar domains, we refer to \S\ref{section:twod}.

\subsection{Helmholtz equation in dimension $n\geq 2$}\label{subsec:Greben}

An alternative approach to obtaining a discrete spectrum, followed by Grebenkov and Chaigneau in~\cite{grebenkov2024steklov}, is to  introduce a fixed parameter $\Lambda>0$ and consider the spectral problem
\begin{equation*}\label{eq:PDEGrebChai}
\begin{cases}
(\Lambda^2-\Delta) u  = 0\qquad&\text{ in } \Omegaext, \\
\partial_\nu u  = \mu u\qquad&\text{ on } \partial\Omega,
\end{cases}\tag{$\mathrm{H}_\Lambda$}
\end{equation*}
where $\mu=\mu\left(\Omegaext,\Lambda\right)$ is the spectral parameter. Since $\Lambda^2>0$, the corresponding bilinear form is $H^1\left(\Omegaext\right)$-coercive. In~\cite[Sections 4 and 8]{auchmuty2012bases},  it is demonstrated that coercivity can be used to show that the spectrum of a Steklov problem on bounded domains is purely discrete.  Similarly, by applying the same modifications to \cite{auchmuty2012bases} as in, e.g., \cite[Theorem 5.1]{auchmuty2013spectral}, we deduce that the spectrum of  \eqref{eq:PDEGrebChai} is purely discrete, with eigenvalues given by
\begin{equation}\label{eq:lambdak}
\mu_k\left(\Omegaext,\Lambda\right) = \min_{\substack{M \subset H^1\left(\Omegaext\right) \\ \dim(M)=k}} \max_{\substack{u \in M\\u\ne 0}} \frac{\int_{\Omegaext} |\nabla u|^2 \, \mathrm{d}x + \Lambda^2\int_{\Omegaext} |u|^2 \,\dr x}{\int_{\partial\Omega} |u|^2 \, \dr S}.
\end{equation}
The corresponding eigenfunctions belong to $H^1\left(\Omegaext\right)$, and, in addition, all eigenfunctions are exponentially decaying as $|x|\to\infty$ due to

\begin{lemma}[\cite{BardosMerigot1977}]\label{lem:Bardos}
Let $\Omega \subset \R^n$, $n \geq 2$, be a bounded domain with a Lipschitz boundary $\partial\Omega$ and $\Lambda>0$. For any function $u\in H^1\left(\Omegaext\right)$ which satisfies $(\Lambda^2-\Delta)u=0$ in $\Omegaext$, there exists a constant $C>0$ such that
\[
|u(x)| \leq C |x|^{\frac{n-1}{2}} \er^{-\Lambda |x|}
\]
for sufficiently large $|x|$.
\end{lemma}

Since $\mu_k\left(\Omegaext,\Lambda\right)$ is monotonically increasing in $\Lambda$ for $\Lambda>0$, the limit 
\[
\mu_k\left(\Omegaext\right) := \lim_{\Lambda\searrow 0}\mu_k\left(\Omegaext,\Lambda\right)
\]
exists. The asymptotic behaviour of $\mu_k\left(\Omegaext,\Lambda\right)$ has been studied in \cite{grebenkov2024steklov}, providing both theoretical results on the convergence rate and numerical examples. 

\begin{example}\label{examp:ball3n}
For a given $\Lambda>0$, functions satisfying $(\Lambda^2-\Delta ) f = 0 $ in the exterior of $\Omega = B_\rho \subset \mathbb{R}^n $ with $n \geq 2$ are of the form
\[
f(r, \theta) = \sum_{\ell=0}^\infty \sum_{i=1}^{d_{n,\ell}} r^{-\frac{n-2}{2}} \left( b_{\ell,i} K_\frac{n+2\ell-2}{2}(\Lambda r)+ c_{\ell,i} I_\frac{n+2\ell-2}{2}(\Lambda r) \right) {Y}_{\ell,i}(\theta),
\]
where $I_m(\cdot)$ and $K_m(\cdot)$ are the modified Bessel functions of the first and second kind, respectively, of order $m$, and $ b_{\ell,i},  c_{\ell,i}$ are constants. In view of the asymptotic behaviour  of the Bessel functions (see, e.g., \cite[Chapter 10]{dlmf}), we have to choose $c_{\ell,i}=0$ in order to ensure the decay at infinity. Since $K_m(\Lambda r)$ decays exponentially for $r \to \infty$, the resulting function belongs to $H^1\left(\Omegaext\right)$. So, the Steklov eigenvalue problem \eqref{eq:PDEGrebChai} in $B_\rho^\ext$ has the eigenfunctions (with arbitrary scaling factors $b_{\ell,i}$)
\begin{equation}\label{eq:EFHelmHoltzBall}
\begin{gathered}
u_{(\ell,i)}( r, \theta,\Lambda) = b_{\ell,i}  r^{1-\frac{n}{2}}  K_{\ell+\frac{n}{2}-1}(\Lambda r) {Y}_{\ell,i}(\theta),\\ 
\ell\in\{0\}\cup\mathbb{N},\quad i=1, \dots, d_{n,\ell}, 
\end{gathered}
\end{equation}
associated to the eigenvalues
\[
\mu_{(\ell)} \left( B_\rho^\ext,\Lambda \right) =\frac{\Lambda K_{\ell+\frac{n}{2}}(\Lambda \rho)}{K_{\ell+\frac{n}{2}-1}(\Lambda \rho)}-\frac{\ell}{\rho},\qquad\ell\in\{0\}\cup\mathbb{N},
\]
of multiplicity $d_{n,\ell}$. 

We now consider the cases of  balls of dimension $n\ge 3$ and disks in the planar case $n=2$ separately.

We start with the case $n\ge 3$. Using $K_m(x) = \frac{\Gamma(m)2^m}{2 x^m}+o\left(x^{-m}\right)$ as $x \to 0$ for $m > 0$, we obtain 
\[
\lim_{\Lambda\searrow 0} \mu_{(\ell)} \left( B_\rho^\ext, \Lambda \right) =  \frac{n+2\ell-2}{\rho}-\frac{\ell}{\rho}=\frac{n+\ell-2}{\rho},
\]
which are exactly the eigenvalues given by \eqref{eq:sigmaextBall}.
Choosing the normalising constants $b_{\ell,i} = \frac{\left(\frac{\Lambda}{2}\right)^{\ell+\frac{n}{2}-1}}{2\Gamma(\ell+\frac{n}{2}-1)}$ in \eqref{eq:EFHelmHoltzBall}, we get
\[
\lim_{\Lambda \searrow 0} u_{(\ell, i)}(r, \theta, \Lambda) =  r^{-(n+\ell-2)}  {Y}_{\ell,i}(\theta),
\] 
see \eqref{eq:usigmaextBall}. So, we recover the original eigenvalues and eigenfunctions as $\Lambda \searrow 0$.

For $n=2$, we have the same form of the  eigenfunctions as in~\eqref{eq:EFHelmHoltzBall}.
For $\ell \geq 1$, we proceed as before. For $\ell=0$, we use additionally $K_0(x) = -\log\frac{x}{2}+O(1)$ as $x \searrow 0$. 
This yields 
\[
\lim_{\Lambda\searrow 0} \mu_{(0)} \left( B_\rho^\ext, \Lambda \right)=0.
\]
Furthermore, with $b_0 = -\frac{1}{2\log\frac{\Lambda}{2}}$, we obtain
\[
\lim_{\Lambda\searrow 0} u_{0}(r, \theta) = 1.
\]
Therefore, as in the case of higher-dimensional balls, the  eigenvalues and eigenfunctions of \eqref{eq:PDEGrebChai} in the exterior of a disk converge to those of \eqref{eq:PDEext}.

We emphasise that in these examples the convergence of eigenfunctions is established only pointwise; in \S\ref{sec:equivhelmn2} we provide a stronger statement.
\end{example}

\subsection{Exterior problem for planar domains via conformal mappings}\label{subsec:conformal}

In this section, we present a new formulation of the Steklov eigenvalue problem in exterior domains in two dimensions, based on conformal mappings.  

Let $\Omega \subset \mathbb{R}^2$ be a  bounded  open set with Lipschitz boundary and with connected $\Omegaext$. 
We assume that the origin $0 \in \Omega$, which can always be achieved by a change of variables, and we identify $\mathbb{R}^2$ with $\mathbb{C}$ by interpreting $x=\left(x_1,x_2\right)\in\mathbb{R}^2$ as   $z=x_1+\ir x_2\in\mathbb{C}$.

Let
\begin{equation}\label{eq:phi}
\phi: \mathbb{C}\setminus \{ 0 \}  \to \mathbb{C}\setminus \{ 0 \},
\qquad \phi(z) = \frac{1}{z},
\end{equation}
and set 
\[
\Omega^* := \phi\left(\Omega^\ext\right)\cup\{0\}.
\]
We remark that $\phi$ is a conformal diffeomorphism between $\Omega^\ext$ and $\Omega^* \setminus \{0\}$. The connectedness of $\Omega^{\mathrm{ext}}$ implies that $\Omega^*$ is connected; furthermore, if $\Omega$ is connected, then $\Omega^*$ is simply connected. 
See Figure \ref{fig:reg2D} for an illustration. 

\begin{figure}[htb]
\centering
\includegraphics{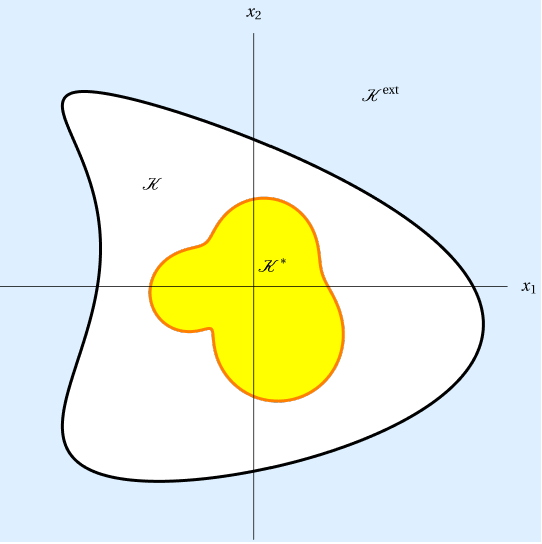}
\caption{The kite domain $\mathcal{K}$ from \S\ref{sec:numerics}, its exterior $\mathcal{K}^\ext$, and the bounded domain $\mathcal{K}^*=\phi\left(\mathcal{K}^\ext\right)\cup \{ 0 \}$.}\label{fig:reg2D}
\end{figure}

As in \cite{girouard2017spectral}, the key idea is to pull back the eigenvalue problem on $\Omegaext$ to the bounded domain $\Omega^*$. This yields the transformed eigenvalue problem 
\begin{equation*}\label{eq:problemtransform}
\begin{cases}
\Delta w = 0 &\text{ in } \Omega^*,\\
\partial_\nu w =  \xi |\phi '(z)| w&\text{ on } \partial\Omega^*,
\end{cases}\tag{$\mathrm{CT}$}
\end{equation*}
where $\xi$ is the new spectral parameter, and  $\left|\phi'(z)\right|=\frac{1}{|z|^2}.$

As always, we understand \eqref{eq:problemtransform} in the weak sense. It may seem natural to formulate~\eqref{eq:problemtransform} on $\Omega^* \setminus \{0\}$, but any bounded harmonic function $u : \Omegaext \to \mathbb{R}$ admits a limit at infinity. Consequently, the transformed function $w = u \circ \phi$ extends smoothly to all of $\Omega^*$. Also, we note that inverting with respect to a point in $\Omega$ other than the origin, does not affect the eigenvalues.

Since $\Omega^*$ is a bounded Lipschitz domain, the eigenvalues of~\eqref{eq:problemtransform} admit the characterisation 
\begin{equation}\label{eq:sigma*}
\xi_k\left(\Omegaext\right) = \inf_{\substack{M \subset H^1(\Omega^*) \\ \dim(M)=k}} \sup_{0 \neq w \in M} \frac{\int_{\Omega^*} | \nabla w |^2 \, \mathrm{d}z}{\int_{\partial\Omega^*} |\phi'|\cdot |  w |^2 \, \mathrm{d}S}.
\end{equation}
We will show later  that~\eqref{eq:sigma*} allows us to adapt classical arguments by Weinstock to derive isoperimetric-type inequalities for the first non-trivial Steklov eigenvalue. Prior to that, we will establish, in \S\ref{section:twod}, the equivalence between the formulation~\eqref{eq:problemtransform}   and the approaches described in \S\ref{subsec:truncatedomains} and \S\ref{subsec:Greben}.

\begin{remark}
This conformal mapping approach is specific to two dimensions. While the Kelvin transform preserves harmonicity in $\mathbb{R}^n$ for $n \geq 3$, the boundary condition becomes significantly more complicated under this transformation.
\end{remark}

\section{Boundary layer potentials}\label{sec:layer}

Numerical methods for solving boundary value problems often involve reformulating the boundary value problem on $\Omega$ as an integral equation over $\partial\Omega$, using single and double layer potentials. This approach is also applicable to boundary value problems on  $\Omega^\ext$, provided that the solutions have a suitable behaviour at infinity, see \cite[pp. 234--239]{mclean2000strongly}, \cite[Section 7.5]{steinbach2007numerical} and also~\cite{ChandlerWilde2012}. 

Recently in \cite{KNP25}, the Dirichlet-to-Neumann map in exterior domains has been studied using layer potentials. Building on these results, we formulate the exterior Steklov problem using layer potentials. While this is closely related to \cite{KNP25}, we also represent the limit at infinity via layer potentials, which makes the resulting formulation especially suitable for the numerical computation of eigenvalues. For completeness, we briefly recall known results below.

In this section, differential equations and harmonic functions are again understood in the weak sense, unless stated otherwise.

\subsection{Bounded domains}\label{subsubsec:extlayerbdd}
On bounded domains the Steklov eigenvalue problem on $\Omega$ can be reformulated in terms of single and double layer potentials; see, for instance, \cite[p.~244--245]{mclean2000strongly}, \cite[Section~6.6.3]{steinbach2007numerical}. In what follows, we express the exterior Steklov problem using the same boundary layer potentials, defined below.

\begin{definition}
Suppose $\Omega \subset \R^n$, $n \geq 2$, is a bounded open set with a Lipschitz boundary $\partial\Omega$. Let $\Phi$ be the fundamental solution for the Laplace equation,
\[
\Phi(x,y) := \begin{cases}  
- \frac{\log(|x-y|)}{2\pi} \quad&\text{ if } n=2,  \\
\frac{1}{(n-2) |\mathbb{S}^{n-1}| |x-y|^{n-2}} \quad&\text{ if } n>2, 
\end{cases} 
\qquad \text{ for } x,y \in \mathbb{R}^n \text{ and } x \neq y,
\]
where $|\mathbb{S}^{n-1}| = \frac{2\pi^{n/2}}{\Gamma(n/2)}$ is the measure of the unit sphere in $\R^n$. 

We define the following layer potentials. Note that, compared to \cite{steinbach2007numerical}, we might use a different sign convention.
\begin{itemize}
\item\textbf{The single layer potential}   $\operatorname{V}: H^{-\frac{1}{2}}(\partial\Omega) \to H^{\frac{1}{2}}(\partial\Omega)$, where $\operatorname{V}\eta: \partial\Omega \to \mathbb{R}$, is given by
\[
(\operatorname{V} \eta)(x) := \int_{\partial\Omega} \Phi(x,y) \eta(y) \, \mathrm{d}S_y.
\]
\item\textbf{The double layer potential} $\operatorname{K}: H^{\frac{1}{2}}(\partial\Omega) \to H^{\frac{1}{2}}(\partial\Omega)$, where $\operatorname{K}\eta: \partial\Omega \to \mathbb{R}$, is given by
\[
(\operatorname{K} \eta) (x):= \int_{\partial\Omega} (\partial_{\nu_y} \Phi)(x,y) \eta(y) \, \mathrm{d}S_y,
\]
and $\nu_y$ is the unit normal vector oriented outwards the domain $\Omega$ and
\[ 
\partial_{\nu_y} \Phi(x,y) = \mydotp{ \nu_y , \nabla_y \Phi(x,y) } = \frac{ \mydotp{ x-y,  \nu_y }}{|\mathbb{S}^{n-1}| |x-y|^n}.
\]
\item\textbf{The adjoint double layer potential} $\operatorname{K}': H^{-\frac{1}{2}}(\partial\Omega) \to H^{-\frac{1}{2}}(\partial\Omega)$, where $\operatorname{K}'\eta: \partial\Omega \to \mathbb{R}$, is given by 
\[
(\operatorname{K}' \eta) (x):= \int_{\partial\Omega} (\partial_{\nu_x} \Phi)(x,y) \eta(y) \, \mathrm{d}S_y.
\]
\item\textbf{The hypersingular boundary integral operator} $\operatorname{D}: H^{\frac{1}{2}}(\partial\Omega) \to H^{-\frac{1}{2}}(\partial\Omega)$, where $\operatorname{D}\eta: \partial\Omega \to \mathbb{R}$, is given by
\[
(\operatorname{D} \eta) (x):= -\partial_{\nu_x}\int_{\partial\Omega} ( \partial_{\nu_y} \Phi)(x,y) \eta(y) \, \mathrm{d}S_y.
\]
\end{itemize}
Moving the normal derivative inside the integral in $(\operatorname{D} \eta)$ may lead to a non-integrable singularity. Also note that $\operatorname{V}$ and $\operatorname{D}$ are self-adjoint, and $\operatorname{K}$ and $\operatorname{K}'$ are formally adjoint on $L^2(\partial \Omega)$. 
\end{definition}
 
For connected $\Omega$, the classical (interior) Dirichlet-to-Neumann operator  $\mathcal{D}: H^{\frac{1}{2}}(\partial\Omega) \to H^{-\frac{1}{2}}(\partial\Omega)$  is considered in \cite[Section 6.6.3]{steinbach2007numerical} (note that  $\sigma(x)$ in \cite{steinbach2007numerical} equals  $\frac{1}{2}$ almost everywhere since $\partial\Omega$ is Lipschitz). Provided that the single layer operator $\operatorname{V}$ is invertible, $\mathcal{D}$ admits the representation
\[
\mathcal{D}= \operatorname{V}^{-1} \left(  \frac{1}{2} \operatorname{I} - K \right) = \operatorname{D} + \left( \frac{1}{2} \operatorname{I} - \operatorname{K}' \right) \operatorname{V}^{-1} \left( \frac{1}{2} \operatorname{I} - \operatorname{K} \right),
\]
where the second representation has the advantage of being symmetric. Note that the double and adjoint double layer operators defined in \cite{steinbach2007numerical} equal $-\operatorname{K}$ and $-\operatorname{K}'$, respectively, because of the orientation of $\nu$.
Consider a pair $(u,\tau) \in H^{\frac{1}{2}}(\partial\Omega)\times \mathbb{R}$ satisfying $\mathcal{D} u = \tau u$, or, equivalently
\begin{equation} \label{eq:Steklov_integral_interior}
\frac12 u - \operatorname{K}u
= \tau \operatorname{V}u.
\end{equation}
Then the harmonic extension $U$ of $u$ into $\Omega$, given by
\[
\begin{split}
U: \Omega \to \mathbb{R}, \quad U(x)&:= 
\tau  \int_{\partial\Omega} \Phi(x,y)   u(y) \, \mathrm{d}S_y + \int_{\partial\Omega} (\partial_{\nu_y} \Phi)(x,y) u(y) \, \mathrm{d}S_y
\end{split}
\]
is an eigenfunction of~\eqref{eq:PDEint} corresponding to the eigenvalue $\tau$.

By Green's third identity, for $x \in \partial \Omega$,
\begin{equation}  \label{eq:Green_3rd}
\int\limits_{\partial\Omega} (\partial_{\nu_y} \Phi)(x,y) \, \mathrm{d}S_y = \frac12 \,,
\end{equation}
so that $u \equiv 1$ solves~\eqref{eq:Steklov_integral_interior} with $\tau = 0$, as expected. 

Finally, let us recall that the single layer operator is always invertible in dimensions $n\ge3$, whereas for $n=2$ we have to assume a suitable scaling of the domain $\Omega$ \cite[Section 6.6]{steinbach2007numerical}. If $n=2$ and  the logarithmic capacity of $\Omega$ is equal to one, the single layer potential has a nontrivial kernel and one cannot use~\eqref{eq:Steklov_integral_interior} to define the eigenvalues of the Dirichlet-to-Neumann operator. 

\subsection{Exterior domains in dimension two}\label{subsec:dimtwo}
While the boundary integral equation~\eqref{eq:Steklov_integral_interior} remains valid for characterising those eigenfunctions whose harmonic extensions vanish at infinity (up to switching the sign), see \S\ref{subsubsec:extlayer}, the situation differs in dimension two: as noted in \S\ref{examp:balln2}, exterior Steklov eigenfunctions in two dimensions are only required to remain bounded and do not need to decay. The Dirichlet-to-Neumann operator in exterior domains, including dimension two, together with its Fredholm properties, has been analysed in \cite[Section 3]{KNP25}. In contrast to \cite{KNP25}, our formulation also incorporates the limit at infinity via layer potentials.

Let us consider a harmonic function $u \in H^1_{\mathrm{loc}}\left(\Omegaext\right)$ which satisfies the far field boundary condition
\begin{equation}\label{eq:farfield}
\left|  u(x)-u_\infty \right| = O \left( \frac{1}{|x|} \right) \qquad \text{ as } |x| \to \infty
\end{equation}
for a given $u_\infty \in \mathbb{R}$. According to \cite[Section 7.5]{steinbach2007numerical} (note that $\gamma_1^\mathrm{ext}$ in \cite{steinbach2007numerical} is $-\partial_{\nu} $ with our convention for $\nu$), we obtain, for $x \in \Omega^\ext$,
\[
u(x) = u_\infty + \int_{\partial\Omega} \Phi(x,y) (\partial_{\nu_y} u)(y) \, \mathrm{d}S_y - \int_{\partial\Omega} \left( \partial_{\nu_y} \Phi \right) (x,y)   u(y) \,  \mathrm{d}S_y.
\]
Similar to the interior case, applying the exterior trace operator yields as in \cite[p. 182]{steinbach2007numerical},
\begin{equation}\label{eq:ext_2D_2.5}
u(x) =u_\infty + (\operatorname{V} (\partial_\nu u))(x) + \frac{1}{2} u(x) - (\operatorname{K} u)(x) \quad \text{ for $x \in \partial\Omega$}.
\end{equation}
We aim to express the constant $u_\infty$ in terms of single and double layer operators. To this end, we write $u_\infty$ in terms of the harmonic measure density. Let $G_{\Dir}$ be the Dirichlet Green's function ``with pole at infinity'' (see, e.g. \cite[p. 41 and p. 74]{garnett2005harmonic}) for the exterior
problem. This means, $G_{\Dir}$ is the solution of
\[
\begin{cases}
\Delta_y G_{\Dir}(y) = 0 &\text{ in } \Omega^\ext, \\
G_{\Dir}(y) = 0 &\text{ on } \partial\Omega, \\
G_{\Dir}(y) - \frac{\log (|y|)}{2\pi}  + \frac{\log(\operatorname{Cap}_\mathrm{log}(\Omega))}{2 \pi} =  O\left( \frac{1}{|y|} \right) &\text{ as } |y|\to \infty,
\end{cases}
\]
where  $ \operatorname{Cap}_\mathrm{log}(\Omega) $ is the logarithmic capacity of $ \Omega $. While \cite{garnett2005harmonic} only states $ {o}(1) $ decay, any harmonic function vanishing at infinity behaves as $ O(1/|y|) $, cf.\ Lemma \ref{lemma:usefulbasic}.

Define the harmonic measure density $\omega: \partial\Omega \to \mathbb{R}$, by 
\[
\omega(y) = - \partial_{\nu_y} G_{\Dir}(y).
\]
Then, the constant $u_\infty$ can be written as
\[
u_\infty = \int_{\partial\Omega} u(y) \omega(y) \, \mathrm{d}S_y.
\] 
Thus, equation~\eqref{eq:ext_2D_2.5} becomes	
\begin{equation}\label{eq:ext_2D_3}
\frac{u(x)}{2} +(\operatorname{K} u)(x) - \int_{\partial\Omega} u(y) \omega(y)  \, \mathrm{d}S_y =  (\operatorname{V} (\partial_\nu u))(x).
\end{equation}
The main disadvantage of~\eqref{eq:ext_2D_3} is its dependence on the generally unknown harmonic measure density $\omega(y)$.  In the following, we express the Dirichlet Green's function $G_{\Dir}$ in terms of a single layer potential to overcome this issue.

As  in \cite[p. 25]{ahlfors2010conformal} or \cite[p. 41]{garnett2005harmonic}, the Dirichlet Green's function can be expressed as 
\[
G_{\Dir}(x) = - \Phi(x,0) - g_0 + g(x),
\]
where $g_0 = \frac{\log(\operatorname{Cap}_\mathrm{log}(\Omega))}{2 \pi}$ and $g$ is the harmonic function satisfying
\[
\Delta g = 0 \quad \textrm{in}~\Omega^\ext, \qquad  g = \Phi(\cdot,0) + g_0 \quad \textrm{on}~\partial\Omega,  
\qquad g(x) = O \left( \frac{1}{|x|} \right) \quad \text{as}~|x|\to { \infty}.
\] 
Then,
\[
\omega(x) = - \partial_{\nu_x} G_{\Dir}(x) = ( \partial_{\nu_x} \Phi ) (x,0) - (\partial_{\nu_x} g)(x) = - \frac{\mydotp{ x,  \nu(x) } }{2\pi |x|^2} -( \partial_{\nu_x} g)(x).
\]

Let $u \in H^\frac{1}{2}(\partial \Omega)$,  and let $U =\mathcal{H}^\ext u\in H^1_\mathrm{loc}\left(\Omegaext\right)$ be the harmonic extension of $u$ satisfying  \eqref{eq:inftyr}.  Green's identities
on $\Omegaext_R$ yield 
\[
\begin{split}
\int_{\partial\Omega} u(y) (\partial_{\nu_y} g)(y) \, \mathrm{d}S_y  &=   \int_{\partial\Omega} g(y) (\partial_{\nu_y} U)(y) \, \mathrm{d}S_y\\
&+ \int_{\partial B_R} \left(g(y) (\partial_{\nu_y} U)(y) -  U(y) (\partial_{\nu_y} g)(y)\right) \, \mathrm{d}S_y
\end{split}
\]
for any $R>R_0(\Omega)$. According to Lemma \ref{lemma:usefulbasic}, the behaviour of $g$ and $U$ at infinity implies that the last integral vanishes as $R \to \infty$ (note that $\partial_\nu U = \partial_r U $ on $\partial B_R$). As a consequence,
\[
u_\infty = \int_{\partial\Omega} u(y) \omega(y) \, \mathrm{d}S_y
= - \int_{\partial\Omega} \left(u(y)  \frac{\mydotp{ y,  \nu_y }}{2\pi |y|^2} + (\partial_{\nu_y} U)(y) g(y)\right) \, \mathrm{d}S_y.
\]
Since $g = \Phi(\cdot,0) + g_0$ on $\partial\Omega$, we have
\begin{equation}\label{eq:u_0}
u_\infty =  - \int_{\partial\Omega} u(y)\frac{\mydotp{ y,  \nu_y }}{2\pi |y|^2} + (\partial_{\nu_y} U)(y) \left( \frac{-\log(|y|)}{2 \pi} + \frac{\log(\operatorname{Cap}_\mathrm{log}(\Omega))}{2 \pi}\right)  \, \mathrm{d}S_y.
\end{equation}
Hence,~\eqref{eq:ext_2D_3} becomes
\[
\begin{split}
\frac{u(x)}{2} &+\int_{\partial\Omega} u(y)   \frac{\mydotp{ y,  \nu_y }}{2\pi |y|^2} \, \mathrm{d}S_y + (\operatorname{K} u)(x)\\ 
&= \int_{\partial\Omega}  (\partial_{\nu_y} U)(y)     \frac{\log(|y|)-\log(\operatorname{Cap}_\mathrm{log}(\Omega))}{2 \pi}   \, \mathrm{d}S_y + (\operatorname{V} (\partial_\nu U))(x).
\end{split}
\]
To derive a representation of the Dirichlet-to-Neumann map, we define the operator $\operatorname{V}_0: H^{-\frac{1}{2}}(\partial\Omega) \to H^{\frac{1}{2}}(\partial\Omega)$, where $\operatorname{V}_0\eta: \partial\Omega \to \mathbb{R}$ is given by 
\[
(\operatorname{V}_0 \eta)(x) := \int_{\partial\Omega} \left( \frac{\log(|y|)}{2 \pi} - \frac{\log(\operatorname{Cap}_\mathrm{log}(\Omega))}{2 \pi} + \Phi(x,y) \right) \eta(y) \, \mathrm{d}S_y.
\]
Additionally, we define the functional
$\operatorname{K}_0: H^{\frac{1}{2}}(\partial\Omega) \to \mathbb{R}$, as 
\[
\operatorname{K}_0 (\eta) :=\int_{\partial\Omega} \eta(y)   \frac{\mydotp{ -y,  \nu_y }}{2\pi |y|^2} \, \mathrm{d}S_y = (\operatorname{K} \eta)(0).
\]
If $\operatorname{V}_0$ is invertible (see Lemma~\ref{lemma:invsinglelayer}), then the Dirichlet-to-Neumann map can be expressed as
\[
\mathcal{D}^\ext:= \operatorname{V}_0^{-1} \left(  \frac{1}{2} \operatorname{I}- \operatorname{K}_0  + \operatorname{K} \right): H^{\frac{1}{2}}(\partial\Omega) \to H^{-\frac{1}{2}}(\partial\Omega).
\]
Accordingly,  $\mathcal{D}^\ext u = \tau u$  is equivalent to
\begin{equation*}\label{eq:ext_2D_2}\tag{$\mathrm{LP}_2$}
\frac{u}{2} -\operatorname{K}_0 (u) + (\operatorname{K} u) =  \tau (\operatorname{V}_0  u), \qquad u \in H^\frac{1}{2}(\partial \Omega).
\end{equation*}
Writing \eqref{eq:ext_2D_2} as an integral equation, we get
\begin{equation}\label{eq:ext_2D_4}
\begin{split}
\frac{u(x)}{2}  &+ \int_{\partial\Omega} \left(\frac{\mydotp{ x-y, \nu_y }}{2\pi|x-y|^2} + \frac{\mydotp{ y , \nu_y }}{2\pi |y|^2}\right) u(y) \, \mathrm{d}S_y \\
&= \, \tau \int_{\partial\Omega} \frac{-\log(|x-y|)+\log (|y|)-\log (\operatorname{Cap}_\mathrm{log}(\Omega))}{2\pi}  u(y) \, \mathrm{d}S_y.
\end{split}
\end{equation}
As a consistency check, consider the case $u=1$. Using the third Green's identity
\eqref{eq:Green_3rd} (with reversed orientation) yields $\operatorname{K} u=\frac{1}{2}$. Since the harmonic extension of $u = 1$ is constant, we have $u_\infty = 1$, which makes the left-hand side vanish. One therefore retrieves $\tau_1 = 0$, as it should be. 

Finally, it remains to verify the invertibility of $\operatorname{V}_0$. The assumption made on the diameter of $\Omega$ can be satisfied by appropriately rescaling the domain.

\begin{lemma}\label{lemma:invsinglelayer}
Suppose $\Omega \subset \mathbb{R}^2$ is a bounded open set with Lipschitz boundary, with connected $\Omegaext$ and $\operatorname{diam}(\Omega) <1$. Then, $\operatorname{V}_0: H^{-\frac{1}{2}}(\partial\Omega) \to H^{\frac{1}{2}}(\partial\Omega)$ is invertible. 
\end{lemma}

\begin{proof}
In \cite[Theorem 6.23]{steinbach2007numerical}, it is shown that $\operatorname{V}:H^{-\frac{1}{2}}(\partial\Omega) \to H^{\frac{1}{2}}(\partial\Omega)$ is invertible if $\operatorname{diam}(\Omega) <1$. Since we build up on this proof, we briefly recall the relevant constructions. Define the subspace
\[
H_*^{-\frac{1}{2}}(\partial\Omega) := \left\{ \eta \in H^{-\frac{1}{2}}(\partial\Omega) :  \myscal{\eta, 1 }_{\partial\Omega}=0 \right\},
\]
where $\myscal{\cdot, \cdot }_{\partial\Omega}$ means the duality pairing $H^{-\frac{1}{2}}(\partial\Omega) \times H^{\frac{1}{2}}(\partial\Omega) \to \mathbb{R}$. Let $w_\mathrm{eq} \in H^{-\frac{1}{2}}(\partial\Omega)$ be the unique function satisfying
\[
\operatorname{V}(w_\mathrm{eq})(x)=\frac{-\log(\operatorname{Cap}_\mathrm{log}(\Omega))}{2 \pi} \quad \text{ for $x \in \partial\Omega$ }  \quad \text{and} \quad \myscal{ w_\mathrm{eq},1 }_{\partial\Omega} = 1.
\]
Any $w \in  H^{-\frac{1}{2}}(\partial\Omega)$ admits  a unique decomposition
\[
w = w_* + \alpha w_\mathrm{eq}, \qquad \text{with} \quad w_* \in H_*^{-\frac{1}{2}}(\partial\Omega), \quad \alpha=\myscal{ w,1 }_{\partial\Omega} \in \mathbb{R}.
\]
The image of $H_*^{-\frac{1}{2}}(\partial\Omega)$  under $\operatorname{V}$, is given by
\[
\left\{ v \in H^{\frac{1}{2}}(\partial\Omega): \myscal{ v, w_\mathrm{eq} }_{\partial\Omega} = 0 \right\}.
\]
The preliminaries from \cite{steinbach2007numerical} being now in place, we turn to the proof for $\operatorname{V}_0$, which can be written as $\operatorname{V}_0(\eta) = \operatorname{V}(\eta) + \operatorname{W}(\eta)$, where $\operatorname{W}: H^{-\frac{1}{2}}(\partial\Omega) \to \mathbb{R}$, is defined by
\[
\begin{split}
\operatorname{W}(\eta) &= \frac{1}{2 \pi} \int_{\partial\Omega} \left( \log(|y|) - \log(\operatorname{Cap}_\mathrm{log}(\Omega))  \right) \eta(y) \, \mathrm{d}S_y\\ 
&= -(\operatorname{V}(\eta)(0) - \frac{\log(\operatorname{Cap}_\mathrm{log}(\Omega))}{2 \pi} \myscal{\eta , 1}_{\partial \Omega}.
\end{split}
\]
To show that $\operatorname{V}_0$ is invertible, we prove  its kernel is trivial. Assume $\operatorname{V}_0(w)=0$. For $w = w_* + \alpha w_\mathrm{eq}$,
\[
\begin{split}
\operatorname{V}_0(w) &= \operatorname{V}(w) + \operatorname{W}(w) = \operatorname{V}(w_*) + \alpha \operatorname{V}(w_\mathrm{eq})+\operatorname{W}(w)\\ 
&= \operatorname{V}(w_*) + \alpha \frac{-\log(\operatorname{Cap}_\mathrm{log}(\Omega))}{2 \pi}+\operatorname{W}(w).
\end{split}
\]
Since $\operatorname{W}(w)$ is a constant, $\operatorname{V}(w_*)$ must be constant. However, $ \operatorname{V}(w_*)$ is orthogonal to $w_\mathrm{eq}$ and since $\myscal{ w_\mathrm{eq},1 }_{L^2(\partial\Omega)} = 1$, it follows that $\operatorname{V}(w_*)=0$. By the invertibility of $\, \operatorname{V}$, we conclude $w_* = 0$. Moreover, 
\[
\begin{split}
\operatorname{W}(w_\mathrm{eq}) &= -\operatorname{V}(w_\mathrm{eq})(0) - \frac{\log(\operatorname{Cap}_\mathrm{log}(\Omega))}{2 \pi} \myscal{w_\mathrm{eq} , 1}_{\partial \Omega}\\ 
&= - \frac{-\log(\operatorname{Cap}_\mathrm{log}(\Omega))}{2 \pi} - \frac{\log(\operatorname{Cap}_\mathrm{log}(\Omega))}{2 \pi} = 0.
\end{split}
\]
Thus, if $\operatorname{V}_0( \alpha w_\mathrm{eq}) = 0$, then $\operatorname{V}( \alpha w_\mathrm{eq}) = 0$ and therefore $\alpha = 0$ because $\operatorname{V}$ is invertible. Thus, $\operatorname{V}_0$ is injective.

On the other hand, for any $v \in H^\frac{1}{2}(\partial \Omega)$ there exists a $w \in  H^{-\frac{1}{2}}(\partial \Omega)$ with $\operatorname{V}(w) = v$ because $\operatorname{V}$ is invertible. Then, since $\operatorname{W}\left( w \right) \in \mathbb{R}$ and $\operatorname{W}\left( w_\mathrm{eq} \right) = 0$,
\[
\begin{split}
\operatorname{V}_0\left( w + \frac{2 \pi \operatorname{W}(w)}{\log(\operatorname{Cap}_\mathrm{log}(\Omega))} w_\mathrm{eq} \right) &= \operatorname{V} ( w )+ \frac{2 \pi \operatorname{W}(w)}{\log(\operatorname{Cap}_\mathrm{log}(\Omega))} \operatorname{V}\left( w_\mathrm{eq} \right)+\operatorname{W} ( w ) \\
&= v+ \frac{2 \pi \operatorname{W}(w)}{\log(\operatorname{Cap}_\mathrm{log}(\Omega))} \frac{-\log(\operatorname{Cap}_\mathrm{log}(\Omega))}{2 \pi}+\operatorname{W} ( w ) = v.
\end{split}
\]
Thus, $\operatorname{V}_0$ is surjective.
\end{proof}

To summarise, we present the following

\begin{theorem}\label{theo:layer2d}
Suppose $\Omega \subset \R^2$ is a bounded open set with Lipschitz boundary and with connected $\Omegaext$, and $\operatorname{diam}(\Omega) < 1$.  The boundary integral equation~\eqref{eq:ext_2D_2} has a sequence of eigenvalues
\[
0 = \tau_1 \leq \tau_2 \leq \ldots
\]
accumulating to $+\infty$, with associated eigenfunctions $u_k \in H^{\frac{1}{2}}(\partial\Omega)$. For any eigenfunction $u_k $, we define $U_k: \Omega^\ext \to \mathbb{R}$ by
\begin{equation}\label{eq:Ukx}
 U_k(x) := (u_k)_\infty + \tau_k \int_{\partial\Omega} \Phi(x,y)  u_k(y) \, \mathrm{d}S_y - \int_{\partial\Omega} \left( \partial_{\nu_y} \Phi\right) (x,y)   u_k(y) \,  \mathrm{d}S_y,
\end{equation}
where $(u_k)_\infty$ is the limit at infinity, cf.\ \eqref{eq:u_0}. Then, $U_k$ is a weak solution of
\[
\begin{cases}
\Delta U_k  = 0\qquad &\text{ in } \Omegaext, \\
\partial_\nu U_k  = \tau_k U_k\qquad&\text{ on } \partial\Omega,\\
\left|U_k(x)-(u_k)_\infty\right| = O( |x|^{-1}) \qquad&\text{ as } |x| \to \infty.
\end{cases}
\]
Additionally, if $(U_k, \tau_k)$ is a solution of the above boundary-value problem, then $\tau_k$ and $U_k |_{\partial \Omega}$  satisfy \eqref{eq:ext_2D_2}.
\end{theorem}

\begin{proof}
The main steps of the proof are given above. However, the asymptotic behaviour of $U_k$ may not be immediately apparent. From the definition of $\Phi$, it is clear that the last integral in \eqref{eq:Ukx} decays. For the first eigenfunction, which is constant, the first integral vanishes. In the case of functions orthogonal to a constant, the decay of the first integral is shown in \cite[Lemma 6.21]{steinbach2007numerical}.
\end{proof}

\begin{remark}\label{remark:layerpot}
When restricting to functions with $\myscal{ u,1 }_{L^2(\partial \Omega)}=0$, the constant $\operatorname{Cap}_\mathrm{log}(\Omega)$ in~\eqref{eq:ext_2D_4}  can be replaced by any other constant without affecting the value of the integral. This could be useful for numerical calculations, as determining the logarithmic capacity, despite the availability of various numerical techniques, introduces an additional computational step. 
\end{remark}

\subsection{Exterior domains in dimensions $n\ge3$}

\subsubsection{Decay at infinity}\label{subsubsec:extlayer}

As  in \cite{KNP25} (see also \cite{Sayas,mclean2000strongly,steinbach2007numerical}), the boundary integral equation~\eqref{eq:Steklov_integral_interior} becomes
\begin{equation*}\label{eq:Steklov_integral}\tag{$\mathrm{LP}_n$}
\frac12 u + \operatorname{K}u
= \tau \operatorname{V}u, \qquad u \in H^\frac{1}{2}(\partial \Omega)
\end{equation*}
for the exterior Steklov problem in $\Omega^{\mathrm{ext}}$ when $n \geq 3$ and we assume decay at infinity. 
Once an eigenfunction $u$ is found on $\partial\Omega$ with eigenvalue $\tau$ , it
can be harmonically extended into $\Omega^\ext$ as $U=\mathcal{H}u$, see \S\ref{sec:statement}. 
As before, this harmonic extension can also be written as
\[
U(x) =  \tau \int_{\partial\Omega} \Phi(x,y)  u(y) \, \mathrm{d}S_y - \int_{\partial\Omega} \left( \partial_{\nu_y} \Phi \right)(x,y)  u(y) \,  \mathrm{d}S_y.
\]
Since $\Phi(x,y)=\frac{1}{(n-2)|\mathbb{S}^{n-1}| |x-y|^{n-2}}$, it is easy to see that $U(x)$  satisfies the decay condition \eqref{eq:inftyr}. 

In summary, we have the following 

\begin{theorem}\label{thm:equivGlobalLayerdim3}
Suppose $\Omega \subset \R^n$, $n \geq3$, is a bounded open set with Lipschitz boundary and  with connected $\Omegaext$.  The boundary integral equation~\eqref{eq:Steklov_integral} has a sequence of eigenvalues
\[
0 < \tau_1 \leq \tau_2 \leq \ldots,
\]
accumulating to $+\infty$.
For any eigenfunction $u_k $, we define $U_k \in E^1\left(\Omegaext\right)$ by
\[
 U_k(x) := \tau_k \int_{\partial\Omega} \Phi(x,y)  u_k(y) \, \mathrm{d}S_y - \int_{\partial\Omega} \left( \partial_{\nu_y} \Phi \right) (x,y)   u_k(y) \,  \mathrm{d}S_y.
\]
Then, $U_k$ is a weak solution of
\[
\begin{cases}
\Delta U_k  = 0 \, &\text{ in } \Omegaext, \\
\partial_\nu U_k  = \tau_k U_k \, &\text{ on } \partial\Omega,\\
U_k \text{ satisfies \eqref{eq:inftyr}}.
\end{cases}
\]
\end{theorem}

As an immediate consequence, we have the equivalence to the finite energy approach.

\begin{corollary}\label{cor:FELP}
For any  bounded  open set $\Omega \subset \mathbb{R}^n$, $n\geq 3$, with Lipschitz boundary and with connected $\Omegaext$, and  any $k \in \mathbb{N}$,
we have
\[
\lambda_k\left(A^\Dir\right) = \tau_k(\Omegaext).
\]
Moreover, if $u_k$ is an eigenfunction associated to $\lambda_k\left(A^\Dir\right)$, then $\left( u_k,\lambda_k\left(A^\Dir\right) \right) $ is a solution of \eqref{eq:Steklov_integral} and  vice versa.
\end{corollary}

\subsubsection{Vanishing flow at infinity}\label{subsec:vanishflow}

In \cite{henrici1970sloshing}, Henrici, Troesch and Wuytack investigate the exterior Steklov problem in the upper half space  with circular or strip-like aperture in the special case $n = 3$.  Instead of a decay property, they impose the condition, which can be generalised to any $n \geq 2$ as
\begin{equation}\label{eq:alternbehav}
|x|^{n-1} | \nabla u(x)| \to 0 \qquad \text{as} \quad |x| \to \infty.
\end{equation}
To establish the well-posedness of the exterior Steklov problem under~\eqref{eq:alternbehav}, they consider an equivalent boundary integral equation. However, their approach relies on the specific geometry of the domain (a disc or a strip) and is not directly applicable to a general domain $\Omega$. By adapting the calculations from \S\ref{subsec:dimtwo} to higher dimensions (we only have to adjust Green's function), we obtain a representation for the corresponding Dirichlet-to-Neumann operator.

Let $\Omega \subset \mathbb{R}^n$, $n \geq 3$, be a bounded open set with Lipschitz boundary and $u \in H^1_{\text{loc}}(\Omegaext)$ a harmonic function satisfying \eqref{eq:alternbehav}.  First, note that \eqref{eq:alternbehav} implies that $|\nabla u| \in L^2\left(\Omegaext\right)$. Hence, the average of $u$ as defined in \eqref{eq:averageu} and denoted by $\langle u \rangle$, is finite. 
Define $w := u - \langle u \rangle$. Then, by definition,  $w \in W^\Dir\left(\Omegaext\right)$ (see  \S\ref{sec:DtNs}) and therefore, by \eqref{eq:WD=E1}, we have  $w \in E^1\left(\Omegaext\right)$. By Proposition \ref{prop:xiongdecay}, $w = O(|x|^{2-n})$ which shows that a higher-dimensional analogue of the far-field condition \eqref{eq:farfield} is satisfied with $u_\infty := \langle u \rangle$.

As before, $u_\infty$ can be expressed in terms of the Dirichlet Green's function with pole at infinity,
\[
\begin{cases}
\Delta_y G_{\Dir}(y) = 0 \qquad&\text{ in } \Omega^\ext, \\
G_{\Dir}(y) = 0 \qquad&\text{ on } \partial\Omega, \\
G_{\Dir}(y)   =  a_0 +  \frac{a_1}{|y|^{n-2}} + O\left( |y|^{1-n} \right)\qquad&\text{ as } |y|\to \infty.
\end{cases}
\]
The constants $a_0$ and $a_1$ are determined as follows. The function $P(y) :=  \frac{a_0-G_{\Dir}(y)}{a_0}$ (if $a_0 = 0$, then $G_{\Dir} \equiv 0$) is the capacitary potential of $\Omega$ and therefore, by definition of the Newtonian capacity,
\[
P(y) = \frac{\operatorname{Cap}(\Omega)}{(n-2)|\mathbb{S}^{n-1}| |y|^{n-2}} + o(|y|^{n-2}) \qquad \text{as } |y| \to \infty,
\]
where $\operatorname{Cap}(\Omega)$ is the Newtonian capacity. This implies $ a_1 = -a_0\frac{\operatorname{Cap}(\Omega)}{(n-2)|\mathbb{S}^{n-1}|}$. On the other hand, we normalise $G_{\Dir}$ by requiring
\[
\int_{\partial \Omega} (-\partial_{\nu_y} G_{\Dir} )(y) \, \mathrm{d}S_y = 1.
\]
Using Green's identities, this implies $1= a_1 (2-n) |\mathbb{S}^{n-1}|$. Hence, we choose 
\[
a_0 = \frac{1}{\operatorname{Cap} (\Omega)} \qquad \text{ and } \qquad a_1 = \frac{1}{|\mathbb{S}^{n-1}|(2-n)}.
\]
 Then, as before,
\[
u_\infty = \int_{\partial \Omega} u(y) \omega(y) \, \mathrm{d}S_y = \int_{\partial \Omega} u(y) \left( -\partial_{\nu_y} G_{\Dir}\right) (y)  \, \mathrm{d}S_y.
\]
Moreover, the Dirichlet Green's function can be written as 
\[
G_{\Dir}(y) = \frac{1}{\operatorname{Cap} (\Omega)} + \frac{1}{|\mathbb{S}^{n-1}|(2-n)|y|^{n-2}}   + g(y),
\]
where $g$ is the harmonic function satisfying
\[
\begin{cases}
\Delta g = 0 \qquad&\text{ in } \Omega^\ext, \\
g =  -\frac{1}{\operatorname{Cap} (\Omega)} - \frac{1}{|\mathbb{S}^{n-1}|(2-n)|y|^{n-2}} \qquad&\text{ on } \partial\Omega, \\
g(y) = O \left( |y|^{1-n} \right)\qquad&\text{ as } |y|\to \infty.
\end{cases}
\] 
If $u$ is harmonic, integration by parts yields 
\[
\begin{split}
u_\infty &=
\int_{\partial\Omega} u(y) ( - \partial_{\nu_y} G_{\Dir})(y) \, \mathrm{d}S_y\\	
&=  \int_{\partial\Omega} - u(y)  \frac{\langle y, \nu_y \rangle}{|\mathbb{S}^{n-1}| |y|^{n}} \, \mathrm{d}S_y - \int_{\partial\Omega} g(y) (\partial_{\nu_y} u)(y) \, \mathrm{d}S_y \\
&\quad- \int_{\partial B_R} \left(g(y) (\partial_{\nu_y} u)(y) -  u(y) (\partial_{\nu_y} g)(y)\right) \, \mathrm{d}S_y
\end{split}
\]
for any $R>R_0(\Omega)$. Again, the last integral vanishes because $u=O(1)$, $\partial_r u = O(r^{1-n})$, $g=O(r^{1-n})$ and $\partial_r g = O(r^{-n})$ (cf.\ the proof of Proposition \ref{pro:decay}).
Hence,~\eqref{eq:ext_2D_2.5} becomes
\[
\frac{u}{2} = \widehat{\operatorname{V}}_0\left(\partial_\nu u\right)  - 	\widehat{\operatorname{K}}_0 u \quad \text{ on $\partial\Omega$},
\]
where
\[ 
\left( \widehat{\operatorname{V}}_0 \eta \right)(x) := \left( {\operatorname{V}} \eta \right)(x) + \int_{\partial \Omega} \eta(y) \left( \frac{1}{\operatorname{Cap} (\Omega)} + \frac{1}{|\mathbb{S}^{n-1}|(2-n)|y|^{n-2}} \right) \, \mathrm{d}S_y, 
\]
and 
\[ 
\left( \widehat{\operatorname{K}}_0 \eta \right)(x) := \left( {\operatorname{K}} \eta \right)(x) + \int_{\partial \Omega} \eta(y) \frac{\langle y, \nu_y \rangle}{|\mathbb{S}^{n-1}| |y|^{n}}\, \mathrm{d}S_y.
\]
As in Lemma \ref{lemma:invsinglelayer}, we can show that $\widehat{\operatorname{V}}_0$ is invertible. The only difference is that the logarithmic Green function and the logarithmic equilibrium density $w_\mathrm{eq}$ are replaced by their Newtonian counterparts; see, for instance, \cite[formula~(6.35)]{steinbach2007numerical}.

\begin{lemma}
Suppose $\Omega \subset \mathbb{R}^n$, $n\geq 3$, is a bounded open set with Lipschitz boundary and with connected $\Omegaext$. Then, $\widehat{\operatorname{V}}_0: H^{-\frac{1}{2}}(\partial\Omega) \to H^{\frac{1}{2}}(\partial\Omega)$ is invertible. 
\end{lemma}

Thus, we have the following representation for the Dirichlet-to-Neumann operator
\[
\mathcal{D}:= \widehat{\operatorname{V}}_0^{-1} \left(  \frac{1}{2} \operatorname{I}+ \widehat{\operatorname{K}}_0 \right): H^{\frac{1}{2}}(\partial\Omega) \to H^{-\frac{1}{2}}(\partial\Omega).
\]

\begin{proposition}
Let $\Omega \subset \mathbb{R}^n$, $n\geq 3$ be a  bounded open set with Lipschitz boundary and with connected $\Omegaext$.
\begin{enumerate}[{\textrm (a)}]
\item Let $u \in H^{\frac{1}{2}}(\partial\Omega)$ , $\chi \in \mathbb{R}$ solve $\frac{u}{2}+ \widehat{\operatorname{K}}_0(u)   = \chi \widehat{\operatorname{V}}_0(u)$, 
then $U \in H^{1}_\mathrm{loc}\left(\Omegaext\right)$, defined as
\[
 U(x) := (u)_\infty + \chi \int_{\partial\Omega} \Phi(x,y)  u(y) \, \mathrm{d}S_y - \int_{\partial\Omega} \left( \partial_{\nu_y} \Phi \right) (x,y) u(y) \,  \mathrm{d}S_y
\]
solves
\begin{equation}\label{eq:henrici}\tag{$\mathrm{VF}$}
\begin{cases}
\Delta U  = 0\qquad&\text{in } \Omegaext,\\
\partial_\nu U  = \chi U\qquad&\text{on } \partial\Omega,\\
|x|^{n-1} | \nabla U(x)| \to 0\qquad&\text{as }  |x| \to \infty.
\end{cases}
\end{equation}
\item Let $\chi \in \R$ and $U \in H^{1}_\mathrm{loc}\left(\Omegaext\right)$ be a weak solution of~\eqref{eq:henrici}. Then, $u = U|_{\partial \Omega} \in H^{\frac{1}{2}}(\partial\Omega)$  solves $\frac{u}{2}+ \widehat{\operatorname{K}}_0(u)   = \chi \widehat{\operatorname{V}}_0(u)$.
\end{enumerate}
\end{proposition}

Since the first eigenvalue vanishes and the first eigenfunction is constant, this formulation is not equivalent to~\eqref{eq:PDEext},~\eqref{eq:PDEArendtElstDirichlet}, or~\eqref{eq:PDEGrebChai}, but to~\eqref{eq:PDEArendtElstNeumann}, see \S\ref{sec:VFNT}.

As in Remark \ref{remark:layerpot}, when restricting to functions with $\myscal{ u,1 }_{L^2(\partial \Omega)}=0$, the constant in $\hat{\operatorname{V}}_0$  can be replaced by any other constant without affecting the value of the integral.

\begin{remark}\label{prop:2D-equivalence}
In dimension $ n = 2$, the vanishing flow condition does not lead to a new problem, since the condition $|x|\cdot|\nabla u(x)| \to 0$ as $|x| \to \infty$  is equivalent to condition \eqref{eq:inftyr} for harmonic functions, which can easily be seen from Lemma \ref{lemma:usefulbasic}.
\end{remark}

\section{Equivalence of different approaches}\label{sec:equivalence}
In this section we compare the different formulations of the exterior Steklov problem introduced above and prove the corresponding equivalence statements. We first treat the case $n\geq 3$, where the finite-energy formulation is available and provides the natural common framework for the truncation, Helmholtz, and layer-potential approaches. We also prove the equivalence between the Neumann truncation and the vanishing flow approach in dimensions $n\ge 3$. The two-dimensional case requires a separate argument, because exterior harmonic functions have different behaviour at infinity. It is therefore handled separately at the end of the section using the conformal formulation.
\subsection{Equivalence in dimensions $n \geq 3$}\label{sec:ngeq3}

In this section, we prove that the approaches introduced in \S\ref{sec:survey} are equivalent for dimensions $n \geq 3$. The arguments are largely straightforward and rely on approximating the eigenfunctions of one formulation by functions from the corresponding alternative function space.

As before, all differential equations are understood in the weak sense.

\subsubsection{Finite-energy functions and Dirichlet truncation}\label{sec:AHvAtE}

One could use the resolvent convergence of $A_R^\Dir$, shown in \cite{arendt2015dirichlet}, to infer that the spectral problems~\eqref{eq:PDEext} and~\eqref{eq:PDEArendtElstDirichlet}  yield the same spectrum as $R \to \infty$. Note that the Neumann formulation~\eqref{eq:PDEArendtElstNeumann} leads to a different spectrum. We instead proceed directly via the variational characterisation. This approach yields a more direct and unified framework which is also applicable to ~\eqref{eq:PDEGrebChai}.

We use that smooth, compactly supported functions are dense in $E^1\left(\Omegaext\right)$.  Specifically, we use the following approximation property of functions in $E^1\left(\Omegaext\right)$. Let $\psi \in C_0^\infty(\mathbb{R}^n)$ be a cut-off function satisfying
\[
\psi(x)=1 \,\text{ for } |x|<1,\quad \psi(x)=0 \,\text{ for } |x|>2,\quad |\psi(x)|\leq 1  \text{ and }  |\nabla \psi(x)|\leq 2 \,\text{ for }  x \in \mathbb{R}^n.
\]
For each $R>0$, define the scaled function $\psi_R(x) := \psi\left( \frac{x}{R} \right)$. Since every $u \in E^1\left(\Omegaext\right)$ has zero average, i.e.\ $\langle u \rangle = 0$, the following lemma is an immediate consequence of~\cite[Proposition 2.2]{LuOu}; see also~\cite[Proof of Theorem 1]{bundrock2024optimizing}.

\begin{lemma}\label{lemma2}
Suppose $\Omega \subset \mathbb{R}^n$, $n \geq 3$, is a  bounded open set with Lipschitz boundary and with  connected $\Omegaext$. 
For any function $u \in E^1\left(\Omegaext\right)$ and any $\varepsilon > 0$, there exists ${R^*}(u, \varepsilon)>0$ such that
\[
\int_{\Omegaext} \left|\nabla u - \nabla (\psi_R u)\right|^2 \, \mathrm{d}x < \varepsilon\qquad\text{for all }R>{R^*}(u, \varepsilon).
\]
\end{lemma}

To establish the convergence of $\sigma_k^\Dir\left(\Omegaext_R\right)$, we separately prove an upper and a lower bound, both of which converge to $\lambda_k\left(A^\Dir\right)$. 

We start with the easier inequality. Extending the eigenfunctions corresponding to $\sigma_k^\Dir\left(\Omegaext_R\right)$ by zero yields the following extension of Lemma \ref{lemma:monotonicitymixed}.

\begin{lemma}\label{lemma1}
Suppose $\Omega \subset \mathbb{R}^n$, $n \geq 3$, is a bounded open set with Lipschitz boundary and with connected $\Omegaext$. For any $R_2>R_1>R_0(\Omega)$ and any $k \in \mathbb{N}$,
\[
\lambda_k\left(A^\Dir\right) \leq \sigma_k^\Dir\left(\Omegaext_{ R_2}\right) \leq \sigma_k^\Dir\left(\Omegaext_{R_1}\right).
\] 
In particular, the limit $\lim_{R \to \infty}\sigma_k^\Dir\left(\Omegaext_R\right)$ exists.
\end{lemma}

\begin{lemma}\label{theo:sigmamuallg}
Suppose $\Omega \subset \mathbb{R}^n$, $n \geq 3$, is a  bounded open set with Lipschitz boundary and with connected $\Omegaext$. For any $\varepsilon > 0$ and any $k \in \N$, there exists a constant ${R_*}(\varepsilon,k)>0$ such that 
\[
\sigma_k^\Dir\left(\Omegaext_{{R_*}}\right) \leq \lambda_k\left(A^\Dir\right)  \left( 1 + \varepsilon (2+ \varepsilon) (  1 +  k   ) \right).
\]
\end{lemma}

\begin{proof}
For each $1 \leq m \leq k$, let $f_m$ be an eigenfunction associated with $\lambda_m\left(A^\Dir\right)$, normalised such that $\int_{\Omegaext} | \nabla f_m |^2 \, \mathrm{d}x=1$. By Lemma~\ref{lemma2}, for any $\varepsilon >0$, there exists ${R^*}>0$ such that the cut-off function satisfies $\operatorname{supp}(\psi_{{R^*}} f_m) \subset B_{2 {{R^*}}}$, $\psi_{{R^*}}(x) = 1$ on $\Omegaext_{{{R^*}}}$, and
\[
\int_{\Omegaext} \left| \nabla f_m - \nabla ( \psi_{{R^*}} f_m )\right|^2 \, \mathrm{d}x < \varepsilon^2.
\]
Note that we can choose a single cut-off function $\psi_{{R^*}}$ that works simultaneously for all $m \leq k$. Define $g_m := \psi_{{R^*}} f_m$ and consider
\[
M := \operatorname{Span} \left\{ g_1, \ldots, g_k\right \}  \subset H^1_D(\Omegaext_{2 {R^*}}).
\]
Since $g_m$ vanishes outside of $\Omegaext_{2 {R^*}}$, we use the same notation for the restriction of $g_m$ to $\Omegaext_{2 {R^*}}$ as for the ``original'' function. Since $\| \cdot \|_{H^1\left(\Omegaext_{2 {R^*}}\right)}$ and $\| \cdot \|_{E^1\left(\Omegaext_{2 {R^*}}\right)}$ are equivalent, \cite[Theorem A.4]{auchmuty2014representations}, $M$ is $k$-dimensional. Consequently,
\begin{equation}\label{eq:upboundsigmaD}
\sigma_k^\Dir(\Omegaext_{2 {R^*}}) \leq \sup_{u\in M, \, u\ne 0} \frac{\int_{\Omegaext_{2 {R^*}}} | \nabla u |^2 \, \mathrm{d}x}{\int_{\partial\Omega} |  u |^2 \, \mathrm{d}S}.
\end{equation}
We now estimate the Rayleigh quotient for each $g_m$ using the triangle inequality, the normalisation of $f_m$, and that $g_m=f_m$ on $\partial\Omega$,
\[
\begin{split}
\frac{\int_{\Omegaext_{2 {R^*}}} | \nabla g_m |^2 \, \mathrm{d}x}{\int_{\partial\Omega} |  g_m |^2 \, \mathrm{d}S}  &\leq \frac{\left( \| \nabla (f_m- g_m)\|_{L^2\left(\Omegaext_{2{R^*}}\right)}+\| \nabla f_m\|_{L^2\left(\Omegaext_{2{R^*}}\right)} \right)^2}{\int_{\partial\Omega} |  f_m |^2 \, \mathrm{d}S}\\  
&\leq \frac{( \varepsilon+1)^2}{\int_{\partial\Omega} |  f_m |^2 \, \mathrm{d}S} = (1+ \varepsilon)^2 \lambda_m\left(A^\Dir\right).
\end{split}
\]
Moreover, the functions  $f_1, \ldots, f_k$ are orthogonal on $\partial\Omega$, hence
\[
\int_{\partial\Omega} g_m g_j \, \mathrm{d}S = \int_{\partial\Omega} f_m f_j \, \mathrm{d}S = 0 \qquad \text{ for }  m \neq j.
\]
For the gradients, we estimate using H\"older's inequality, for $m \neq j$,
\[
\begin{split}
&\left| \int_{\Omegaext} \mydotp{ \nabla g_m ,  \nabla g_j } \, \mathrm{d}x  \right| \leq \int_{\Omegaext} |  \nabla f_m - \nabla g_m| | \nabla g_j | + | \nabla f_m | | \nabla f_j - \nabla g_j |  \, \mathrm{d}x  \\
\leq &\left( \varepsilon^2  \int_{\Omegaext} |  \nabla g_j |^2 \, \mathrm{d}x \right)^\frac{1}{2} + \left( \varepsilon^2  \int_{\Omegaext} |  \nabla f_m |^2 \, \mathrm{d}x \right)^\frac{1}{2} \leq \varepsilon  (1+ \varepsilon)+ \varepsilon.
\end{split}
\]
Now, inserting $u = \sum_{m=1}^k c_m g_m$ into~\eqref{eq:upboundsigmaD}, 
\[
\begin{split}
\frac{\int_{\Omegaext_{2 {R^*}}} | \nabla u |^2 \, \mathrm{d}x}{\int_{\partial\Omega} u^2 \, \mathrm{d}S} &= 
 \frac{  \int_{\Omegaext_{2 {R^*}}} \left( \sum_{m=1}^k c_m^2 | \nabla  g_m|^2+ \sum_{j,m=1, \, j \neq m}^k  c_m c_j   \mydotp{ \nabla  g_j, \nabla  g_m }\right) \, \mathrm{d}x}{\int_{\partial\Omega}  \left(\sum_{m=1}^k c_m^2  g_m^2+  \sum_{j,m=1, \, j \neq m}^k c_m c_j   g_j  g_m\right)  \, \mathrm{d}S} \\
&\leq  \frac{ (1+\varepsilon)^2 \sum_{m=1}^k c_m^2 +  \varepsilon (\varepsilon+2) \sum_{j,m=1, \, j \neq m}^k  |c_m c_j|  }{\sum_{m=1}^k \frac{c_m^2}{\lambda_m\left(A^\Dir\right)} }\\ 
&\leq \lambda_k\left(A^\Dir\right) \left( 1 + \varepsilon (2+ \varepsilon) (  1 +  k   ) \right),
\end{split}
\]
where in the last step we used the Cauchy--Schwarz inequality, and the fact that $\lambda_m\left(A^\Dir\right) \leq \lambda_k\left(A^\Dir\right)$ for $1 \leq m \leq k$. Hence, with $R_*  = 2 R^*$, this completes the proof.
\end{proof}

Combining Lemma~\ref{lemma1} and Lemma~\ref{theo:sigmamuallg}, we obtain the following convergence result. 

\begin{theorem}\label{thm:truncation.equiv.finiteenergy}
Suppose $\Omega \subset \mathbb{R}^n$, $n \geq 3$, is a  bounded open set with Lipschitz boundary and with connected $\Omegaext$. Then, for every $k \in \mathbb{N}$,
\[
\lim_{R \to \infty} \sigma_k^\Dir\left(\Omegaext_R\right) = \lambda_k\left(A^\Dir\right).
\]
\end{theorem}

For the Dirichlet-type problem, convergence of  eigenfunctions holds in a suitable sense. The proof follows a similar strategy to that used for the eigenvalues. Note that convergence cannot be expected for an arbitrary sequence of eigenfunctions. For instance, if $\Omega = B_1 \subset \mathbb{R}^3$, then $\sigma_2^\Dir\left(\Omegaext_R\right)$ has multiplicity three, so a sequence of eigenfunctions $\left(u_{k, R_m}\right)_{m \in \mathbb{N}}$ associated with $\sigma_2^\Dir \left(\Omegaext_{R_m} \right)$  could ``jump'' between different eigenmodes as $m$ varies.

\begin{corollary}\label{lemma3}
Suppose $f_k \in E^1\left(\Omegaext\right)$ is an eigenfunction associated with $\lambda_k\left(A^\Dir\right)$. Then, there exists a sequence $\left( R_m \right)_{m \in \mathbb{N}}$ with $\lim_{m \to \infty} R_m =\infty$ and a sequence $ \left( g_{k,R_m} \right)_{m \in \mathbb{N}} \subset H^1\left(\Omegaext\right)$, such that
\begin{itemize}
\item $g_{k,R_m}(x) = 0$  for $|x|>R_m$ and $g_{k,R_m}|_{\Omegaext_{R_m}}$ is a linear combination of eigenfunctions associated with the eigenvalues $\sigma_j^\Dir(\Omegaext_{R_m})$ for indices $j$ such that $\lambda_j\left(A^\Dir\right) = \lambda_k\left(A^\Dir\right)$; in particular, if $\lambda_k\left(A^\Dir\right)$ is simple, then $g_{k,R_m}|_{\Omegaext_{R_m}}$ is an eigenfunction associated with $\sigma_k^\Dir(\Omegaext_{R_m})$,
\item $\lim_{m \to \infty} \| g_{k,R_m} -f_k \|_{E^1\left(\Omegaext\right)} = 0$.
\end{itemize}
\end{corollary}

\begin{proof}
We assume, without loss of generality, that $\int_{\Omegaext} \mydotp{ \nabla {f}_{j}, \nabla {f}_{k}}  \, \mathrm{d}x = \delta_{j,k} $ (where $\delta_{j,k}$ is the Kronecker delta). For any $R>R_0(\Omega)$, let $u_{k,R} \in H^1\left(\Omegaext_R\right)$ denote an eigenfunction associated with $\sigma_{k}^\Dir\left(\Omegaext_R\right)$ (respectively its extension by zero outside $\Omegaext_R$), normalised by
\[
\int_{\Omegaext_R}  \mydotp{ \nabla {u}_{j,R}, \nabla {u}_{k,R}}  \, \mathrm{d}x = \int_{\Omegaext}  \mydotp{ \nabla {u}_{j,R}, \nabla {u}_{k,R}}  \, \mathrm{d}x = \delta_{j,k}. 
\]
By Theorem~\ref{thm:truncation.equiv.finiteenergy},
\[
\lim_{R \to \infty} \int_{\partial\Omega} | u_{k,R}|^2 \, \mathrm{d}S = \lim_{R \to \infty}\frac{1}{\sigma_k^\Dir\left(\Omegaext_R\right)} = \frac{1}{\lambda_k\left(A^\Dir\right)} = \int_{\partial\Omega} | f_k|^2 \, \mathrm{d}S.
\]
Hence, the norm $\| {u_{k,R}} \|_{E^1\left(\Omegaext\right)}$ is bounded. Thus, there exists a weakly convergent subsequence in $E^1\left(\Omegaext\right)$, which we again denote by $u_{k,R}$, with weak limit $\widehat{u}_k \in E^1\left(\Omegaext\right)$. By compactness of the trace operator, $u_{k,R}\to\widehat{u}_k$ in $L^2(\partial\Omega)$ and hence
\[
\int_{\partial\Omega}|\widehat{u}_k|^2\,\mathrm{d}S=\frac{1}{\lambda_k\left(A^\Dir\right)}>0.
\]
For an arbitrary $v\in E^1\left(\Omegaext\right)$, Lemma~\ref{lemma2} provides a cutoff $\psi_{R'}v$ which agrees with $v$ on $\partial\Omega$, is an admissible test function on $\Omegaext_R$ for all sufficiently large $R$, and approximates $v$ in Dirichlet energy. Using $\|\nabla u_{k,R}\|_{L^2}=1$, we may therefore pass to the limit in the truncated weak equation to obtain
\[
\int_{\Omegaext}\mydotp{\nabla\widehat{u}_k,\nabla v}\,\mathrm{d}x
=\lambda_k\left(A^\Dir\right)\int_{\partial\Omega}\widehat{u}_kv\,\mathrm{d}S,
\qquad v\in E^1\left(\Omegaext\right).
\]
Thus, $\widehat{u}_k$ is a Steklov eigenfunction with associated eigenvalue $\lambda_k\left(A^\Dir\right)$. Taking $v=\widehat{u}_k$ gives $\int_{\Omegaext}|\nabla\widehat{u}_k|^2\,\mathrm{d}x=1$. Hence weak convergence, equality of the Dirichlet-energy norms, and strong convergence of the traces imply $u_{k,R}\to\widehat{u}_k$ strongly in $E^1\left(\Omegaext\right)$.

To obtain the approximation of the prescribed eigenfunction $f_k$, apply the preceding argument simultaneously to all indices $j$ such that
\[
\lambda_j\left(A^\Dir\right)=\lambda_k\left(A^\Dir\right).
\]
There are only finitely many such indices, so a diagonal argument gives a sequence $R_m\to\infty$ for which the corresponding $u_{j,R_m}$ converge in $E^1\left(\Omegaext\right)$ to limits $h_j$. By the first part of the proof, each $h_j$ is an eigenfunction associated with $\lambda_k\left(A^\Dir\right)$. Since $u_{i,R_m}$ is an admissible test function in the weak equation for $u_{j,R_m}$, the normalisation above gives
\[
\int_{\partial\Omega}u_{i,R_m}u_{j,R_m}\,\mathrm{d}S
=\frac{\delta_{i,j}}{\sigma_j^\Dir\left(\Omegaext_{R_m}\right)}.
\]
Passing to the limit,
\[
\int_{\partial\Omega} h_i h_j\,\mathrm{d}S
= \frac{\delta_{i,j}}{\lambda_k\left(A^\Dir\right)}.
\]
Hence these limits are linearly independent; since their number is the multiplicity of $\lambda_k\left(A^\Dir\right)$, they form a basis of the corresponding eigenspace. Write $f_k$ as a linear combination of these $h_j$ and define $g_{k,R_m}$ by taking the same linear combination of the corresponding $u_{j,R_m}$. Then $g_{k,R_m}\to f_k$ in $E^1\left(\Omegaext\right)$, and the required properties follow.
\end{proof}

\subsubsection{Finite-energy functions and Helmholtz equation}\label{sec:AHvGC}

In this section, similarly to \S\ref{sec:AHvAtE}, we show that as $\Lambda \searrow 0$, the eigenvalues $\mu_k\left(\Omegaext,\Lambda\right)$ converge to $\lambda_k\left(A^\Dir\right)$ and that the associated eigenfunctions converge in the sense of Corollary~\ref{lemma5}.

Since $H^1\left(\Omegaext\right) \subset E^1\left(\Omegaext\right)$, the variational characterisations~\eqref{eq:mu_k} and~\eqref{eq:lambdak} immediately imply 
\[
\lambda_k\left(A^\Dir\right) < \mu_k(\Omegaext,\Lambda_1) < \mu_k(\Omegaext,\Lambda_2).
\]
for any $k \in \mathbb{N}$ and any $0<\Lambda_1<\Lambda_2$. Thus, to conclude the convergence of  $\mu_k\left(\Omegaext,\Lambda\right)$, it suffices to establish an upper bound which converges to $\lambda_k\left(A^\Dir\right)$ as $\Lambda \searrow 0$.

\begin{theorem}\label{theo:lambdamuallg}
Suppose $\Omega \subset \mathbb{R}^n$, $n \geq 3$, is a  bounded open set with Lipschitz boundary and with  connected $\Omegaext$, and let $k \in \mathbb{N}$. Then,
\[
\lim_{\Lambda \searrow 0} \mu_k\left(\Omegaext,\Lambda\right) = \lambda_k\left(A^\Dir\right).
\]
\end{theorem}

\begin{proof}
For $1 \leq m \leq k$, let $f_m \in E^1\left(\Omegaext\right)$ be the eigenfunction associated with $\lambda_m\left(A^\Dir\right)$ normalised to $\int_{\Omegaext} | \nabla f_m |^2 \, \mathrm{d}x=1$. By Lemma~\ref{lemma2}, for any $\varepsilon >0$, there exist a constant ${R^*}(\varepsilon,m)>0$ and functions $g_{m}$ with $\operatorname{supp}(g_m) \subset B_{2 {R^*}}$ such that $g_m(x) = f_m(x)$ in $B_{{R^*}} \setminus {\Omega}$ and
\[
\int_{\Omegaext} | \nabla f_m - \nabla g_m|^2 \, \mathrm{d}x < \varepsilon^2.
\]
Define $M := \operatorname{Span} \left\{ g_1, \ldots, g_k  \right\} \subset H^1\left(\Omegaext\right)$. Since the restrictions of $g_m$ to $\partial \Omega$ are part of a basis, $\operatorname{dim}(M) = k$. 
Using the variational characterisation of $\mu_k\left(\Omegaext,\Lambda\right)$, we obtain
\[
\mu_k\left(\Omegaext,\Lambda\right) \leq \sup_{\substack{u \in M\\u\ne 0}} \frac{\Lambda^2 \int_{\Omegaext} |u|^2 \,\mathrm{d}x + \int_{\Omegaext} |\nabla u|^2 \, \mathrm{d}x}{	\int_{\partial\Omega} |u|^2 \, \mathrm{d}S}.
\]
We now show that the right-hand side tends to $\lambda_k\left(A^\Dir\right)$ as $\Lambda \searrow 0$. Since $\| \nabla f_m - \nabla g_m\|_{L^2\left(\Omegaext\right)} < \varepsilon$ and $g_m = f_m$ on $\partial\Omega$, we estimate
\[
\begin{split}
\frac{\Lambda^2 \int_{\Omegaext} |g_m|^2 \,\mathrm{d}x + \int_{\Omegaext} |\nabla g_m|^2 \, \mathrm{d}x}{	\int_{\partial\Omega} |g_m|^2 \, \mathrm{d}S} &\leq \frac{\Lambda^2 \int_{\Omegaext} |g_m|^2 \,\mathrm{d}x + \left( \| \nabla f_m \|_{L^2\left(\Omegaext\right)} +\varepsilon \right)^2}{	\int_{\partial\Omega} |f_m|^2 \, \mathrm{d}S } \\
&= \lambda_m\left(A^\Dir\right) \left( \Lambda^2 \int_{\Omegaext} |g_m|^2 \,\mathrm{d}x + \left( 1 +\varepsilon \right)^2 \right).
\end{split}
\]
Since the $g_m$ are independent of $\Lambda$, it follows that
\[
\lim_{\Lambda \searrow 0} \frac{\Lambda^2 \int_{\Omegaext} |g_m|^2 \,\mathrm{d}x + \int_{\Omegaext} |\nabla g_m|^2 \, \mathrm{d}x}{	\int_{\partial\Omega} |g_m|^2 \, \mathrm{d}S} \leq \lambda_m\left(A^\Dir\right) (1+ \varepsilon)^2.
\]
As in the proof of Lemma~\ref{theo:sigmamuallg}, 
\[
 \int_{\partial\Omega} g_i g_j \, \mathrm{d}S =0 \quad \text{ and } \quad
 \left| \int_{\Omegaext} \mydotp{ \nabla g_i ,  \nabla g_j } \, \mathrm{d}x  \right|  < 3\varepsilon \qquad \text{for} \, i \neq j.
\]
Now let $u = \sum_{m=1}^k c_m g_m \in M$. Then, as $\Lambda \searrow 0$,
\[
\begin{split}
&\lim_{ \Lambda \searrow 0 }\frac{\Lambda^2 \int_{\Omegaext} | u |^2 \, \mathrm{d}x + \int_{\Omegaext} | \nabla u |^2 \, \mathrm{d}x}{\int_{\partial\Omega} u^2 \, \mathrm{d}S} =  \frac{\sum_{m,j=1}^k c_m c_j \int_{\Omegaext}  \mydotp{ \nabla  g_m, \nabla g_j } \, \mathrm{d}x}{\sum_{m,j=1}^k c_m c_j \int_{\partial\Omega}    g_m g_j \, \mathrm{d}S} \\
 \leq & \frac{\sum_{m=1}^k c_m^2 (1+ \varepsilon)^2 + 3 \varepsilon \sum_{\substack{m, j=1, \, m \neq j}}^k |c_m c_j| }{\sum_{m=1}^k c_m^2 \lambda_m\left(A^\Dir\right)^{-1} } \leq  \lambda_k\left(A^\Dir\right) \left( (1 + \varepsilon)^2 + 3k \varepsilon \right).
 \end{split}
\]
Since $\varepsilon > 0$ was arbitrary, the claim follows.
\end{proof}

Using the convergence of $\mu_k\left(\Omegaext,\Lambda\right)$, we can now establish the convergence of the associated eigenfunctions, similar to Corollary~\ref{lemma3}.

\begin{corollary}\label{lemma5}
Suppose $f_k \in E^1\left(\Omegaext\right)$ is an eigenfunction associated with $\lambda_k\left(A^\Dir\right)$. Then, there exist a sequence $\left( \Lambda_m \right)_{m \in \mathbb{N}} \subset (0, \infty)$ with $\lim_{m \to \infty} \Lambda_m = 0$ and a sequence $ \left( g_{k,\Lambda_m} \right)_{m \in \mathbb{N}} \subset H^1\left(\Omegaext\right)$, such that
\begin{itemize}
\item $g_{k,\Lambda_m}$ is a linear combination of eigenfunctions associated with the eigenvalues $\mu_j(\Omegaext,\Lambda_m)$ for indices $j$ such that $\lambda_j\left(A^\Dir\right) = \lambda_k\left(A^\Dir\right)$ (in particular, if $\lambda_k\left(A^\Dir\right)$ is simple, then $g_{k,\Lambda_m}$ is an eigenfunction associated with $\mu_k(\Omegaext,\Lambda_m)$), and $\lim_{m \to \infty}\mu_k(\Omegaext,\Lambda_m) = \lambda_k\left(A^\Dir\right)$,
\item $\lim_{m \to \infty} \| g_{k,\Lambda_m} -f_k \|_{E^1\left(\Omegaext\right)} = 0$.
\end{itemize}
\end{corollary}

\begin{proof}
For each index $j$ such that $\lambda_j\left(A^\Dir\right)=\lambda_k\left(A^\Dir\right)$ and each $\Lambda>0$, let $u_{j,\Lambda}\in H^1\left(\Omegaext\right)$ be an eigenfunction associated with $\mu_j\left(\Omegaext,\Lambda\right)$. For a given $\Lambda$, we choose these eigenfunctions pairwise orthogonal in $L^2(\partial\Omega)$ and normalise them by
\[
\int_{\Omegaext}|\nabla u_{j,\Lambda}|^2\,\mathrm{d}x=1.
\]

Since $H^1\left(\Omegaext\right)\subset E^1\left(\Omegaext\right)$, the variational characterisation~(3.5) of $\lambda_1\left(A^\Dir\right)>0$ gives
\[
\int_{\partial\Omega}|u_{j,\Lambda}|^2\,\mathrm{d}S
\leq \frac{1}{\lambda_1\left(A^\Dir\right)}.
\]
Thus the family $u_{j,\Lambda}$ is bounded in $E^1\left(\Omegaext\right)$. Choose any sequence $\Lambda_m\searrow0$. Since there are only finitely many indices $j$ with $\lambda_j\left(A^\Dir\right)=\lambda_k\left(A^\Dir\right)$, we may successively extract subsequences for these finitely many indices. After relabelling the final subsequence as $\Lambda_m$, we have simultaneously, for every such $j$,
\[
u_{j,\Lambda_m}\rightharpoonup h_j \quad\text{weakly in }E^1\left(\Omegaext\right),
\qquad
u_{j,\Lambda_m}\longrightarrow h_j \quad\text{in }L^2(\partial\Omega),
\]
where the second convergence follows from the compactness of the trace operator.

For every $m$ and every $v\in H^1\left(\Omegaext\right)$, the weak eigenvalue equation is
\[
\int_{\Omegaext}\mydotp{\nabla u_{j,\Lambda_m},\nabla v}\,\mathrm{d}x
+\Lambda_m^2\int_{\Omegaext}u_{j,\Lambda_m}v\,\mathrm{d}x
=\mu_j\left(\Omegaext,\Lambda_m\right)\int_{\partial\Omega}u_{j,\Lambda_m}v\,\mathrm{d}S.
\]
Taking $v=u_{j,\Lambda_m}$ gives
\[
1+\Lambda_m^2\int_{\Omegaext}|u_{j,\Lambda_m}|^2\,\mathrm{d}x
=\mu_j\left(\Omegaext,\Lambda_m\right)\int_{\partial\Omega}|u_{j,\Lambda_m}|^2\,\mathrm{d}S.
\]
By Theorem~5.6, $\mu_j\left(\Omegaext,\Lambda_m\right)\to\lambda_j\left(A^\Dir\right)=\lambda_k\left(A^\Dir\right)$. Together with the preceding boundary estimate, this shows that the right-hand side of the last identity is uniformly bounded in $m$, and hence $\Lambda_m\|u_{j,\Lambda_m}\|_{L^2\left(\Omegaext\right)}$ is bounded. Therefore, for every fixed $v\in H^1\left(\Omegaext\right)$,
\[
\Lambda_m^2\left|\int_{\Omegaext}u_{j,\Lambda_m}v\,\mathrm{d}x\right|
\leq \Lambda_m\bigl(\Lambda_m\|u_{j,\Lambda_m}\|_{L^2\left(\Omegaext\right)}\bigr)\|v\|_{L^2\left(\Omegaext\right)}\longrightarrow0,
\]
because $\Lambda_m\to0$ and the factor in parentheses is bounded.

We may now pass to the limit in the preceding weak equation: the first term converges by weak convergence in $E^1$, the $\Lambda_m^2$-term tends to zero as just shown, the boundary term converges by the strong convergence of the traces, and the eigenvalues converge by Theorem~5.6. Thus
\[
\int_{\Omegaext}\mydotp{\nabla h_j,\nabla v}\,\mathrm{d}x
=\lambda_k\left(A^\Dir\right)\int_{\partial\Omega}h_jv\,\mathrm{d}S,
\qquad v\in H^1\left(\Omegaext\right).
\]
Using Lemma~5.1 to approximate an arbitrary function in $E^1\left(\Omegaext\right)$ by compactly supported $H^1$ functions extends this identity to every $v\in E^1\left(\Omegaext\right)$. Thus $h_j$ is an eigenfunction associated with $\lambda_k\left(A^\Dir\right)$ once we know that it is nonzero.

The preceding energy identity also gives
\[
\frac{1}{\displaystyle\int_{\partial\Omega}|u_{j,\Lambda_m}|^2\,\mathrm{d}S}
\leq \mu_j\left(\Omegaext,\Lambda_m\right).
\]
Hence the strong convergence of the traces and Theorem~5.6 yield
\[
\int_{\partial\Omega}|h_j|^2\,\mathrm{d}S
\geq \frac{1}{\lambda_k\left(A^\Dir\right)}>0,
\]
so $h_j\neq0$. Taking $v=h_j$ in the limiting weak equation and using weak lower semicontinuity of the $L^2$ norm, together with the normalisation $\int_{\Omegaext}|\nabla u_{j,\Lambda_m}|^2\,\mathrm{d}x=1$, gives
\[
1\leq \lambda_k\left(A^\Dir\right)\int_{\partial\Omega}|h_j|^2\,\mathrm{d}S
=\int_{\Omegaext}|\nabla h_j|^2\,\mathrm{d}x\leq1.
\]
Therefore
\[
\int_{\Omegaext}|\nabla h_j|^2\,\mathrm{d}x=1,
\qquad
\int_{\partial\Omega}|h_j|^2\,\mathrm{d}S=\frac{1}{\lambda_k\left(A^\Dir\right)}.
\]
The gradients therefore converge strongly in $L^2\left(\Omegaext\right)$, since they converge weakly and their norms converge. Together with the strong convergence of the traces, this implies
\[
u_{j,\Lambda_m}\longrightarrow h_j \quad\text{strongly in }E^1\left(\Omegaext\right).
\]

Finally, the functions $u_{j,\Lambda_m}$ were chosen pairwise orthogonal in $L^2(\partial\Omega)$, and hence their limits $h_j$ are pairwise orthogonal on $\partial\Omega$. They are nonzero eigenfunctions associated with $\lambda_k\left(A^\Dir\right)$, and their number is precisely the multiplicity of this eigenvalue. Thus they form a basis of the corresponding eigenspace. Write the prescribed eigenfunction $f_k$ as a linear combination of these $h_j$ and define $g_{k,\Lambda_m}$ by taking the same linear combination of the corresponding $u_{j,\Lambda_m}$. Then
\[
\|g_{k,\Lambda_m}-f_k\|_{E^1\left(\Omegaext\right)}\longrightarrow0,
\]
and $g_{k,\Lambda_m}$ has the required linear-combination property. In the simple-eigenvalue case only the index $j=k$ occurs, so $g_{k,\Lambda_m}$ is itself an eigenfunction associated with $\mu_k\left(\Omegaext,\Lambda_m\right)$.
\end{proof}

In summary, the four different available approaches to the exterior problem, introduced in \S\ref{subsec:AuchmutyHan}, \S\ref{subsec:truncatedomains}, \S\ref{subsec:Greben}, and~\S\ref{subsubsec:extlayer},  are equivalent in dimensions $n \geq 3$. 

\subsubsection{Vanishing flow and Neumann truncation}\label{sec:VFNT}

In \S\ref{sec:AHvAtE} we have seen that the finite energy approach is equivalent to approximating the exterior domain by truncated domains equipped with Dirichlet boundary conditions on the outer boundary. By contrast, if Neumann boundary condition is imposed on the truncation boundary, one obtains a different limiting spectrum in dimensions $n \geq 3$, see Remark~\ref{rem:other}. In dimension $n=2$, however, the Dirichlet and Neumann truncations lead to the same spectrum in the limit; this will be shown in \S\ref{subusbsec:CMtrunc}.

In \S\ref{subsec:vanishflow} we introduced, via layer potentials, an alternative formulation of the exterior Steklov problem. In this setting, functions are not required to vanish at infinity, but their gradients must decay faster than $|x|^{1-n}$. We now prove that this formulation is equivalent to the Neumann truncation. Recall that the eigenvalues $\lambda_k\left(A^\Neu\right)$ of the operator $A^\Neu$, defined in \S\ref{subsec:truncatedomains}, satisfy 
\[
\lambda_k\left(A^\Neu\right)=\lim_{R \to \infty} \sigma^\Neu_k\left(\Omegaext_R\right), \qquad k\in\mathbb{N}.
\] 

\begin{theorem}\label{thm:FVTN}
Suppose $\Omega \subset \mathbb{R}^n$, $n \geq 3$, is a bounded open set with Lipschitz boundary and with connected $\Omegaext$. Then the spectrum of  \eqref{eq:henrici} is discrete, and eigenvalues of  \eqref{eq:henrici} coincide with those of $A^\Neu$ with account of multiplicities.
\end{theorem}

\begin{proof}
Let $\lambda^\Neu \in \mathbb{R} $ be an eigenvalue of $A^\Neu$, and let $u^\Neu \in W\left( \Omegaext \right)$ be such that the trace of $u^\Neu$ is an eigenfunction associated with $\lambda^\Neu$, that is,  $\Delta u^\Neu =0$ and $\partial_\nu u^\Neu = A^\Neu \left(\left.u^\Neu\right|_{\partial \Omega} \right)$. Let $\left\langle u^\Neu \right\rangle$ denote the average of $u^\Neu$, as defined in \eqref{eq:averageu}. Then, by definition, $w := u^\Neu - \left\langle u^\Neu \right\rangle$ is in $W^\Dir\left(\Omegaext\right)$, and therefore, by \eqref{eq:WD=E1},  $w \in E^1\left(\Omegaext\right)$. Without loss of generality, assume $\Omega \subset B_1$. As in the proof of Proposition~\ref{pro:decay}, which will be presented later, we expand $w$ in spherical harmonics
\[ 
w(x) = \sum_{k=0}^\infty \sum_{i=1}^{d_{n,k}} \myscal{ w, {u}_{(k,i)} }_{L^2\left(\partial B_1\right)} {u}_{(k,i)}(x), 
\]
where ${u}_{(k,i)}\in E^1\left(B_1^\ext\right)$, given by ${u}_{(k,i)}(r, \theta) =  r^{2-n-k} Y_{k,i}(\theta)$, are Steklov eigenfunctions in the exterior of the unit ball whose traces form an orthonormal basis of $L^2\left(\partial B_1\right)$, see \S\ref{examp:balln3}.  Following the proof of Proposition~\ref{pro:decay} and using the decay at infinity, we obtain that this expansion converges uniformly in $\Omegaext$.  

From the explicit form of ${u}_{(k,i)}$ it follows that $w$ (and thus $u^\Neu$) satisfies the vanishing flow condition \eqref{eq:alternbehav} if and only if $\myscal{ w, {u}_{(0,1)} }_{L^2\left(\partial B_1\right)} = 0$.
To establish this, note first that since $w$ is harmonic and $\partial_\nu w = \partial_\nu u^\Neu = \lambda^\Neu u^\Neu$ on $\partial\Omega$, we have
\[
\begin{split}
0 &= \int_{\Omegaext \cap B_1} \Delta w \, \mathrm{d}x = \int_{\partial\Omega} \partial_\nu w \, \mathrm{d}S + \int_{\partial B_1} \partial_\nu w \, \mathrm{d}S\\
&= \lambda^\Neu \int_{\partial\Omega}  u^\Neu \, \mathrm{d}S + \int_{\partial B_1} \partial_\nu w \, \mathrm{d}S=\int_{\partial B_1} \partial_\nu w \, \mathrm{d}S.
\end{split}
\]
In the last step we used that the first eigenvalue vanishes and the associated eigenfunction is a constant, and, therefore, any other eigenfunction is orthogonal to a constant. Furthermore, $\partial_\nu$ coincides with the radial derivative on $\partial B_1$. Hence,
\[
0 = \int_{\partial B_1} \partial_\nu w \, \mathrm{d}S = \sum_{k=0}^\infty \sum_{i=1}^{d_{n,k}} \myscal{ w, {u}_{(k,i)} }_{L^2\left(\partial B_1\right)}  (2-n-k) \int_{\partial B_1}    Y_{k,i}(\theta)  \, \mathrm{d}S.
\]
Since the spherical harmonics of degree $k \neq 0$ integrate to zero, only the $k=0$ term remains. Therefore $0 =\myscal{ w, {u}_{0,1} }_{L^2\left(\partial B_1\right)} c_{0,1} (2-n) \int_{\partial B_1}    Y_{0,1}(\theta)  \, \mathrm{d}S $. Hence, $ \myscal{ w, {u}_{0,1} }_{L^2\left(\partial B_1\right)} = 0$. Thus $u^\Neu$ satisfies the vanishing flow condition \eqref{eq:alternbehav}, and $\left(\sigma^\Neu, u^\Neu\right)$ is an eigenpair of \eqref{eq:henrici}.

Conversely, if $(\chi, U)$ is an eigenpair of \eqref{eq:henrici}, then $U \in W\left(\Omegaext\right)$ by definition, hence $\chi$ and $U|_{\partial \Omega}$ are an eigenvalue and eigenfunction of $A^\Neu$.
\end{proof}

\subsection{Equivalence in dimension two}\label{section:twod}

We show that formulation \eqref{eq:problemtransform} is equivalent to the approach via the Helmholtz equation, the truncated domain method and the approach via layer potentials. Recall that the approach of Auchmuty and Han does not extend to two dimensions, and Arendt and ter Elst did not study that case either, see Remark \ref{rem:not2d}.

Moreover, the conformal mapping framework leads to several results analogous to those known for bounded domains.

In the following, when discussing \eqref{eq:problemtransform}, we switch from the complex formulation to real variables. We now view $\Omega^*$ as a subset of $\mathbb{R}^2$ and represent points as $x = (x_1, x_2)$ instead of $z = x_1 + \ir x_2$, so that $\phi: \left(x_1,x_2\right)\mapsto \left(\frac{x_1}{|x|^2},  -\frac{x_2}{|x|^2}\right)$, and $\left|\phi'(x)\right|=\frac{1}{|x|^2}$. Accordingly, the derivatives and normal vectors are interpreted in the standard real sense.

\subsubsection{Conformal mapping and  truncated domains}\label{subusbsec:CMtrunc}

In contrast to higher dimensions, the choice of boundary condition at the outer boundary of the truncated domain plays no role in two dimensions when $R \to \infty$. Loosely speaking, this is because under $\phi$, defined in \eqref{eq:phi}, the outer boundary $\partial B_R$ becomes a small hole in $\Omega^*$ whose influence vanishes as $R \to \infty$.

\begin{theorem}\label{theo:sigmaD*}
Let $\Omega \subset \mathbb{R}^2$ be a bounded open set with Lipschitz boundary and with connected $\Omegaext$,  and such that $\Omega$ contains the origin. Then, for any $k \in \mathbb{N}$,
\[
\lim_{R \to \infty} \sigma_k^\Dir\left(\Omegaext_R\right) = \xi_k\left(\Omegaext\right),
\]
where $\sigma_k^\Dir\left(\Omegaext_R\right)$ is defined as in \S\ref{subsec:truncatedomains}.
\end{theorem}

\begin{proof}
Let $u_{k,R} \in H^1_D\left(\Omegaext_R\right)$ denote an eigenfunction associated with $\sigma_k^\Dir\left(\Omegaext_R\right)$. By slightly abusing the notation, $u_{k,R}$ also denotes the extension to $\Omegaext$ by zero. In addition, we define $f_{k,R} \in H^1(\Omega^*)$ by
\[
f_{k,R} : \Omega^* \to \mathbb{R}, \qquad f_{k,R}(x) = 
\begin{cases}
u_{k,R} \left( \phi(x) \right)\qquad &\text{ if } |x|\geq\frac{1}{R},\\
0 \qquad&\text{ if } |x|<\frac{1}{R}.
\end{cases} 
\]
Since $\phi$ is a conformal map,
\[
\frac{\int_{\Omega^* } | \nabla f_{k,R}  |^2 \, \mathrm{d}x}{\int_{\partial\Omega^*} |\phi'| \cdot |  f_{k,R}  |^2 \, \mathrm{d}S} =\frac{\int_{\Omega^*\setminus \overline{B_{1/R}}} | \nabla f_{k,R}  |^2 \, \mathrm{d}x}{\int_{\partial\Omega^*} |\phi'| \cdot |  f_{k,R}  |^2 \, \mathrm{d}S} = \frac{\int_{\Omegaext_R} | \nabla u_{k,R}  |^2 \, \mathrm{d}x}{\int_{\partial\Omega}  |  u_{k,R}  |^2 \, \mathrm{d}S}.
\]
For fixed $R > R_0(\Omega)$, define the $j$-dimensional subspaces
\[
M(j) := \operatorname{Span} \left\{ f_{k,R} \mid 1 \leq k \leq j  \right\} \subset H^1(\Omega^*) 
\]
and
\[
N(j) := \operatorname{Span} \left \{ u_{k,R} \mid 1 \leq k \leq j  \right\} \subset  H^1_D\left(\Omegaext_R\right).
\]
Then,
\[
\xi_j\left(\Omegaext\right) \leq \sup_{0 \neq u \in M(j)} \frac{\int_{\Omega^* } | \nabla u  |^2 \, \mathrm{d}x}{\int_{\partial\Omega^*} |\phi'| \cdot |  u  |^2 \, \mathrm{d}S} = \sup_{0 \neq u \in N(j)} \frac{\int_{\Omegaext_R} | \nabla u  |^2 \, \mathrm{d}x}{\int_{\partial\Omega}  |  u  |^2 \, \mathrm{d}S} = \sigma_j^\Dir\left(\Omegaext_R\right).
\]
For the reverse inequality, let $f_k \in H^1(\Omega^*)$ denote eigenfunctions associated with $\xi_k\left(\Omegaext\right)$, and define
\[
g_{k,R}:\Omegaext \to \R, \qquad g_{k,R}(x) = 
\begin{cases}
f_k\left( \phi^{-1}(x) \right)  \quad&\text{ if } x \in \Omegaext_R,\\
f_k\left( \phi^{-1}(x) \right) \left( \frac{\log(R^2)- \log(|x|)}{\log(R^2)- \log(R)} \right) \quad&\text{ if } R < |x|\leq R^2, \\
0 \quad&\text{ if } |x|>R^2.
\end{cases}
\]
Since $|\nabla f_k| \in L^2(\Omega^*)$, a straightforward calculation shows $\lim_{R \to \infty} \int_{B_{R^2} \setminus \overline{B_R}} |\nabla g_{k,R}|^2 \, \mathrm{d}x = 0$. Furthermore,
\[
\int_{\partial\Omega} | g_{k,R}|^2 \, \mathrm{d}S = \int_{\partial\Omega^*}  |\phi'| \cdot |  f_k  |^2 \, \mathrm{d}S  \quad \text{and} \quad 
\int_{\Omegaext_R} |\nabla g_{k,R}|^2 \, \mathrm{d}x = \int_{\Omega^* \setminus \overline{B_{1/R}}} |\nabla f_k|^2 \, \mathrm{d}x.
\]
Hence,
\[
\lim_{R \to \infty} \frac{\int_{\Omegaext_R} |\nabla g_{k,R}|^2 \, \mathrm{d}x}{\int_{\partial\Omega} | g_{k,R}|^2 \, \mathrm{d}S} = \frac{\int_{\Omega^* } |\nabla f_k|^2 \, \mathrm{d}x}{\int_{\partial\Omega^*}  |\phi'| \cdot |  f_k  |^2 \, \mathrm{d}S} = \xi_k\left(\Omegaext\right).
\]
So, in order to find an upper bound for $\sigma_k^\Dir(\Omegaext_{R^2})$, we consider
\[
L(j) := \operatorname{Span} \left\{ g_{k,R},  1 \leq k \leq j  \right\} \subset   H^1_D(\Omegaext_{R^2}).
\]
For any $u = \sum_{p=1}^j c_p g_{p,R}$, one estimates
\[
\lim_{R \to \infty} \sqrt{\int_{\Omegaext_{R^2} \setminus \overline{B_R}} |\nabla u|^2 \, \mathrm{d}x} \leq  \lim_{R \to \infty} \sum_{p=1}^j |c_p| \sqrt{\int_{\Omegaext_{R^2} \setminus \overline{B_R}} |\nabla g_{p,R}|^2 \, \mathrm{d}x} = 0
\]
and
\[
\lim_{R \to \infty}\int_{\Omegaext_{R}} |\nabla u|^2 \, \mathrm{d}x  = \int_{\Omega^* } \left| \nabla \left( \sum_{p=1}^j c_p  f_p \right) \right|^2 \mathrm{d}x,\quad \int_{\partial\Omega} | u|^2 \, \mathrm{d}S =  \int_{ \partial  \Omega^* } |\phi'| \cdot \left| \sum_{p=1}^j c_p  f_p\right|^2  \mathrm{d}S.
\]
Therefore,
\[
\begin{split}
\lim_{R \to \infty} \sigma_j^\Dir(\Omegaext_{R^2})&\leq \lim_{R \to \infty} \sup_{0 \neq u \in L(j)} \frac{\int_{\Omegaext_{R^2}} | \nabla u  |^2 \, \mathrm{d}x}{\int_{\partial\Omega}  |  u  |^2 \, \mathrm{d}S}\\ 
&\leq \sup_{0 \neq u \in \operatorname{Span} \{f_1, \ldots, f_j\}} \frac{\int_{\Omega^*} |\nabla u|^2 \, \mathrm{d}x }{\int_{\partial\Omega^*} |\phi'| \cdot | u|^2 \, \mathrm{d}S} = \xi_j\left(\Omegaext\right).
\end{split}
\]
Combining both inequalities gives the desired limit.
\end{proof}

While $\xi_j\left(\Omegaext\right) \leq  \sigma_j^\Dir\left(\Omegaext_R\right)$ for any $R > R_0(\Omega)$, the Neumann-truncated eigenvalues will approximate $\xi_j\left(\Omegaext\right)$ from below.

\begin{theorem}\label{theo:convEVNeum}
Let $\Omega \subset \mathbb{R}^2$ be a bounded open set with Lipschitz boundary and with connected $\Omegaext$,  and such that $\Omega$ contains the origin. 
Then, for any $k \in \mathbb{N}$,
\[
\lim_{R \to \infty} \sigma_k^\Neu\left(\Omegaext_R\right) = \xi_k\left(\Omegaext\right),
\]
where $\sigma_k^\Neu\left(\Omegaext_R\right)$ is defined as in \S\ref{subsec:truncatedomains}.
\end{theorem}

\begin{proof}
It is straightforward to verify that $\sigma_k^\Neu\left(\Omegaext_R\right)$ is monotonically increasing in $R$ and bounded above by $\xi_k(\Omegaext)$; see Lemma~\ref{lemma:monotonicitymixed}.

The case $k=1$ is immediate, since $\sigma_1^\Neu\left(\Omegaext_R\right)=\xi_1\left(\Omegaext\right)=0$. Hence assume $k\geq2$. For the reverse inequality, for each $R$ choose eigenfunctions $u_{1,R},\dots,u_{k,R}$ so that $u_{1,R}$ is constant with $\|u_{1,R}\|_{L^2(\partial\Omega)}=1$, while for $j\geq2$ they are normalised by
\[
\int_{\Omegaext_R}\mydotp{\nabla u_{i,R},\nabla u_{j,R}}\,\mathrm{d}x=\delta_{i,j},\qquad 2\leq i,j\leq k.
\]
Define, for $j=1,\dots,k$,
\[
f_{j,R}:\Omega^*\setminus\overline{B_{1/R}}\to\mathbb{R},\qquad f_{j,R}(x)=u_{j,R}(\phi(x)).
\]
Since $\phi$ is conformal,
\[
\int_{\Omega^*\setminus\overline{B_{1/R}}}|\nabla f_{j,R}|^2\,\mathrm{d}x
=\int_{\Omegaext_R}|\nabla u_{j,R}|^2\,\mathrm{d}x
=\begin{cases}0,&j=1,\\1,&2\leq j\leq k.\end{cases}
\]
We now choose $\rho=\rho(R)$ before applying the extension lemma. Let $m_0\geq0$ be such that $\Omega\subset B_{2^{m_0}}$, let $M=M(R)$ be the largest integer with $2^M\leq R$, and consider the disjoint annuli
\[
A_m:=B_{2^m}\setminus B_{2^{m-1}}\subset\Omegaext_R,\qquad m=m_0+1,\dots,M.
\]
Since the total Dirichlet energy of each $u_{j,R}$ is at most one, there is an index $m:=m(R)$ such that
\[
\sum_{j=1}^k\int_{A_{m(R)}}|\nabla u_{j,R}|^2\,\mathrm{d}x\leq\frac{k}{M-m_0}.
\]
Set $\rho(R):=2^{m(R)}$ and restrict $f_{j,R}$ to $\Omega^*\setminus B_{1/\rho(R)}$. Applying \cite[Lemma~4.5]{PostAnne2021} with $\delta=1/\rho(R)$, we extend this restriction harmonically into $B_{1/\rho(R)}$, obtaining a function
\[
\widetilde f_{j,R}\in H^1(\Omega^*)
\]
such that
\[
\widetilde f_{j,R}=f_{j,R}
\qquad\text{on }\Omega^*\setminus B_{1/\rho(R)}.
\] 
We have
\begin{equation}\label{eq:AP21}
\int_{B_{1/\rho(R)}}|\nabla\widetilde f_{j,R}|^2\,\mathrm{d}x
\leq C_1\int_{B_{\rho(R)}\setminus B_{\rho(R)/2}}|\nabla u_{j,R}|^2\,\mathrm{d}x
\leq\frac{C_1k}{M(R)-m_0}\xrightarrow[R\to\infty]{}0,
\end{equation}
for $j=1,\dots,k$, where $C_1$ is independent of $R$.

To use the $\widetilde{f}_{k,R}$'s as trial functions, we define the $k$-dimensional space
\[
M(k,R) := \operatorname{Span} \left\{ \widetilde{f}_{1,R}, \ldots, \widetilde{f}_{k,R} \right \} \subset H^1(\Omega^*).
\] 
Let $w = \sum_{j=1}^k c_j \widetilde{f}_{j,R} \in M(k,R)$. Using that the $u_{j,R}$'s are orthogonal with respect to $\myscal{ \nabla f, \nabla g  }_{L^2\left(\Omegaext_R\right)}$ and with respect to $\myscal{f,  g}_{L^2(\partial\Omega)}$, and that $\widetilde{f}_{1,R}$ is constant, so that it contributes to neither term below, we obtain, using Cauchy--Schwarz inequality, 
\[
\int_{\Omega^*} |\nabla w|^2 \, \mathrm{d}x \leq \sum_{j=2}^k kc_j^2 \int_{B_{1/\rho(R)}}  |\nabla \widetilde{f}_{j,R}|^2 \, \mathrm{d}x + \sum_{j=2}^k c_j^2
\]
and
\[
\int_{\partial\Omega^*} |\phi'| \cdot |  w |^2 \, \mathrm{d}S =   \sum_{j=1}^k c_j^2 \int_{\partial\Omega}  |   u_{j,R} |^2 \, \mathrm{d}S {}= c_1^2 + \sum_{j=2}^k c_j^2 \left(\sigma_j^\Neu\left(\Omegaext_R\right)\right)^{-1}.
\]
With $\varepsilon(R) := \max_{1 \leq j \leq k} \int_{B_{1/\rho(R)}} |\nabla \widetilde{f}_{j,R}|^2  \, \mathrm{d} x$, we obtain
\[
\begin{split}
\frac{\int_{\Omega^*} | \nabla w |^2 \, \mathrm{d}x}{\int_{\partial\Omega^*} |\phi'| \cdot |  w |^2 \, \mathrm{d}S} &\leq 
\frac{ \left( 1+ k \varepsilon (R) \right) \sum_{j=2}^k c_j^2}{c_1^2 +{} \sum_{j=2}^k c_j^2 \left(\sigma_j^\Neu\left(\Omegaext_R\right)\right)^{-1} } \leq \sigma_k^\Neu\left(\Omegaext_R\right) \left( 1+ k \varepsilon (R) \right).
\end{split}
\]
Thus, $\xi_k\left(\Omegaext\right) \leq (1+k \varepsilon(R)) \sigma_k^\Neu\left(\Omegaext_R\right)$ and taking the limit $R \to \infty$ yields the desired convergence.
\end{proof}

We proceed, by analogy with Corollary~\ref{lemma3}, to infer the convergence of the eigenfunctions from Theorem~\ref{theo:sigmaD*}.

\begin{corollary}\label{coro:convEF}
Suppose $f_k \in H^1(\Omega^*)$ is an eigenfunction associated with $\xi_k\left(\Omegaext\right)$. Then, there exists a  sequence $\left( R_m \right)_{m \in \mathbb{N}}$ with $\lim_{m \to \infty} R_m =\infty$ and a sequence $ \left( g_{k,R_m} \right)_{m \in \mathbb{N}} \subset H^1\left(\Omegaext\right)$, such that
\begin{itemize}
\item $g_{k,R_m}(x) = 0$  for $|x|>R_m$ and  $g_{k,R_m}|_{\Omegaext_{R_m}}$ is a linear combination of eigenfunctions associated with the eigenvalues $\sigma_j^\Dir(\Omegaext_{R_m})$ for indices $j$ such that $\xi_j\left(\Omegaext\right) = \xi_k\left(\Omegaext\right)$; in particular, if $\xi_k\left(\Omegaext\right)$ is simple, then $g_{k,R_m}|_{\Omegaext_{R_m}}$ is an eigenfunction associated with $\sigma_k^\Dir(\Omegaext_{R_m})$,
\item $\lim_{m \to \infty} \| g_{k,R_m} \circ \phi -f_k \|_{H^1(\Omega^*)} = 0$.
\end{itemize}
\end{corollary}

\begin{proof}
We distinguish two cases. First, assume $k>1$ and hence $\xi_k(\Omegaext) \neq 0$. 
Let $u_{k,R} \in H^1\left(\Omegaext_R\right)$ be an eigenfunction associated with $\sigma_{k}^\Dir(\Omegaext_R)$, extended by zero to $\Omegaext$ and define
\[
f_{k,R}: \Omega^* \to \R, \qquad f_{k,R}(x) = \begin{cases}
u_{k,R} \left( \phi(x) \right) \qquad& \text{if } \, |x|\geq \frac{1}{R},\\
0 \qquad& \text{if } \, |x|< \frac{1}{R}.
\end{cases}
\]
We may normalise the approximating eigenfunctions so that
\[
\int_{\Omegaext_R} \mydotp{\nabla u_{j,R},\nabla u_{k,R}} \, \mathrm{d}x=\delta_{j,k} \quad \text{ for all }R.
\]
Since $\phi$ is conformal, we also have $\int_{\Omega^*} |\nabla f_{k,R}|^2 \,\mathrm{d}x = 1$. From Theorem~\ref{theo:sigmaD*}, we know
\[
\xi_k\left(\Omegaext\right) = \lim_{R \to \infty} \sigma_k^\Dir(\Omegaext_R) = \lim_{R \to \infty} \frac{1}{\int_{\partial\Omega} | u_{k,R}|^2 \, \mathrm{d}S} = \lim_{R \to \infty} \frac{1}{\int_{\partial\Omega^*}|\phi'| \cdot | f_{k,R}|^2 \, \mathrm{d}S}.
\]
Moreover, as $ \inf_{x \in \partial\Omega^*}|\phi'(x)|>0$, the norms $\|f_{k,R}\|_{H^1(\Omega^*)}$ are bounded. Consequently, there exists a subsequence, again denoted by ${f_{k,R}}$ converging weakly in $H^1(\Omega^*)$ to some function ${g_{k}} \in H^1(\Omega^*)$. Using the logarithmic cutoff from the proof of Theorem~\ref{theo:sigmaD*}, we may pass to the limit in the weak eigenvalue equation, exactly as in the first part of the proof of Corollary~\ref{lemma3}; hence $g_k$ is an eigenfunction associated with $\xi_k\left(\Omegaext\right)$. The convergence of the boundary norms and the limiting eigenvalue equation then imply equality of the Dirichlet-energy norms, so the convergence is strong in $H^1(\Omega^*)$. Applying this simultaneously to all indices $j$ with $\xi_j=\xi_k$ and repeating the final linear-combination step of Corollary~\ref{lemma3} proves the claim.

Consider now the case $k=1$. Then $\xi_1=0$ and the corresponding eigenfunction is constant.
Let $u_{k,R}$ be an eigenfunction associated with $\sigma_{k}^\Dir(\Omegaext_R)$ normalised to $\| u_{k,R}\|_{L^2(\partial\Omega)} = 1$.  The convergence of the eigenvalue implies
\[
\lim_{R \to \infty} \int_{\Omega^*} | \nabla u_{k,R} \circ \phi|^2 \, \mathrm{d}x = 0.
\]
By the Poincar\'e inequality, this implies that $\left(u_{k,R} \circ \phi\right) - \frac{\int_{\Omega^*}u_{k,R} \circ \phi \, \mathrm{d}x}{|\Omega^*|}$ converges strongly to zero in $H^1(\Omega^*)$. Moreover, since 
\[
\begin{split}
\int_{\partial\Omega^*}  |  u_{k,R} \circ \phi|^2 \, \mathrm{d}S &\leq \frac{1}{\inf_{x \in \partial\Omega^*}|\phi'(x)|} \int_{\partial\Omega^*} |\phi'| \cdot |  u_{k,R} \circ \phi|^2 \, \mathrm{d}S  \\
&= \frac{1}{\inf_{x \in \partial\Omega^*}|\phi'(x)|} \int_{\partial\Omega} |  u_{k,R}|^2 \, \mathrm{d}S = \frac{1}{\inf_{x \in \partial\Omega^*}|\phi'(x)|}   < \infty,
\end{split}
\]
the average value of $u_{k,R} \circ \phi$ remains bounded. Therefore, there exists a subsequence that converges strongly in $H^1(\Omega^*)$ to a constant $c \in \mathbb{R}$. The boundary normalisation gives $c\neq0$, and rescaling this subsequence if necessary gives convergence to the prescribed eigenfunction $f_1$.
\end{proof}

\begin{remark}
An alternative approach to the proof of Theorem~\ref{theo:convEVNeum} would be to recognise that \eqref{eq:PDEArendtElstNeumann} is isospectral to a mixed Steklov-Neumann problem on $\Omega^* \setminus B_{1/R}$ with density $|\phi'|$ on $\partial \Omega^*$. Since the density does not depend on $R$, and we have the uniform bound \eqref{eq:AP21} on the harmonic extension in the small hole,  \cite[Propositions 4.8 and 4.11]{girouard2021continuity} can be used to establish the convergence of the eigenvalues and to prove Corollary~\ref{coro:convEF}.
\end{remark}

As in the Dirichlet case, we next prove that the eigenfunctions associated with $\sigma_k^\Neu\left(\Omegaext_R\right)$, when pulled back via the conformal map $\phi$, converge in the sense of Corollary~\ref{coro:convEFNeumann} to the eigenfunctions of the limiting problem on $\Omega^*$.

\begin{corollary}\label{coro:convEFNeumann}
Suppose $f_k \in H^1(\Omega^*)$ is an eigenfunction associated with $\xi_k\left(\Omegaext\right)$. Then, there exists a sequence $\left( R_m \right)_{m \in \mathbb{N}}$ with $\lim_{m \to \infty} R_m =\infty$ and a sequence $ \left( g_{k,R_m} \right)_{m \in \mathbb{N}}$ with $g_{k,R_m} \in  H^1(\Omega^*)$, such that
\begin{itemize}
\item $g_{k,R_m} \circ \phi:\Omega^\ext \to \mathbb{R}$, restricted to $\Omegaext_{R_m}$ is a linear combination of eigenfunctions associated with the eigenvalues $\sigma_j^\Neu(\Omegaext_{R_m})$ for indices $j$ such that $\xi_j\left(\Omegaext\right) = \xi_k\left(\Omegaext\right)$; in particular, if $\xi_k\left(\Omegaext\right)$ is simple, it is an eigenfunction associated with $\sigma_k^\Neu(\Omegaext_{R_m})$,
\item $\lim_{m \to \infty} \| g_{k,R_m}  -f_k \|_{H^1(\Omega^*)} = 0$.
\end{itemize}
\end{corollary}

\begin{proof}
For $k=1$ the assertion is immediate, since both the truncated and limiting first eigenfunctions are constant. For $k>1$, let $q\geq k$ be the largest index such that $\xi_q\left(\Omegaext\right)=\xi_k\left(\Omegaext\right)$. Let $u_{j,R}$, $j=2,\dots,q$, be eigenfunctions associated with $\sigma_j^\Neu\left(\Omegaext_R\right)$, chosen so that
\[
\int_{\Omegaext_R}\mydotp{\nabla u_{i,R},\nabla u_{j,R}}\,\mathrm{d}x=\delta_{i,j},\qquad 2\leq i,j\leq q.
\]
Set $v_{1,R}:=u_{1,R}$ with $\|u_{1,R}\|_{L^2(\partial\Omega)}=1$, and $v_{j,R}:=\sqrt{\sigma_j^\Neu(\Omegaext_R)}\,u_{j,R}$ for $2\leq j\leq q$. The $v_{j,R}$ are orthonormal on $\partial\Omega$, so the variational principle for sums of eigenvalues gives
\[
\sum_{j=1}^q\sigma_j^\Neu(\Omegaext_{R/2})
\leq \sum_{j=1}^q\sigma_j^\Neu(\Omegaext_R)
-\sum_{j=2}^q\sigma_j^\Neu(\Omegaext_R)\int_{B_R\setminus B_{R/2}}|\nabla u_{j,R}|^2\,\mathrm{d}x.
\]
Rearranging the preceding inequality and using
Theorem~\ref{theo:convEVNeum}, we obtain
\[
0\leq
\sum_{j=2}^q \sigma_j^\Neu(\Omegaext_R)
\int_{B_R\setminus B_{R/2}}|\nabla u_{j,R}|^2\,\mathrm{d}x
\leq
\sum_{j=1}^q
\left(
\sigma_j^\Neu(\Omegaext_R)
-
\sigma_j^\Neu(\Omegaext_{R/2})
\right)
\longrightarrow0.
\]
Since $\sigma_j^\Neu(\Omegaext_R)\to\xi_j(\Omegaext)>0$ for
$2\leq j\leq q$, it follows that
\[
\int_{B_R\setminus B_{R/2}}|\nabla u_{j,R}|^2\,\mathrm{d}x
\longrightarrow0,
\qquad 2\leq j\leq q.
\]
Applying \cite[Lemma 4.5]{PostAnne2021} with radius $\delta=1/R$
to the pullbacks $f_{j,R}:=u_{j,R}\circ\phi$, we obtain harmonic
extensions $\widetilde f_{j,R}\in H^1(\Omega^*)$ such that
\[
\int_{B_{1/R}}|\nabla\widetilde f_{j,R}|^2\,\mathrm{d}x
\longrightarrow0.
\]
Hence, for every index $j$ with $\xi_j(\Omegaext)=\xi_k(\Omegaext)$,
\[
\int_{\Omega^*}|\nabla\widetilde f_{j,R}|^2\,\mathrm{d}x
\longrightarrow1,
\qquad
\int_{\partial\Omega^*}|\phi'|\,|\widetilde f_{j,R}|^2\,\mathrm{d}S
\longrightarrow\frac{1}{\xi_k(\Omegaext)}.
\]
For every $v\in H^1(\Omega^*)$, its restriction to
$\Omega^*\setminus B_{1/R}$ is an admissible test function for the
Neumann problem. Hence the weak eigenvalue equation, together with
the preceding estimate, passes directly to the limit. As in the proof
of Corollary~\ref{coro:convEF}, after passing to a common subsequence,
the $\widetilde f_{j,R}$ converge strongly in $H^1(\Omega^*)$, for all
indices $j$ such that $\xi_j(\Omegaext)=\xi_k(\Omegaext)$, to
eigenfunctions associated with $\xi_k(\Omegaext)$. Their boundary
traces remain pairwise orthogonal with respect to the weighted boundary
inner product, and hence the limiting eigenfunctions form a basis of
the eigenspace corresponding to $\xi_k(\Omegaext)$. Writing the
prescribed eigenfunction $f_k$ in this basis and taking the same linear
combination of the $\widetilde f_{j,R}$ gives the required sequence
$g_{k,R_m}$. Since $\widetilde f_{j,R}=f_{j,R}$ on
$\Omega^*\setminus B_{1/R}$, the restriction of
$\widetilde f_{j,R}\circ\phi$ to $\Omegaext_R$ equals $u_{j,R}$.
\end{proof}

\subsubsection{Conformal mapping and  Helmholtz equation}\label{sec:equivhelmn2}

In this subsection, we establish the equivalence between the formulation with conformal mappings and the formulation introduced in \S\ref{subsec:Greben}, based on the Helmholtz equation. Theorem ~\ref{theo:lambda*} mirrors the convergence result in Theorem~\ref{theo:sigmaD*} for the Dirichlet problem.

\begin{theorem}\label{theo:lambda*}
Let $\Omega \subset \mathbb{R}^2$ be a bounded open set with Lipschitz boundary and with connected $\Omegaext$,  and such that $\Omega$ contains the origin. 
Then, for any $k \in \mathbb{N}$,
\[
\lim_{\Lambda \searrow 0} \mu_k\left(\Omegaext,\Lambda\right) = \xi_k\left(\Omegaext\right),
\]
where $\mu_k\left(\Omegaext,\Lambda\right)$ is defined as in \S\ref{subsec:Greben}.
\end{theorem}
\begin{proof}
Let $u_{k,\Lambda} \in H^1\left(\Omegaext\right)$ be an eigenfunction associated with $\mu_k\left(\Omegaext,\Lambda\right)$. Define
\[
f_{k,\Lambda} : \Omega^*\to \mathbb{R}, \qquad f_{k,\Lambda}(x) = 
\begin{cases} 
u_{k,\Lambda} \left( \phi(x) \right) \qquad&\text{ if } x \neq 0,\\
0, \qquad&\text{ if } x=0,
\end{cases}
\]
where $\Omega^*$ and $\phi$ are defined as in \S\ref{subsec:conformal}. Since $u_{k,\Lambda}$ vanishes at infinity for any $\Lambda>0$,  it follows that $f_{k,\Lambda} \in H^1(\Omega^*)$. For any given $\Lambda>0$ we define the $j$-dimensional subspaces
\[
M(j) := \operatorname{Span} \left\{ f_{k,\Lambda} | 1 \leq k \leq j  \right\} \subset H^1(\Omega^*)
\] 
and
\[
N(j) := \operatorname{Span} \left\{ u_{k,\Lambda} | 1 \leq k \leq j \right\} \subset  H^1\left(\Omegaext\right). 
\]
Since $\frac{\int_{\Omega^* } | \nabla f_{k,\Lambda}  |^2 \, \mathrm{d}x}{\int_{\partial\Omega^*} |\phi'| \cdot |  f_{k,\Lambda}  |^2 \, \mathrm{d}S}  = \frac{\int_{\Omegaext} | \nabla u_{k,\Lambda}  |^2 \, \mathrm{d}x}{\int_{\partial\Omega}  |  u_{k,\Lambda}  |^2 \, \mathrm{d}S}$ and $\Lambda>0$, we have
\[
\begin{split}
\xi_j\left(\Omegaext\right) &\leq \sup_{0 \neq u \in M(j)} \frac{\int_{\Omega^* } | \nabla u  |^2 \, \mathrm{d}x}{\int_{\partial\Omega^*} |\phi'| \cdot |  u  |^2 \, \mathrm{d}S}\\   
&\leq \sup_{0 \neq u \in N(j)} \frac{\Lambda^2 \int_{ \Omegaext} |u|^2 \, \mathrm{d}x + \int_{\Omegaext } | \nabla u  |^2 \, \mathrm{d}x}{\int_{\partial\Omega}  |  u  |^2 \, \mathrm{d}S} = \mu_j\left(\Omegaext,\Lambda\right).
\end{split}
\]
For the reverse inequality, let $f_{k} \in H^1(\Omega^*)$ be an eigenfunction associated with $\xi_k\left(\Omegaext\right)$, normalised so that $\int_{\partial\Omega^*} |\phi'| f_j f_k \, \mathrm{d}S=\delta_{j,k}$. We construct a cut-off function by
\[
g_{k,R}:\Omegaext \to \R,\quad g_{k,R}(x) = \begin{cases}
f_k\left( \phi^{-1}(x) \right) \quad&\text{ if } x \in \Omegaext_R,\\
f_k\left( \phi^{-1}(x) \right) \left( \frac{\log(R^2)- \log(|x|)}{\log(R^2)- \log(R)} \right) \quad&\text{ if } R < |x|\leq R^2, \\
0 \quad&\text{ if } |x|>R^2.
\end{cases}
\]
Now define the subspaces
\[
L(j,R):=\operatorname{Span} \left\{g_{k,R}, 1 \leq k \leq j  \right\} \subset  H^1\left(\Omegaext\right)
\]
and
\[
K(j):=\operatorname{Span} \left\{ f_{k} | 1 \leq k \leq j \right\} \subset  H^1(\Omega^*).
\]
Since $L(j,R)$ is independent of $\Lambda$, we may pass to the limit and obtain
\begin{equation}\label{eq:pto0}
\mu_j\left(\Omegaext,\Lambda\right) \leq  \sup_{0 \neq u \in L(j,R)} \frac{\Lambda^2 \int_{ \Omegaext} |u|^2 \, \mathrm{d}x + \int_{\Omegaext } | \nabla u  |^2 \, \mathrm{d}x}{\int_{\partial\Omega}  |  u  |^2 \, \mathrm{d}S} \xrightarrow{\Lambda \searrow 0} \sup_{0 \neq u \in L(j,R)} \frac{ \int_{\Omegaext } | \nabla u  |^2 \, \mathrm{d}x}{\int_{\partial\Omega}  |  u  |^2 \, \mathrm{d}S}. 
\end{equation}
Let $u = \sum_{k=1}^j c_k g_{k,R}$ and $v = \sum_{k=1}^j c_k f_{k}$. For any $R > R_0(\Omega)$, we have $\int_{\partial\Omega}  |  u  |^2 \, \mathrm{d}S = \int_{\partial\Omega^*} |\phi'| \cdot  |  v  |^2 \, \mathrm{d}S$ by construction. Moreover, using that
\[
\begin{split}
\int_{B_{R^2}\setminus B_R } \left| \nabla  \left( \frac{\log(R^2)- \log(|x|)}{\log(R^2)- \log(R)} \right) \right|^2 \, \mathrm{d}x &= \int_{B_{R^2}\setminus B_R } \frac{1}{|x|^2 \log(R)^2} \, \mathrm{d}x\\
&= \int_{R}^{R^2} \frac{1}{r \log(R)^2} \, \mathrm{d}x = \frac{1}{\log(R)},
\end{split}
\]
a straightforward computation yields
\[
\lim_{R \to \infty }\int_{\Omegaext } \mydotp{ \nabla g_{k,R}, \nabla g_{i,R} } \, \mathrm{d}x = \int_{\Omegaext } \mydotp{ \nabla (f_k \circ \phi^{-1}), \nabla (f_i \circ \phi^{-1}) } \, \mathrm{d}x = \int_{\Omega^*} \mydotp{ \nabla f_k, \nabla f_i } \, \mathrm{d}x.
\]
Hence, 
\[
\lim_{R \to \infty }\int_{\Omegaext } | \nabla u  |^2 \, \mathrm{d}x = \int_{\Omega^* } | \nabla v  |^2 \, \mathrm{d}x.
\] 
Taking the limit in \eqref{eq:pto0}, we conclude
\[
\lim_{\Lambda \searrow 0} \mu_j\left(\Omegaext,\Lambda\right) \leq\sup_{0 \neq u \in K(j)} \frac{ \int_{\Omega^* } | \nabla u  |^2 \, \mathrm{d}x}{\int_{\partial\Omega^*} |\phi'| \cdot |  u  |^2 \, \mathrm{d}S}= \xi_j\left(\Omegaext\right).
\]
\end{proof}

\begin{corollary}\label{coro:convEF2}
Suppose $f_k \in H^1(\Omega^*)$ is an eigenfunction associated with $\xi_k\left(\Omegaext\right)$. Then, there exists a sequence $\left( \Lambda_m \right)_{m \in \mathbb{N}} \subset (0,1) $ with $\lim_{m \to \infty} \Lambda_m =0$ and a sequence $ \left( g_{k,\Lambda_m} \right)_{m \in \mathbb{N}} \subset H^1\left(\Omegaext\right)$, such that
\begin{itemize}
\item $g_{k,\Lambda_m}$  is a linear combination of eigenfunctions associated with the eigenvalues $\mu_j(\Omegaext,\Lambda_m)$ for indices $j$ such that $\xi_j\left(\Omegaext\right) = \xi_k\left(\Omegaext\right)$ (in particular, if $\xi_k\left(\Omegaext\right)$ is simple, then $g_{k,\Lambda_m}$ is an eigenfunction associated with $\mu_k(\Omegaext,\Lambda_m)$) and $\lim_{m \to \infty} \mu_k(\Omegaext,\Lambda_m) = \xi_k\left(\Omegaext\right)$,
\item $\lim_{m \to \infty} \| g_{k,\Lambda_m} \circ \phi -f_k \|_{H^1(\Omega^*)} = 0$.
\end{itemize}
\end{corollary}

\begin{proof}
The case $k=1$ is handled exactly as in the last part of the proof of Corollary~\ref{coro:convEF}, using Theorem~\ref{theo:lambda*} in place of Theorem~\ref{theo:sigmaD*}. Hence assume $k\geq2$.

Without restriction we may assume  $\| \nabla f_{k} \|_{L^2(\Omega^*)} =1$. Let $u_{k,\Lambda} \in H^1\left(\Omegaext\right)$ denote the eigenfunctions associated with $\mu_k\left(\Omegaext,\Lambda\right)$ such that $\| \nabla u_{k,\Lambda} \|_{L^2\left(\Omegaext\right)} =1$ and define
\[
f_{k,\Lambda} : \Omega^* \to \mathbb{R}, \qquad f_{k,\Lambda} (x) = \begin{cases}
u_{k,\Lambda} \left( \phi(x) \right) \qquad&\text{ if } x \neq 0,\\
0 \qquad&\text{ if } x =0.
\end{cases}
\]
Then,
\[
\frac{1}{\int_{\partial\Omega^*} |\phi'| \cdot |f_{k,\Lambda}|^2 \, \mathrm{d}S} = \frac{1}{\int_{\partial\Omega}  |u_{k,\Lambda}|^2 \, \mathrm{d}S} \leq \mu_k\left(\Omegaext,\Lambda\right).
\]
For $k\geq2$, the claim follows by the same direct compactness and linear-combination argument as in Corollary~\ref{lemma5}, after conformal pullback. The only additional point is the $H^1(\Omega^*)$ bound: Poincar\'e's inequality controls each pullback modulo its average, and boundary orthogonality to the first Helmholtz eigenfunction (whose pullback converges to a non-zero constant by the $k=1$ case) controls these averages. This also bounds the boundary traces, and hence $\Lambda\|u_{k,\Lambda}\|_{L^2\left(\Omegaext\right)}$ by the eigenvalue identity used in Corollary~\ref{lemma5}. The logarithmic cutoff from the proof of Theorem~\ref{theo:lambda*} then permits passage to the limiting weak equation. This proves the claim.
\end{proof}

\subsubsection{Conformal mapping and layer potentials}

In \S\ref{subsec:dimtwo}, we introduced a formulation of the exterior Steklov eigenvalue problem in two dimensions using boundary layer potentials. Using Theorem~\ref{theo:layer2d}, we can readily verify the equivalence with the conformal mapping formulation.

\begin{theorem}\label{thm:CMLP2}
Let $\Omega \subset \mathbb{R}^2$ be a bounded open set with Lipschitz boundary and with connected $\Omegaext$,  and such that $\Omega$ contains the origin. 
Let $\tau_k(\Omegaext)$ be defined as in Theorem \ref{theo:layer2d}. Then, $\tau_k(\Omegaext) = \xi_k\left(\Omegaext\right)$ for any $k \in \mathbb{N}$. 

Let $u_k \in H^{\frac{1}{2}}(\partial\Omega)$ be an eigenfunction associated with $\tau_k(\Omegaext)$, and let $U_k$ be its harmonic extension satisfying \eqref{eq:inftyr}.  Then, $U_k \circ \phi$ coincides, up to a constant factor, with an eigenfunction associated with $\xi_k\left(\Omegaext\right)$.
\end{theorem}

\begin{proof}
As $\tau_1  = \xi_1  = 0$ with constant eigenfunctions, we only have to consider $k \geq 2$.

Let $w_k \in H^1(\Omega^*)$ be an eigenfunction corresponding to $\xi_k  \neq 0$. Define
\[
f_k: \Omega^\ext \to \mathbb{R}, \quad f_k(x):= ( w_k \circ \phi^{-1})(x).
\]
Since $\lim_{|x| \to \infty} f_k(x) = \lim_{x \to 0} w_k(x) = w_k(0)$, the mean value theorem yields
\[
\left| f_k(x)-w_k(0) \right| = \left|w_k \left( \frac{x}{|x|^2}\right)-w_k(0)\right| =  \frac{| \nabla w_k(x_0)|}{|x|}
\]
for some $x_0$ between $0$ and $\frac{x}{|x|^2}$. By regularity of $w_k$ we know that $|\nabla w_k(x_0)|$ is bounded, so
\[
\left| f_k(x)-w_k(0) \right| = O\left( \frac{1}{|x|} \right) \quad \text{ as } |x| \to \infty.
\]
Therefore, $f_k$ satisfies the far field condition \eqref{eq:farfield} and by Theorem~\ref{theo:layer2d}, the pair $(f_k, \xi_k)$ is also a solution of~\eqref{eq:ext_2D_4}. 

Conversely, let $u_k \in  \left\{ u \in H^{\frac{1}{2}}(\partial\Omega) : \myscal{ u, 1 }_{\partial\Omega}=0 \right\}$ be an eigenfunction corresponding to $\tau_k  \neq 0$. Define  $(u_k)_\infty \in \mathbb{R}$ as in~\eqref{eq:u_0} and let $U_k: \Omega^\ext \to \mathbb{R}$  be the harmonic extension of  $u_k$, as given in Theorem~\ref{theo:layer2d}. Then, $W_k \in H^1(\Omega^*)$, given by 
\[
W_k(x) = \begin{cases}
U_k ( \phi(x))\qquad&\text{ for } x \neq 0,\\
(u_k)_\infty\qquad&\text{ for } x = 0,
\end{cases} 
\]
is a solution of \eqref{eq:problemtransform}. Note that $W_k  \in H^1(\Omega^*)$, even though $U_k$ is not necessarily in $H^1(\Omegaext)$.
\end{proof}

\section{Properties of exterior Steklov eigenvalues and eigenfunctions}\label{sec:properties}

In the following, we no longer distinguish between the different formulations of the Steklov problem as they are equivalent,  and we denote the Steklov eigenvalue by $\sigma$.

\subsection{Basic properties}\label{subsec:basic}

\begin{proof}[Proof of Theorem \ref{thm:basic}]
Let $n \geq 3$, let $u$ be a solution of  \eqref{eq:PDEext} and  let $f$ be the restriction of $u$ to $\partial \Omega$. It follows from \cite[Section 5]{auchmuty2014representations} that there is a unique $U\in E^1\left(\Omegaext\right)$ such that $U=f$ on $\partial \Omega$. Moreover, $U=O(|x|^{2-n})$ by \cite[Proposition 12]{xiong2023sharp} (see also Proposition \ref{prop:xiongdecay}). It follows that $U-u$ is harmonic, vanishing on $\partial\Omega$ and decaying at $\infty$, whence $u=U\in E^1\left(\Omegaext\right)$. So, any solution of \eqref{eq:PDEext} is in $E^1\left(\Omegaext\right)$. It then follows by integration by parts against a test function that $u$ is a solution of the finite energy weak Steklov eigenvalue equation \eqref{WeakFormulationExterior}. On the other hand, Theorem \ref{thm:spectral} provides a sequence of weak Steklov eigenfunction $u_k\in E^1\left(\Omegaext\right)$ whose restriction to $\partial \Omega$ form an orthogonal basis of $L^2(\partial \Omega)$. It follows from \cite[Proposition 12]{xiong2023sharp}  that these functions satisfy \eqref{eq:inftyr}, and hence they are solutions of \eqref{eq:PDEext}.

Let $n = 2$ and let $u$ be a solution of  \eqref{eq:PDEext}. Because $u$ is harmonic and bounded, it converges to a finite value at infinity, denoted by $u_\infty$, see \cite[Theorem 4.9]{ABR01}. With $\phi$ as in  \eqref{eq:phi}, we consider the function $v:\Omega^*  \to\R$, $v = u\circ\phi$ in $\Omega^*\setminus \{ 0 \}$ and $v(0)=u_\infty$. Then, $v$ is harmonic and satisfies \eqref{eq:problemtransform} and so it is a Steklov eigenfunction in the conformal sense.  Similarly, any solution of \eqref{eq:problemtransform} will give a solution of \eqref{eq:PDEext} once it is inverted through $\phi$.
\end{proof}

\subsubsection{Courant-type bound for the nodal count}\label{subsec:Courant}

\begin{proof}[Proof of Theorem \ref{thm:Courant}]
For $n \geq 3$, we can use the finite-energy approach and proceed verbatim as in the interior case, see, e.g., \cite{kuttler1969inequality}.  Let us briefly recall this argument. Indeed, assume that an eigenfunction $u_k$ associated with an eigenvalue $\sigma_k$ has  $K \ge k+1$ nodal domains $D_i$, $i=1,\dots, K$. Choose a non-zero linear combination of the restrictions $\left.u_k\right|_{D_i}$, $i=1,\dots,k$, whose trace is orthogonal to the traces of the first $k-1$ exterior Steklov eigenfunctions on $\partial\Omega$. The disjointness of the nodal domains keeps the Rayleigh quotient of this linear combination equal to $\sigma_k$, and the min--max principle then forces it to be an eigenfunction. However, the restriction of this linear combination to $D_{k+1}$ is identically zero, and we get  a contradiction with the unique continuation principle.

For $n=2$, if we show that the number of nodal domains in $\Omegaext$ is preserved under the conformal transformation $\phi$, defined in \eqref{eq:phi}, then the claim follows from Courant's nodal domain theorem for the weighted Steklov problem on a bounded domain $\Omega^*$ (see, for instance, \cite[Section 2.1]{KKP14}).
Let $\Omega \subset \mathbb{R}^2$ be a bounded Lipschitz domain, and let $u$ be a solution of \eqref{eq:PDEext}. Since $\phi$ is a homeomorphism between $\Omegaext$ and $\Omega^* \setminus \{0\}$, it preserves connected components,  and hence the number of nodal domains of $u \in \Omegaext$ equals the number of nodal domains of $w:= u \circ \phi$ in $\Omega^* \setminus \{0\}$. Since $u$ is bounded, it has a limit at infinity, and therefore $w$ extends continuously to zero. If $w(0)=0$, then the origin is not part of any nodal domain,  so extending $w$ to the origin does not merge or split any connected components. If $w(0) \neq 0$, there exists an $R>0$ such that $u(x) \neq 0$ for any $|x|>R$. Then, $w(x) \neq 0$ for any $|x|<\frac{1}{R}$. So again, extending $w$ to the origin does not merge or split any connected components.
\end{proof}

\begin{proof}[Proof of Corollary \ref{corol:simple}]
Suppose, for contradiction, that there are two linearly independent eigenfunctions $u_1, u_2$ associated with $\sigma_1\left(\Omegaext\right)$. By Theorem \ref{thm:Courant}, we may choose them  so that $u_1(x), u_2(x) \ge 0$ for any $x \in \Omegaext$. Take some  $x_0 \in \Omegaext$, and choose  $c\le 0$ such that $u_1(x_0)+c u_2(x_0)=0$. Then, $w(x):=u_1(x)+c u_2(x)$ is a nontrivial eigenfunction associated with $\sigma_1\left(\Omegaext\right)$, and $w(x_0)=0$. By Theorem \ref{thm:Courant} and the maximum principle, $w$ cannot vanish at an interior point unless it is identically zero. This contradiction implies that the first eigenvalue is simple.
\end{proof}

\subsubsection{Decay rate of the eigenfunctions}

First, we  recall a result from \cite[Proposition 12]{xiong2023sharp}, which describes a decay rate of the eigenfunctions. 
\begin{proposition}\label{prop:xiongdecay}
Suppose $\Omega \subset \mathbb{R}^n$, $n \geq 3$, is a bounded domain with Lipschitz boundary and connected $\Omegaext$. Any harmonic function $u \in E^1\left(\Omegaext\right)$ satisfies
\[
u = O\left( |x|^{2-n} \right) \quad \text{as} \quad |x| \to \infty.
\]
In particular, any eigenfunction associated with $\sigma_k\left(\Omegaext\right)$, $k \in \mathbb{N}$, behaves asymptotically as $O\left( |x|^{2-n} \right)$.
\end{proposition}

To get a more precise understanding of the decay rate, we show below that the first eigenfunction cannot decay  faster than $|x|^{2-n}$.

\begin{proposition}\label{pro:decay}
Suppose $\Omega \subset \mathbb{R}^n$, $n \geq 3$, is a bounded domain with Lipschitz boundary and connected $\Omegaext$, and let $u_1 \in E^1\left(\Omegaext\right)$ be an eigenfunction associated with $\sigma_1\left(\Omegaext\right)$. Then, there exists a $c \in \mathbb{R} \setminus \{ 0 \}$ such that
\[
\lim_{|x| \to \infty} \frac{u_1(x)}{|x|^{2-n}} = c.
\]
For $n\in \{3,4\}$ this shows, in particular, that the first eigenfunction is not in $L^2\left(\Omegaext\right)$.
\end{proposition}

\begin{proof}
Outside a ball, $u_1$ can be expanded in terms of Steklov eigenfunctions of the ball, whose explicit form is given in \S\ref{examp:balln3}. The asymptotics then follow from the leading term of this expansion.

Without loss of generality, assume that $\Omega \subset B_1$. With the notation from \S\ref{examp:balln3}, we define for $r \geq 1$, $\theta \in \mathbb{S}^{n-1}$, and $k \in \mathbb{N}_0$,  $1 \leq i \leq d_{n, k}$, 
\[
u_{(k,i)}(r,\theta) =  r^{2-n-k} Y_{k,i}(\theta).
\] 
Then ${u}_{(k,i)}\in E^1\left(B_1^\ext\right)$ are Steklov eigenfunctions associated to $\sigma_{(k)}\left(B_1^\ext\right) = n+k-2 $, and their restrictions to $\partial B_1$ form an orthonormal basis in $L^2\left(\partial B_1\right)$, see Theorem~\ref{thm:basic}. Let $f := u_1|_{\partial B_1}$. In $L^2\left(\partial B_1\right)$ we have the expansion
\begin{equation}\label{eq:seriesexpf}
f(s) = \sum_{k=0}^\infty \sum_{i=1}^{d_{n,k}} \myscal{ f, {u}_{(k,i)} }_{L^2\left(\partial B_1\right)} {u}_{(k,i)}(s), \qquad s \in \partial B_1.
\end{equation}
In the first step, we show uniform convergence of \eqref{eq:seriesexpf}. Since $u_1$ is harmonic, $u_1 \in C^\infty\left(\Omegaext\right)$ and hence $f \in C^\infty\left(\partial B_1\right)$. Using the fact that the spherical harmonics $Y_{k,i}$ are eigenfunctions of the Laplace--Beltrami operator $-\Delta_{\mathbb S}$  with eigenvalue $k(k+n-2)$, 
one obtains 
\[
\begin{split}
\left| \myscal{ f, {u}_{(k,i)}}_{L^2\left(\partial B_1\right)}\right| &= \frac{1}{k(k +n-2)} \left| \myscal{ f, -\Delta_\mathbb{S} Y_{k,i}}_{L^2\left(\partial B_1\right)}\right|\\
&=\frac{1 }{k(k +n-2)} \left| \myscal{-\Delta_\mathbb{S} f,  Y_{k,i}}_{L^2\left(\partial B_1\right)}\right|.
\end{split}
\]
Repeating this process $m$ times and using the Cauchy--Schwarz inequality,
\[
\begin{split}
\left| \myscal{ f, {u}_{(k,i)}}_{L^2\left(\partial B_1\right)}\right| &= \frac{1}{k^m(k +n-2)^m} \left| \myscal{ (-\Delta_\mathbb{S})^m f,  Y_{k,i}}_{L^2\left(\partial B_1\right)}\right|\\ 
&\leq \frac{\left\| (-\Delta_\mathbb{S})^m f \right\|_{L^2\left(\partial B_1\right)}}{k^m(k +n-2)^m}.
\end{split}
\]
Because $f$ is smooth, $\left\|  (-\Delta_\mathbb{S})^m f\right\|_{L^2\left(\partial B_1\right)}$ is bounded. Choosing $m$ large enough and noting that $d_{n,k} = O(k^{n-2})$, we see that \eqref{eq:seriesexpf} converges uniformly. 

By uniqueness of the harmonic extension, the same coefficients give the expansion of $u_1$ in $\Omegaext$,
\[
u_1(x) =  \sum_{k=0}^\infty \sum_{i=1}^{d_{n,k}} \myscal{ f, {u}_{(k,i)} }_{L^2\left(\partial B_1\right)} {u}_{(k,i)}(x), \qquad x \in \Omegaext.
\]
Since ${u}_{(k,i)}(x) = O(|x|^{2-n-k})$, it remains to show that $ \myscal{ f, {u}_{0,1} }_{L^2\left(\partial B_1\right)} \neq 0$. Since spherical harmonics of degree $0$ are constants, this is equivalent to $\int_{\partial B_1} f \, \mathrm{d}S \neq 0$. 

If $\int_{\partial B_1} f \, \mathrm{d}S = 0$, then $f=0$ on $\partial B_1$ since $u_1$ does not change sign by Theorem \ref{thm:Courant}. Then,  $u_1 \equiv 0$ in $\Omegaext$.  Therefore, the coefficient of degree $k=0$ in the spherical expansion of $u_1$ has to be nonzero.
\end{proof}

\subsection{Lower bounds for the first eigenvalue}\label{sec:lowerbounds}

The only established lower bound so far, to the best of our knowledge, is given in \cite[Theorem 1]{xiong2023sharp}. In this section we prove a new lower bound depending on the mean curvatures of the boundary. In Figure~\ref{fig:LowerBounds}, we compare the bound of \cite{xiong2023sharp} with our result for different spheroids, showing that neither inequality implies the other.

\begin{proof}[Proof of Theorem \ref{prop:low}]
We adapt the method from \cite[Proof of Theorem 1.1]{kovavrik2017p}, see also \cite{krejcirik2016optimisation,krejvcivrik2020optimisation}, to derive a lower bound. Since $\Omega$ is convex, the mapping
\[
\Psi: \partial\Omega \times (0, \infty) \to \Omega^\ext, \qquad (s,t) \mapsto s- t \nu(s)
\]
is bijective and locally bi-Lipschitz.  For a discussion of the properties and further background on $\Psi$, we refer to \cite[Chapter 3]{kovavrik2017p} and the references therein (note that we will switch the sign in front of $\kappa_j(s)$ in~\eqref{eq:varphi} compared to \cite{kovavrik2017p} because of the orientation of $\nu$).  Any $u \in E^1\left(\Omegaext\right)$ can be approximated by functions with compact support in $\Omegaext$, allowing us to perform a change of variables, yielding
\[
\int_{\Omega^\ext} |\nabla u|^2 \, \mathrm{d}x = \int_{\partial\Omega \times (0,\infty)} |(\nabla u) \circ \Psi(s,t)|^2 \operatorname{det}(D \Psi) \, \mathrm{d}S_s \, \mathrm{d}t,
\]
where the Jacobian determinant is given by
\begin{equation}\label{eq:varphi}
\zeta(s,t):=\operatorname{det}(D \Psi)  =  \prod_{j=1}^{n-1} \left( 1 + \kappa_j(s) t \right).
\end{equation} 
Define $w(s,t) := u(\Psi(s,t))$. By the chain rule, $\partial_{t} w(s,t) = -\mydotp{ (\nabla u) (\Psi(s,t)), \nu(s) }$ and hence
\[
|(\nabla u) \circ \Psi(s,t)| \geq |\mydotp{ (\nabla u) \circ \Psi(s,t), \nu(s) }| = |\partial_{t} w(s,t)|.
\]
It follows that
\[
\begin{split}
\int_{\Omega^\ext} |\nabla u|^2 \, \mathrm{d}x &= \int_{\partial\Omega \times (0,\infty)} |(\nabla u) \circ \Psi(s,t)|^2 \zeta(s,t) \, \mathrm{d}S_s \, \mathrm{d}t \\
&\geq \int_{\partial\Omega \times (0,\infty)} |\partial_t w(s,t)|^2  \zeta(s,t) \, \mathrm{d}S_s \, \mathrm{d}t.
\end{split}
\]
Thus, we obtain the lower bound
\[
\begin{split}
\sigma_1\left(\Omegaext\right) &=
\inf_{ 0 \neq u \in E^1\left(\Omegaext\right)} \frac{\int_{\Omega^\ext} |\nabla u|^2 \, \mathrm{d}x}{\int_{\partial\Omega} |u|^2 \, \mathrm{d}S}\\ 
&\geq  \inf_{\substack{0 \neq u \in E^1\left(\Omegaext\right)\\w=u \circ \Psi}} \frac{\int_{\partial\Omega \times (0,\infty)} |\partial_t w(s,t)|^2  \zeta(s,t) \, \mathrm{d}S_s \, \mathrm{d}t}{\int_{\partial\Omega} |w(s,0)|^2 \, \mathrm{d}S_s}.
\end{split}
\]
Since $w=u \circ \Psi$ decays at infinity for any $u \in E^1\left(\Omegaext\right)$, we define, for each $s \in \partial\Omega$,
\begin{equation}\label{eq:kappalowerbound}
K\left(s,\Omegaext\right)  :=  \inf_{\substack{f\in {H}_\mathrm{loc}^1\left( (0, \infty) \right)\\\lim_{x \to \infty}f(x)=0}} \frac{\int_0^\infty |f'(t)|^2 \zeta(s,t) \, \mathrm{d}t}{|f(0)|^2},
\end{equation}
which gives
\begin{equation}\label{eq:ineqKs}
\sigma_1\left(\Omegaext\right)\geq \inf_{\substack{0 \neq u \in E^1\left(\Omegaext\right)\\w=u \circ \Psi}} \frac{\int_{\partial\Omega} K(s, \Omega^\ext) |w(s,0)|^2 \, \mathrm{d}S_s}{\int_{\partial\Omega} |w(s,0)|^2 \, \mathrm{d}S_s} \geq \inf_{s \in \partial\Omega} K\left(s,\Omegaext\right) .
\end{equation}
Fix $s\in\partial\Omega$. We first assume that $\kappa_j(s)>0$ for all $j=1,\dots,n-1$. %; vanishing curvatures will be treated at the end.}
To show that the infimum in \eqref{eq:kappalowerbound} is attained, consider a minimising sequence $(f_k)_{k \in \mathbb{N}} \subset H^1_\mathrm{loc}((0,\infty))$ with $f_k(0)=1$, $f_k(x) \to 0$ as $x \to \infty$ and $\lim_{k \to \infty }\int_0^\infty |f_k'(t)|^2 \zeta(s,t) \, \mathrm{d}t = K\left(s,\Omegaext\right) $. It is straightforward to show that, for each  $R>0$, the sequence is bounded in $H^1((0,R))$. Hence, by standard compactness results and a diagonal argument, there exists a subsequence, again denoted by $f_k$, that converges weakly in $H^1((0,R))$ and strongly in $L^2((0,R))$ to a limit function $f \in H^1((0,R))$. Specifically, $f(0)=1$.

It remains to show that $f$ decays at infinity. Since $n \geq 3$ and $\kappa_j(s)>0$ for all $j$, there exists $D_1(s)>0$ with $ \zeta(s,t) \geq D_1(s) t^{n-1}$. Thus,
\[
\int_x^\infty  \frac{1}{\zeta(s,t)} \, \mathrm{d}t < \frac{x^{2-n}}{(n-2)D_1(s)}.
\]
Then, there exists a constant $D_2(s)$ such that for any $x \in (0,\infty)$,
\[
\begin{split}
|f_k(x)| &= \left| \int_x^\infty f_k'(t) \, \mathrm{d}t \right| \leq  \int_x^\infty | f_k'(t) | \, \mathrm{d}t\\ 
&\leq \sqrt{ \int_x^\infty | f_k'(t) |^2 \zeta(s,t) \, \mathrm{d}t} \sqrt{ \int_x^\infty \frac{1}{ \zeta(s,t)} \, \mathrm{d}t} \leq D_2(s) x^{\frac{2-n}{2}},
\end{split}
\]
where we use that the first integral in the product 
 is uniformly bounded for any minimising sequence. This uniform decay implies $f(x)\to 0$ as $x \to \infty$. Finally, by the lower semicontinuity  of the functional, $f$ minimises~\eqref{eq:kappalowerbound}. Moreover, $K\left(s,\Omegaext\right)  >0$ since $f$ is nonconstant. The associated Euler--Lagrange equation for this minimisation problem reads
\begin{equation}\label{eq:PDEw2}
\left( f'(t) \zeta(s,t) \right)' = 0 \, \text{ for } t \in (0, \infty),\quad -f'(0) = K\left(s,\Omegaext\right)  f(0), \quad f(t) \to 0 \, \text{ as } t \to \infty.
\end{equation}

Assume now that the principal curvatures $\kappa_j(s)$ are distinct. Then the general solution of the differential equation in~\eqref{eq:PDEw2} is given by
\[
\begin{split}
f(t) &=C_1 + C_2 \int_0^t \frac{1}{\zeta(s,r)} \mathrm{d}r = C_1 + C_2\int_0^t \left(\sum_{j=1}^{n-1} \frac{c_j(s)}{1+\kappa_j(s) r}\right)\, \mathrm{d}r \\
& = C_1 + C_2\sum_{j=1}^{n-1}\frac{c_j(s)  \log(1+\kappa_j(s)t)}{\kappa_j(s)} = C_1 + C_2 \log\left( \prod_{j=1}^{n-1} (1+\kappa_j(s)t)^\frac{c_j(s)}{\kappa_j(s)}\right),
\end{split}
\]
where
\begin{equation}\label{eq:cjs}
c_j(s)=\frac{\kappa_j(s)^{n-2}}{\prod_{\substack{i=1, \, i \neq j}}^{n-1}  (\kappa_j(s)-\kappa_i(s))},\qquad j=1,\dots, n-1,
\end{equation} 
are the coefficients of the partial fraction decomposition
\[
1 = \sum_{j=1}^{n-1} c_j(s) \prod_{\substack{i=1\\i \neq j}}^{n-1}  (1+\kappa_i(s) r).
\]
It is easy to check that $\sum_{j=1}^{n-1}  \frac{c_j(s)}{\kappa_j(s)} = 0$, which ensures that $f(t)$ is bounded. Thus, the asymptotic behaviour of $f(t)$ simplifies to
\[
\lim_{t \to \infty}f(t) = C_1 + C_2 \log\left( \prod_{j=1}^{n-1} \kappa_j(s)^\frac{c_j(s)}{\kappa_j(s)}\right),
\]
and the solution decays at infinity if and only if  $C_1 = - C_2\log\left( \prod_{j=1}^{n-1} \kappa_j(s)^\frac{c_j(s)}{\kappa_j(s)}\right)$. 
Because $\sum_{j=1}^{n-1}c_j(s)=1$, this yields
\[
K\left(s,\Omegaext\right)  = \frac{-f'(0)}{f(0)} = \frac{1}{\log\left( \prod_{j=1}^{n-1} \kappa_j(s)^\frac{c_j(s)}{\kappa_j(s)}\right)}= \frac{1}{\sum_{j=1}^{n-1} \frac{c_j(s)}{\kappa_j(s)} \log\left(  \kappa_j(s) \right)}.
\]
Substituting the expressions \eqref{eq:cjs} for $c_j(s)$, and taking into account \eqref{eq:LM}, we arrive at
\[
K\left(s,\Omegaext\right)   = \frac{1}{\sum_{j=1}^{n-1} \frac{\kappa_j(s)^{n-3}}{\prod_{\substack{i=1, \, i \neq j}}^{n-1}  (\kappa_j(s)-\kappa_i(s))} \log\left( \kappa_j(s) \right)} 
= (n-2)L\left(\kappa_1(s),\dots,\kappa_{n-1}(s)\right),
\]
which completes the proof when the curvatures $\kappa_j(s)$ are distinct and positive. Solving~\eqref{eq:PDEw2} directly gives
\[
K\left(s,\Omegaext\right)=\left(\int_0^\infty \zeta(s,t)^{-1}\,\mathrm{d}t\right)^{-1}.
\]
In fact the inequality
\[
K\left(s,\Omegaext\right)\geq\left(\int_0^\infty \zeta(s,t)^{-1}\,\mathrm{d}t\right)^{-1},
\]
which is the direction required in~\eqref{eq:ineqKs}, holds without any positivity assumption, with the convention $1/\infty=0$. Indeed, let $f$ be admissible in~\eqref{eq:kappalowerbound}. For every $T>0$ the Cauchy--Schwarz inequality gives
\[
|f(T)-f(0)|^2\leq\int_0^T |f'(t)|^2\zeta(s,t)\,\mathrm{d}t
\int_0^T \zeta(s,t)^{-1}\,\mathrm{d}t,
\]
and letting $T\to\infty$, so that $f(T)\to0$, yields the claim. 
The cases of coinciding or vanishing curvatures then follow, since $(n-2)L\left(\kappa_1(s),\dots,\kappa_{n-1}(s)\right)$ is given by the same integral,  where $L$ is defined by continuity.
\end{proof}

\begin{proof}[Proof of Corollary \ref{coro:lowbound}]
For distinct, positive $\alpha_1, \ldots, \alpha_{n-1}$, by \cite[Proposition 4.1 and Theorem 3.8]{pittenger1985logarithmic}, 
\[
L(\alpha_1, \ldots, \alpha_{n-1}) \geq   \sqrt[n-1]{\prod_{j=1}^{n-1}  \alpha_j}.
\]
So \eqref{eq:logmean} already implies \eqref{eq:Ksgeommean} in the case of distinct principal curvatures. Situations where multiple $\kappa_i$'s  share the same value can be treated as limiting cases of the above computation.

It remains to prove that equality holds only for balls.   Suppose equality holds in \eqref{eq:Ksgeommean}. Then both inequalities in \eqref{eq:ineqKs} must be equalities, and
\begin{equation}\label{eq:equalityescobar}
\inf_{s \in \partial \Omega} K\left(s,\Omegaext\right)  = (n-2) \inf_{s \in \partial \Omega}\prod_{j=1}^{n-1}  \kappa_j(s)^\frac{1}{n-1}.
\end{equation}
 The first inequality in \eqref{eq:ineqKs} can be an equality only if $u$ is an eigenfunction. Since the first eigenfunction does not change sign (see \S\ref{subsec:Courant}), the second inequality in  \eqref{eq:ineqKs} can only be an equality if $K\left(s,\Omegaext\right) $ is constant. By \eqref{eq:Ksgeommean} and \eqref{eq:equalityescobar}, then the geometric mean of the curvatures must also be a constant. By \cite[Theorem 2]{mulero1987compact}, this forces $\partial\Omega$ to be a sphere.
\end{proof}

We now compare our lower bound~\eqref{eq:logmean} with the one from \cite[Theorem~1]{xiong2023sharp}, considering prolate and oblate spheroids.

\begin{example}\label{examp:ellipsoid}
For $a \in (0,1)$, we consider two families of spheroids,
\[
\mathscr{p}_a := \left\{ x \in \mathbb{R}^3: \frac{x_1^2}{a^2}+\frac{x_2^2}{a^2}+x_3^2 =1 \right\} \quad \text{ and } \quad \mathscr{o}_a := \left\{ x \in \mathbb{R}^3: \frac{x_1^2}{a^2}+x_2^2+x_3^2 =1 \right\},
\]
and denote the interior and exterior, with respect to $\mathscr{p}_a$ and $\mathscr{o}_a$, domains in $\mathbb{R}^3$  by $\mathcal{P}_a$ and $\mathcal{P}^\ext_a$, and $\mathcal{O}_a$ and $\mathcal{O}^\ext_a$, respectively.

The principal curvatures of $\mathscr{p}_a$ and $\mathscr{o}_a$ are given by
\begin{align*}
\kappa_1^{\mathscr{p}_a}(x) &= \frac{a}{\left(1-(1-a^2)x_3^2 \right)^\frac{3}{2}}, \quad  &&\kappa_2^{\mathscr{p}_a}(x) = \frac{1}{a \left(1-(1-a^2)x_3^2 \right)^\frac{1}{2}}, \\
\kappa_1^{\mathscr{o}_a}(x) &= \frac{a^4}{\left(a^4   +(1-a^2)x_1^2 \right)^\frac{3}{2}}, \quad  &&\kappa_2^{\mathscr{o}_a}(x) = \frac{a^2}{\left(a^4+(1-a^2)x_1^2 \right)^\frac{1}{2}}.
\end{align*}
Since $\kappa_1^{\mathscr{p}_a}(x)$ and $\kappa_2^{\mathscr{p}_a}(x) $ become minimal when $|x_3|=0$, we have for prolate spheroids
\begin{equation}\label{eq:weprolate}
\beta\left(\mathscr{p}_a\right) =\min_{x \in  \mathscr{p}_a } L\left( \kappa_1^{\mathscr{p}_a}(x) ,\kappa_2^{\mathscr{p}_a}(x)\right)=L\left( a ,a^{-1} \right)= \frac{1-a^2  }{-2a \log(a) }
\end{equation}
and analogously  for oblate spheroids,
\begin{equation}\label{eq:weoblate}
\beta\left(\mathscr{o}_a\right) =\min_{x \in  \mathscr{o}_a } L\left( \kappa_1^{\mathscr{o}_a}(x) ,\kappa_2^{\mathscr{o}_a}(x)\right) = L\left( a , a \right) = a.
\end{equation}

We now compare our eigenvalue bound \eqref{eq:logmean} with Xiong's bound \eqref{eq:Xiongbound}.
A straightforward calculation yields
\begin{equation}\label{eq:xiongsph}
\beta_\mathrm{X}\left(\mathscr{p}_a\right)= \beta_\mathrm{X}\left(\mathscr{o}_a\right) 
= \begin{cases}
 \frac{3\sqrt{3}}{2}\frac{a}{(1+a^2)^\frac{3}{2}}\qquad&\text{if }0<a\le \frac{1}{\sqrt{2}},\\
 1\qquad&\text{if }\frac{1}{\sqrt{2}}<a<1.
\end{cases}
\end{equation}

The comparison of the bounds \eqref{eq:weprolate}--\eqref{eq:xiongsph} and the  first Steklov eigenvalues for prolate and oblate spheroids, computed numerically using the separation of variables technique \cite{grebenkov2024spectral},  is presented in Figure \ref{fig:LowerBounds}.

\begin{figure}[htb]
\centering
\includegraphics[width=0.95\textwidth]{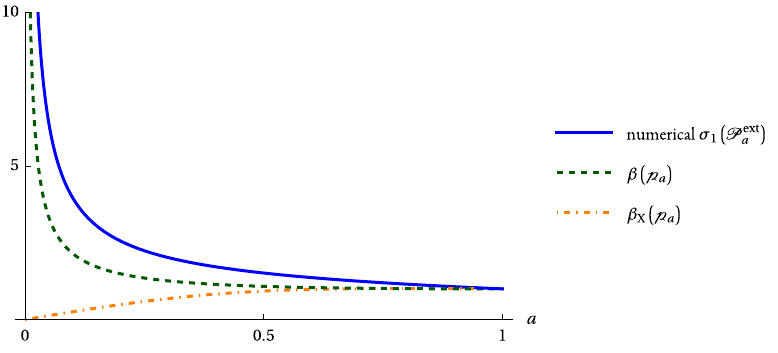}\\
\includegraphics[width=0.95\textwidth]{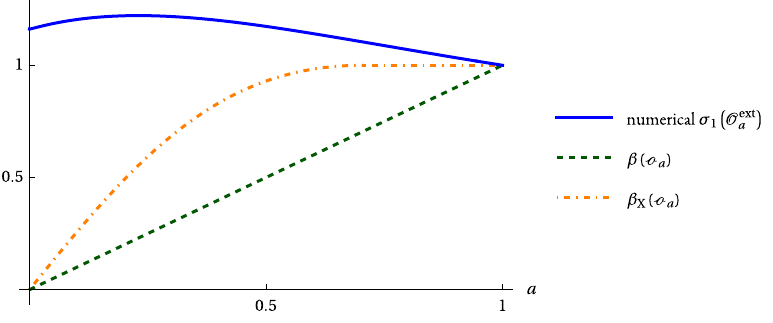}
\caption{Comparison between the numerical results, our bounds  \eqref{eq:weprolate},  \eqref{eq:weoblate}, and Xiong's bound \eqref{eq:xiongsph} for the exterior of prolate spheroids (top figure) and oblate spheroids (bottom figure).}%
\label{fig:LowerBounds}%
\end{figure}

Note that for oblate spheroids, our lower bound is weaker than the one from \cite{xiong2023sharp}, however, neither bound appears to be effective as $a \to 0$. In contrast, for prolate spheroids, the lower bound from \cite{xiong2023sharp} vanishes as $a \to 0$, whereas our lower bound  \eqref{eq:weprolate} exhibits an asymptotic behaviour that, up to a constant factor, agrees with the numerical results observed in \cite[Section~II.F]{grebenkov2024spectral}, see also Figure \ref{fig:LowerBoundsscaled} below. In addition, Remark \ref{rem:uppbound} will provide, later on, an upper bound that differs from the numerically observed asymptotic behaviour only by a logarithmic factor.

In order to better compare the behaviour of bounds for prolate spheroids for  different values of $a$ and to illustrate Theorem \ref{theo:noisopineq} in dimension three, we normalise the solid prolate spheroids by volume, by setting 
\[
\widetilde{\mathscr{p}}_a :=a^{-\frac{2}{3}}\mathscr{p}_a, \qquad  \widetilde{\mathcal{P}}_a:=a^{-\frac{2}{3}} \mathcal{P}_a, \qquad \text{and} \qquad  
\widetilde{\mathcal{P}}^\ext_a:=a^{-\frac{2}{3}} \mathcal{P}^\ext_a,
\]
so that $\left| \widetilde{\mathcal{P}}_a \right| = \frac{4\pi}{3}$ is constant for all $a\in(0,1)$. Then, bound \eqref{eq:weprolate} and the scaling $\sigma_1(\alpha\Omegaext) = \alpha^{-1} \sigma_1( \Omegaext)$ yield 
\begin{equation}\label{eq:lowboundprol}
\sigma_1\left(\widetilde{\mathcal{P}}_a^\ext\right)  \geq \beta\left(\widetilde{\mathscr{p}}_a\right)=a^\frac{2}{3} \beta\left(\mathscr{p}_a\right) = \frac{1-a^2}{-2a^\frac{1}{3} \log(a) } \xrightarrow[]{a \searrow 0} +\infty,
\end{equation}
and bound \eqref{eq:xiongsph} becomes
\begin{equation}\label{eq:xiongscaled}
\sigma_1\left(\widetilde{\mathcal{P}}_a^\ext\right)  \geq \beta_\mathrm{X}\left(\widetilde{\mathscr{p}}_a\right)=a^\frac{2}{3} \beta_\mathrm{X}\left(\mathscr{p}_a\right) 
= \begin{cases}
 \frac{3\sqrt{3}}{2}\frac{a^\frac{4}{3}}{(1+a^2)^\frac{3}{2}}\qquad&\text{if }0<a\le \frac{1}{\sqrt{2}},\\
a^\frac{2}{3}\qquad&\text{if }\frac{1}{\sqrt{2}}<a<1.
\end{cases}
\end{equation}
In Figure~\ref{fig:LowerBoundsscaled}, we compare both lower bounds and numerical computations of~$\sigma_1\left(\widetilde{\mathcal{P}}_a^\ext\right)$.

\begin{figure}[htb]
\centering
\includegraphics[width=0.95\textwidth]{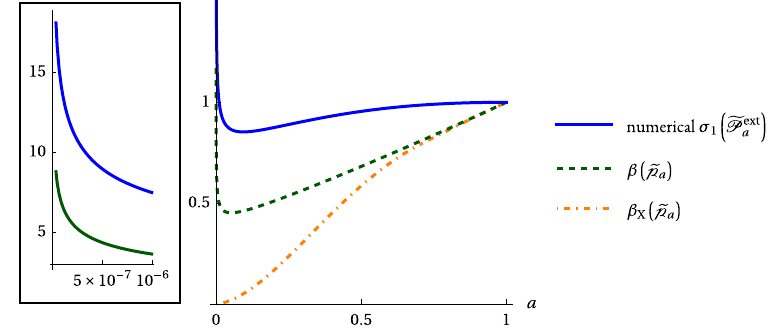}
\caption{Comparison between the numerical results, our bound  \eqref{eq:lowboundprol}, and Xiong's bound \eqref{eq:xiongscaled} for the exterior of rescaled prolate spheroids. The inset on the left zooms onto very small values of $a$.}%
\label{fig:LowerBoundsscaled}%
\end{figure}
\end{example}

If we consider an arbitrary, not necessarily convex or star-shaped, domain, there are no nontrivial lower bounds, as the following example shows.

\begin{definition}
Let $\Omega\subset\mathbb{R}^n$, $n\ge 2$,  be a bounded open set with a connected $\Omegaext$. We say that $\Omegaext$ has a \emph{passage of width} $\eps>0$ and \emph{base} $\Gamma\subset\mathbb{R}^{n-1}$ if, subject to a rigid motion  change of coordinates, there exists a bounded cylinder $\mathcal{C}_{\eps, \Gamma}:=\Gamma\times (-\epsilon,\epsilon)$ such that $\mathcal{C}_{\eps, \Gamma}\subset \Omegaext$ and
$\partial\mathcal{C}_{\eps, \Gamma}\cap\partial\Omega=\Gamma\times \{\pm\eps\}$, see Figure \ref{fig:thin}.
\end{definition}

\begin{figure}[htb]
\centering
\includegraphics{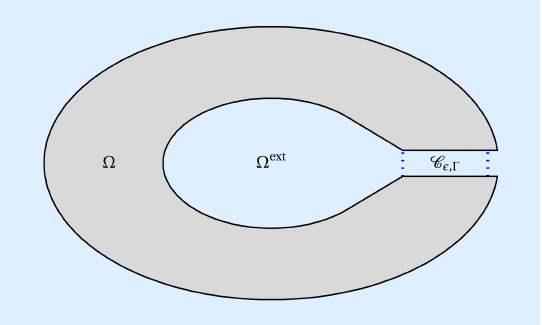}
\caption{An example of $\Omegaext$ with a passage.}%
\label{fig:thin}%
\end{figure}

\begin{proposition}\label{prop:thinpass}
Let $\Omega\subset\mathbb{R}^n$ be a bounded open set with Lipschitz boundary. If the exterior domain $\Omegaext$ has a passage of width $\eps>0$ and base $\Gamma\subset\mathbb{R}^{n-1}$, then its exterior Steklov eigenvalues satisfy
\[
\sigma_k\left(\Omegaext\right)<\Lambda^\Dir_k(\Gamma)\eps,\qquad  k\in\mathbb{N},
\]
where $\Lambda^\Dir_k(\Gamma)$ is the $k$th Dirichlet eigenvalue of $\Gamma$.
\end{proposition}
 
\begin{proof}
It follows from the variational characterisation~\eqref{eq:mu_k}, that
\[
\sigma_k\left(\Omegaext\right)\leq \sigma^\Dir_k\left(\mathcal{C}_{\eps, \Gamma}\right),
\]
where $\sigma^\Dir_k\left(\mathcal{C}_{\eps, \Gamma}\right)$ is the $k$th eigenvalue of the mixed Steklov--Dirichlet problem
\begin{equation}\label{eq:SDC}
\begin{cases}
\Delta u=0\qquad&\text{in }\mathcal{C}_{\eps, \Gamma},\\
u=0\qquad&\text{on }\partial\Gamma\times[-\eps,\eps],\\
\partial_\nu u = \sigma^\Dir u &\text{on } \Gamma\times\{\pm\eps\}.
\end{cases}
\end{equation}

It is easily seen by separation of variables that the set of eigenvalues of \eqref{eq:SDC} is given by
\[
\left\{\sqrt{\Lambda_k^\Dir(\Gamma)}\tanh\left(\eps \sqrt{\Lambda_k^\Dir(\Gamma)}\right), k\in\mathbb{N}\right\}
\cup
\left\{\sqrt{\Lambda_k^\Dir(\Gamma)}\coth\left(\eps \sqrt{\Lambda_k^\Dir(\Gamma)}\right), k\in\mathbb{N}\right\}.
\]
Ignoring the second subset, we conclude that $\sigma^\Dir_k\left(\mathcal{C}_{\eps, \Gamma}\right)\le \sqrt{\Lambda_k^\Dir(\Gamma)}\tanh\left(\eps \sqrt{\Lambda_k^\Dir(\Gamma)}\right) < \Lambda_k^\Dir(\Gamma) \eps$, and the result follows.
\end{proof}

If we now consider a family of $\eps$-dependent exterior domains, each containing a passage with the same fixed base and of width $\eps\searrow 0$, then the $k$th Steklov eigenvalue will tend to zero as well. For an illustration in dimension two, let $\Gamma:=\left(-\frac{\pi}{2},\frac{\pi}{2}\right)$, and let $\Omega_\eps\subset\mathbb{R}^2$ be a disjoint union of two identical squares of side $\pi$, $\Omega_\eps:=\Gamma\times (-\pi-\eps, -\eps)\sqcup \Gamma\times (\eps, \pi+\eps)$.  Then Proposition \ref{prop:thinpass} implies that 
\[
\sigma_k\left(\Omega_\eps^\ext\right) < \eps k^2,\qquad k\in\mathbb{N}.
\]

\begin{remark} The same result holds for the interior Steklov problem: if $\Omega\subset\mathbb{R}^n$ is a bounded open set with Lipschitz boundary, which contains a passage of width $\eps>0$ and base $\Gamma\subset\mathbb{R}^{n-1}$, then $\sigma_k\left(\Omega\right)<\Lambda^\Dir_k(\Gamma)\eps$ for all $k\in\mathbb{N}$.

\end{remark}

\subsection{Shape optimisation and spectral asymptotics}\label{sec:shapeopt}

\subsubsection{Proof of the isoperimetric inequality in dimension two}

\begin{proof}[Proof of Theorem \ref{theo:weinstock}]
Assume for simplicity that $ \left|\partial\Omega\right|=2\pi$.  As in \S\ref{subsec:conformal}, we  identify $\R^2$ with the complex plane $\C$ and consider the holomorphic map $\phi:\C\setminus\{0\}\to\C\setminus\{0\}$ defined by $\phi(z)=1/z$.  Consider, as before,  the bounded simply-connected domain  $\Omega^*$ such that $\phi\left(\Omega^*\setminus\{0\}\right)=\Omegaext$.  By the Riemann mapping theorem there exists a biholomorphic map  $\psi: B_1\to\Omega^*$. Therefore, the exterior problem \eqref{eq:PDEext} is equivalent to 
\begin{equation}\label{eq:SteklovDiskDensity}
\begin{cases}
\Delta w=0\qquad&\text{ in }B_1,\\
\partial_{\nu}w=\sigma \omega(z)w\qquad&\text{ on }\partial B_1,
\end{cases}
\end{equation}
where the function $\omega:\partial B_1\to\R$ is given by
\[
\omega(z)=\left|\left(\phi(\psi(z))\right)'\right|= \left| \left(\frac{1}{\psi(z)} \right)' \right|= \left|\frac{\psi'(z)}{\psi(z)^2} \right|,
\]
and satisfies $\int_{\partial B_1}\omega\,\dr S=|\partial\Omega|=2\pi$.
Following \cite[Theorem 4.12]{weinstock1954inequalities}, we plan to use the coordinate functions $x_1, x_2: B_1\to\R$ as trial functions for this problem. 
Indeed, using the Hersch trick, it is possible to choose the conformal map $\psi$ so that $\int_{\partial B_1}x_i \omega\,\dr S =0$, leading to
\[
\sigma_2\int_{\partial B_1}x_i^2\omega\,\dr S  \leq \int_{B_1}\left|\nabla x_i\right|^2\,\dr x =|B_1|=\pi.
\]
Summing over $i=1,2$, and using the fact that $x_1^2+x_2^2=1$ on $\partial B_1$, leads to 
\[
\sigma_2\int_{\partial B_1}\omega\,\dr S \leq 2\pi = \sigma_2\left|\partial\Omega\right|.
\] 

In case of equality,  arguing in the same way as in \cite[Section 2.3]{GP2010} and \cite[pp. 1038-1039]{FreitasLaugesen2020},  we get that  each $x_i$ is an eigenfunction of the weak version of \eqref{eq:SteklovDiskDensity} corresponding to eigenvalue one, and it follows that $\omega(z)=1$ for almost all $z\in\partial B_1$. Let $a:=\psi^{-1}(0)\in B_1$, and consider the conformal automorphism $\delta_a: B_1\to B_1$ of the unit disk defined by 
\[
\delta_{a}(z):=\frac{z-a}{1-\overline{a}z}.
\]
Because $\psi$ is a biholomorphic map, $\psi'(z)$ does not vanish for $z\in B_1$. Moreover, $a$ is the unique zero of $\psi$, and it is simple. It follows that $\frac{1}{\psi(z)}$ is a meromorphic function with a simple pole at $a$. Therefore, the function  $\upsilon(z):=\frac{\psi'(z)}{\psi(z)^2}\delta_a(z)^2$ has a removable singularity at $a$, and hence it can be extended to a holomorphic function on $B_1$. Note that $\upsilon$ does not vanish in $B_1$. Moreover, $|\upsilon(z)|=1$ for almost all $z\in\partial B_1$. Hence, the function $z\mapsto\log\left(|\upsilon(z)|\right)$ is harmonic on $B_1$ and it vanishes almost everywhere on $\partial B_1$. Since $\Omega$ is Lipschitz, the conformal map $\psi$ is regular enough for $\log |\upsilon|$ to admit a representation as a Poisson integral on $\partial B_1$, hence $\log |\upsilon|$ vanishes everywhere on $B_1$,  see \cite{FreitasLaugesen2020}. This implies the existence of a constant $c$ with $|c|=1$, such that $\upsilon(z)=c$ for all $z\in B_1$.  As a consequence, for any $z\in B_1\setminus\{a\}$ we have,
\begin{equation}\label{eq:caK}
\left( \frac{1}{\psi(z)}\right)'=c\left(\frac{1-\overline{a}z}{z-a}\right)^{2}.
\end{equation}
Let us compute the residues at $a$ of both sides of \eqref{eq:caK}. A direct computation yields
\[
c\left(\frac{1-\overline{a}z}{z-a}\right)^{2}=c\left(\frac{(1-|a|^2)^2}{(z-a)^2}-2\overline{a}\frac{1-|a|^2}{z-a}+\overline{a}^2\right).
\]
It follows that the residue of the right-hand side of  \eqref{eq:caK}  at $a$ is equal to $-2c\overline{a}(1-|a|^2)$. On the other hand, the left-hand side of  \eqref{eq:caK}  is the derivative of a meromorphic function with a pole at $a$. Using the Laurent series representation of $1/\psi$  we find that the residue  of the left-hand side at $a$ is equal to zero. Combining these two calculations implies that $a=0$. Hence,
\[
\left( \frac{1}{\psi(z)} \right)'=\frac{c}{z^2}.
\]
Integrating,  we obtain
\[
\frac{1}{\psi(z)}=\frac{-c}{z}+C,
\]
where $C\in\C$ is a constant of integration.
Therefore,
\[
\psi(z)=\frac{z}{Cz-c}=\frac{1}{C-\frac{c}{z}},
\]
and $|C|<1$ since $\psi$ is holomorphic. It is now easy to see that $\Omega^{*}=\psi(B_1)$ is a disk, and hence  $\Omega$ is a disk as well. 
\end{proof}

\subsubsection{Higher dimensions}

In higher dimensions the exterior of a ball is known to be a local maximiser for the first eigenvalue among domains with given measure \cite{bundrock2023robin}, cf.\ the eigenvalue curves in Figure~\ref{fig:LowerBounds} in the vicinity of $a=1$. However, even within the class of convex domains, no global isoperimetric bound holds. In this regard, the exterior problem differs significantly from the interior problem. 

Recall Example~\ref{examp:ellipsoid}:   the normalised  solid spheroids $\widetilde{\mathcal{P}}_a$ have constant volume, while the lower bound \eqref{eq:lowboundprol} shows that $\lim_{a \searrow 0} \sigma_1(\widetilde{\mathcal{P}}_a^\ext) =+ \infty$. Hence, the ball does not maximise the first Steklov exterior eigenvalue among domains of fixed volume. Moreover, if the ball were the maximiser among domains with fixed surface area, then by the classical isoperimetric inequality it would also maximise among domains with fixed volume, which contradicts the above example.

An analogous behaviour can be observed for higher-dimensional spheroids. The following example provides an immediate proof of Theorem \ref{theo:noisopineq} for all $n\ge 3$.

\begin{example}
Let $n\ge 3$. For $a \in (0,1)$, and $k \in \mathbb{N}$, $2 \leq k \leq n-1$, we consider a spheroid
\[
\mathscr{e}_a^{(k)} := \left\{ x \in \mathbb{R}^n: \sum_{j=1}^k \frac{x_j^2}{a^2} + \sum_{j=k+1}^n x_j^2 = 1 \right\}.
\]
We write $x = \begin{pmatrix}x^{(1)}\\x^{(2)}\end{pmatrix}$, with $x^{(1)} \in \mathbb{R}^k$ and $x^{(2)} \in \mathbb{R}^{n-k}$, and define  $D(x):=1-(1-a^2)\left|x^{(2)}\right|^2$. Then, the principal  curvatures of $\mathscr{e}_a^{(k)}$ are given by
\[
\begin{split}
\kappa_1(x) &= \ldots = \kappa_{k-1}(x) = \frac{1}{a\sqrt{D(x)}},\\
\kappa_k(x) &= \ldots = \kappa_{n-2}(x) = \frac{a}{\sqrt{D(x)}}, \quad \kappa_{n-1}(x) = \frac{a}{D(x)^\frac{3}{2}}.
\end{split}
\]
Each principal curvature becomes minimal at $\left|x^{(2)}\right|=0$, so the logarithmic mean does also become minimal when $\left|x^{(2)}\right|=0$. Note that we have to understand $L(\kappa_1, \ldots,\kappa_{n-1})$ as a limit. For $k=2$, this is
\begin{multline*}
\min_{x \in  \mathscr{e}_a^{(2)}} L(\kappa_1(x), \ldots, \kappa_{n-1}(x)) \\= 
\frac{1}{(n-2)\left(-\log(a) \frac{a}{(1-a^2)^{n-2}}-\frac{a}{(n-3)!}\left.\frac{\mathrm{d}^{n-3}}{\mathrm{d}t^{n-3}}\left( \frac{t^{n-3} \log(t)}{1-at} \right)\right|_{t=a}\right)}. 
\end{multline*}
For $1 \leq i \leq m$,
\[
\frac{\mathrm{d}^{i}}{\mathrm{d}t^{i}} \log(t) = (-1)^{i-1}\frac{(i-1)!}{t^i}, \quad \frac{\mathrm{d}^{i}}{\mathrm{d}t^{i}} t^m = \frac{m!}{(m-i)!}t^{m-i}, \quad \frac{\mathrm{d}^{i}}{\mathrm{d}t^{i}} \frac{1}{1-at} = \frac{a^i i!}{(1-at)^{i+1}},
\]
and applying the Leibniz rule, we obtain
\[
\frac{\mathrm{d}^{m}}{\mathrm{d} t^{m}}
\left(\frac{t^{m}\log t}{1-at}\right) = \frac{m!}{(1-at)^{m+1}} \log(t) + O(t)\qquad\text{as }t\searrow 0.
\]
Therefore,
\[
\min_{x \in  \mathscr{e}_a^{(2)}} L(\kappa_1(x), \ldots, \kappa_{n-1}(x)) 
= \frac{(1-a^2)^{n-2}}{(n-2)\left(-2a \log(a) + O(a)\right)}\qquad\text{as }a\searrow 0.
\]
So,  bound \eqref{eq:logmean} implies that for the normalised solid spheroid 
\[
\widetilde{\mathcal{E}}_a^{(2)} := a^{-\frac{2}{n}} \left\{ x \in \mathbb{R}^n:  \frac{x_1^2+x_2^2}{a^2} +x_3^{2}+\ldots +x_n^2 \leq 1 \right\},
\] 
with volume $\left|\widetilde{\mathcal{E}}_a^{(2)}\right|=\omega_n$, we have
\[
\sigma_1\left( \left( \widetilde{\mathcal{E}}_a^{(2)}\right)^\ext \right)  \geq  \frac{(1-a^2)^{n-2}}{-2a^{\frac{n-2}{n}}  \log(a) +O\left( a^{\frac{n-2}{n}}\right)} \to +\infty \qquad \text{as } a \searrow 0.
\]

For $k >2$, the argument carries over in the same way, leading again to $\sigma_1\left( \left(\tilde{\mathcal{E}}_a^{(k)}\right)^\ext \right) \to \infty$ as $a \searrow 0$. In fact, the divergence is faster, which can be understood informally from the fact that the principal curvature  $\frac{1}{a\sqrt{D(x)}}$ is weighted more heavily in the logarithmic mean as $k$ increases.
\end{example}

Another interesting non-convex example demonstrating the absence of an isoperimetric inequality is presented below.

\begin{example}
Let $M$ be a smooth compact Riemannian manifold of dimension $m\geq 3$, without boundary. Let $N\subset M$ be a closed submanifold of positive codimension and for $\eps>0$, consider the tubular neighbourhood $N_\eps=\{x\in M: \operatorname{dist}(x,N)<\eps\}$.
Brisson  \cite{BrissonJ2022} studied the Steklov spectrum of $N_\eps^\ext:=M\setminus N_\eps$. 
She proved that
\[
\lim_{\eps\to 0}\sigma_1\left(N_\eps^\ext\right)\left|\partial N_\eps^\ext \right|^{\frac{1}{m-1}}=\infty.
\]
Because the approach in \cite{BrissonJ2022} is based on quasi-isometries and mixed Steklov--Neumann eigenvalue problems, it can be adapted verbatim in the setting of this paper to obtain the following result, via the mixed approach of Theorem \ref{thm:resolvconv}.

\begin{proposition}
Let $m\geq 3$. Consider a compact smooth submanifold $N \subset \R^m$ of dimension $ n\leq m-2$. Then, the exterior Steklov eigenvalues of the tubular neighbourhood $N_\eps$ satisfy
\[
\lim_{\eps\to 0}\sigma_1\left(N_\eps^\ext\right)\left|\partial N_\eps\right|^{\frac{1}{m-1}}=\infty.
\]
\end{proposition}
This applies, in particular, to  neighbourhoods of closed curves in $\R^3$.
\end{example}

\begin{remark}\label{rem:uppbound} 
In \cite[Theorem~4]{xiong2023sharp}, several upper bounds for $\sigma_1$ are established. 
Another such bound was obtained in \cite[Theorem~3]{bundrock2024optimizing}, following the approach of \cite[Section~2]{giorgi2005monotonicity}, which shows the following.
Suppose $\Omega \subset \mathbb{R}^n$, $n \geq 3$, is a bounded Lipschitz domain, and $B_\rho \subset \Omega$ for some $\rho>0$. Then $\sigma_1\left(\Omegaext\right) \leq \sigma_1\left(B_\rho^\ext\right)$. As an immediate consequence of Theorem~\ref{theo:weinstock}, the same inequality also holds for the first nontrivial Steklov eigenvalue in dimension two: $\sigma_2\left(\Omegaext\right) \leq \sigma_2(B_\rho^\ext) = \frac{1}{\rho}$.

The same idea as in \cite[Theorem~3]{bundrock2024optimizing} can be used to prove, for any $\Lambda >0$, that
\[
B_\rho \subset \Omega \quad \implies \quad 
\mu_1\left(\Omegaext,\Lambda\right) \leq \mu_1\left(B_\rho^\mathrm{ext},\Lambda\right).
\]
Let $u(x) = f(|x|)$ denote the eigenfunction corresponding to $\mu_1\left(B_\rho^\ext,\Lambda\right)=\mu_{(0)}\left(B_\rho^\ext,\Lambda\right)$, see \eqref{eq:EFHelmHoltzBall}. We define the function
\[
\widehat{\mu}: \partial \Omega \to \R, \quad y \mapsto \frac{\mydotp{ \nabla f(|y|), \nu_y}}{f(|y|)}.
\]
Note that $f'(\rho) < 0$, and that $\frac{f'(r)}{f(r)} = \frac{-\Lambda K_{\frac{n}{2}}(\Lambda r)}{K_{\frac{n-2}{2}}(\Lambda r)}$ is monotonically increasing in $r$, see  \cite[proof of Proposition 3]{krejvcivrik2020optimisation}. Thus,
\[
\widehat{\mu}(y) \leq \frac{-f'(|y|)}{f(|y|)} \leq \frac{-f'(\rho)}{f(\rho)} = \mu_1(B_\rho^\ext,\Lambda)\qquad\text{for all }y\in \partial \Omega.
\]
Furthermore, $u$ solves 
\[
\begin{cases}
(\Lambda^2 -\Delta) u = 0 \qquad&\text{ in } \Omega^{\text{ext}} \subset B_\rho^{\ext}, \\
 \partial_\nu u  = \widehat{\mu} u \qquad&\text{ on } \partial \Omega.
\end{cases}
\]
Thus, integration by parts  yields 
\[
\begin{split}
\Lambda^2 \int_{\Omegaext} u^2 \, \mathrm{d}x &= -\int_{\Omegaext} |\nabla u|^2 \, \mathrm{d}x + \int_{\partial \Omega} \widehat{\mu} u^2 \, \mathrm{d}S\\ 
&\leq -\int_{\Omegaext} |\nabla u|^2 \, \mathrm{d}x + \mu_1(B_\rho^\ext,\Lambda) \int_{\partial \Omega}  u^2 \, \mathrm{d}S.
\end{split}
\]
Therefore, 
\[
\mu_1(B_\rho^\ext,\Lambda) \geq\frac{\Lambda^2 \int_{\Omegaext} u^2 \, \mathrm{d}x+\int_{\Omegaext} |\nabla u|^2 \, \mathrm{d}x}{\int_{\partial \Omega}  u^2 \, \mathrm{d}S} \geq \mu_1\left(\Omegaext,\Lambda\right),
\]
which completes the proof.
\end{remark}

\subsubsection{Proofs of  eigenvalue asymptotics}\label{sec:asymptproofs}

\begin{proof}[Proof of Proposition~\ref{thm:Asymptotics}]
Let $\Omega\subset\R^n$ be a bounded open set  with $C^{2,\alpha}$ boundary for some 
$\alpha>0$. For $R>R_0(\Omega)$, it easily follows from Lemma \ref{lemma:monotonicitymixed} (see also  Lemma \ref{lemma1})  that 
\begin{equation}\label{eq:extbrack}
\sigma_k^\Neu\left(\Omegaext_R\right)\leq\sigma_k\left(\Omegaext\right)\leq\sigma_k^\Dir\left(\Omegaext_R\right)
\end{equation}
for every $k\ge 1$.
Inspecting the proof of \cite[Theorem 1.11]{GKLP2022},  we 
conclude that it can be applied verbatim to
the mixed boundary problems on $\Omegaext_R$ with Steklov condition on $\partial \Omega$ and either Dirichlet or Neumann condition on $\partial B_R$. 
Indeed, the  argument is local near the Steklov part of the boundary: the vector field used in the H\"{o}rmander--Pohozhaev identity \cite[Theorem 1.5]{GKLP2022} can be chosen to be supported in a small neighbourhood of $\partial\Omega$, away from $\partial B_R$. Therefore the Dirichlet or Neumann condition imposed on $\partial B_R$ does not affect the comparison with 
$\sqrt{-\Delta_{\partial\Omega}}$, and only $\partial\Omega$ contributes to the leading Weyl term.
Therefore, 
\[
\#\left\{k\in\N\,:\,\sigma_k^{\Dir,\Neu}\left(\Omegaext_R\right)<\sigma\right\}=\frac{\omega_{n-1}}{(2\pi)^{n-1}}  |\partial\Omega| \sigma^{n-1} + O\left(\sigma^{n-2}\right)\qquad\text{as } \sigma\to+\infty.
\]
Combining this with \eqref{eq:extbrack} yields the result.
\end{proof}

\begin{proof}[Proof of Proposition~\ref{prop:Weyl2n}]
Let us use the conformal approach and consider the spectrum of the weighted interior problem \eqref{eq:problemtransform}. Note that $\Omega^*$ is a Lipschitz domain, since 
an inversion with the centre away from the boundary preserves the Lipschitz regularity of the boundary. Therefore,  it follows from \cite[Theorem 1.8]{KLP2023} that the eigenvalue counting function for  the problem \eqref{eq:problemtransform} satisfies the  asymptotics \eqref{eq:Weyl2}. 
\end{proof}

\begin{ack}
This project originated from questions that have been raised in conversations with Bruno Colbois. The authors are very grateful to him for numerous stimulating exchanges related to this work. We would also like to thank Dorin Bucur,  Yakar Kannai, Yuri Latushkin, and Nilima Nigam for useful discussions,
as well as the anonymous referee for helpful suggestions.
\end{ack}

\begin{funding}
Research of AG  and IP was partially supported by NSERC and FRQNT.  Research of ML was partially supported by the EPSRC. DSG acknowledges the support of the Simons Foundation for his sabbatical stay in 2024 at the CRM in  Montr\'eal, Canada.
\end{funding}

\section*{Conflicts of interest statement}
There are no conflicts of interest.

\section*{Data availability statement}\label{sec:data}
The accompanying \texttt{Mathematica} and \texttt{FreeFEM} scripts, their printouts, and additional images  are available for download at \url{https://www.michaellevitin.net/exteriorsteklov.html} or at \url{https://github.com/michaellevitin/exteriorsteklov}. 

\section*{AI usage disclosure}
The authors used different LLMs, in particular ChatGPT 5.6 and Claude Opus 5, for galley proof checks and minor corrections in the final version of the manuscript.

\end{document}